\documentclass[11 pt,a4paper,twoside]{article}
\usepackage{amsmath,amssymb,amsthm,epsfig,graphics,amsfonts}
\usepackage[latin1]{inputenc}
\usepackage[francais,english]{babel}
\usepackage[all]{xypic}
\usepackage{bm}
\usepackage{pstricks,pst-plot}

\newcommand{\DD}{\mathbb{D}}

\newcommand{\NN}{\mathbb{N}}

\newcommand{\QQ}{\mathbb{Q}}
\newcommand{\RR}{\mathbb{R}}
\renewcommand{\SS}{\mathbb{S}}
\newcommand{\TT}{\mathbb{T}}

\newcommand{\ZZ}{\mathbb{Z}}

\def\cA{{\cal A}}  \def\cG{{\cal G}} \def\cM{{\cal M}} 
\def\cB{{\cal B}}  \def\cH{{\cal H}}  
    
\def\cD{{\cal D}}   \def\cP{{\cal P}} 
\def\cE{{\cal E}}    
\def\cF{{\cal F}}

\def\bPsi{\bm{\Psi}}
\def\bM{\bm{M}}
\def\bcE{\bm{\cE}}
\def\bcA{\bm{\cA}}
\def\bcG{\bm{\cG}}
\def\bE{\bm{E}}

\def\ra{\rightarrow}
\renewcommand{\phi}{\varphi}
\renewcommand{\epsilon}{\varepsilon}
\newcommand{\La}{\Lambda}

\renewcommand{\t}{\widetilde}
\def\wt{\widetilde}
\def\wh{\widehat}

\def\adhe{\mathrm{Cl}}

\def\inte{\mathrm{Int}}

\def\diam{\mathrm{Diam}}

\def\supp{\mathrm{Supp}}

\newcommand{\id}{\operatorname{Id}}

\newtheorem{theo}{Theorem}[section]
\newtheorem*{theo*}{Theorem}
\newtheorem{prop}[theo]{Proposition}
\newtheorem{coro}[theo]{Corollary}
\newtheorem*{conj*}{Conjecture}
\newtheorem{lemma}[theo]{Lemma}

\newtheorem{fact}[theo]{Fact}

\newtheorem{adde}[theo]{Addendum}

\theoremstyle{definition}
\newtheorem{defi}[theo]{Definition}

\newtheorem{rema}[theo]{Remark}

\newtheorem*{rema*}{Remark}
\newtheorem*{remas*}{Remarks}

\newtheorem{ques}{Question}

\title{Construction of curious minimal uniquely ergodic homeomorphisms on manifolds: the Denjoy-Rees technique}
\author{F. B\'eguin\footnote{Laboratoire de math\'ematiques, Univ. Paris Sud, 91405 Orsay Cedex, France.}, S. Crovisier\footnote{CNRS - Laboratoire Analyse, G\'eom\'etrie et Applications, Univ. Paris 13, 93430 Villetaneuse, France.}
and F. Le Roux\footnote{Laboratoire de math\'ematiques, Univ. Paris Sud,  91405 Orsay Cedex, France.}}

\begin{document}
\renewcommand{\thepage}{\roman{page}}

\maketitle

\sloppy

\begin{abstract}
In~\cite{Ree}, Mary Rees has constructed a minimal homeomorphism of the  $2$-torus with positive topological entropy. This homeomorphism $f$ is obtained by enriching the dynamics of an irrational rotation $R$. We improve Rees construction, allowing to start with any homeomorphism $R$ instead of an irrational rotation and to control precisely the measurable dynamics of $f$. This yields in particular the following result: \textit{Any compact manifold of dimension $d\geq 2$ which carries a minimal uniquely  ergodic homeomorphism also carries a  minimal uniquely ergodic homeomorphism with positive topological entropy.}

More generally, given some homeomorphism $R$ of a (compact) manifold and some
homeomorphism $h_{C}$ of a Cantor set, we construct a homeomorphism $f$
which ``looks like'' $R$ from the topological viewpoint and ``looks like'' $R\times h_{C}$
from the measurable viewpoint. This construction can be seen as a partial answer to the following realisability question: which measurable dynamical systems are represented by homeomorphisms on manifolds ?

\end{abstract}

\newpage
\tableofcontents

\newpage
\setcounter{page}{1}
\setcounter{section}{0}
\renewcommand{\thepage}{\arabic{page}}

\section{Introduction}
\label{s.intro}

\subsection{Denjoy-Rees technique}
\label{ss.Denjoy-Rees-technique}

Twenty-five years ago, M. Rees has constructed a homeomorphism of the torus $\TT^d$ ($d \geq 2$) which is minimal and has positive topological entropy (see \cite{Ree}). The existence of such an example is surprising for several reasons: 
\begin{itemize}
\item[--] classical examples of minimal homeomorphisms (irrational rotations, time $t$ maps of horocyclic flows, etc.) are also typical examples of zero entropy maps.
\item[--] a classical way for proving that a map $f$ has positive topological entropy is to show that the number of periodic orbits of period $\leq n$ for $f$ grows exponentially fast when $n\to\infty$. So, in many situations, ``positive topological entropy" is synonymous of ``many periodic orbits". But a minimal homeomorphism do not have any periodic orbit. 
\item[--] a beautiful theorem of A. Katok states that, if $f$ is a $C^{1+\alpha}$ diffeomorphism of a compact surface $S$ with positive topological entropy, then there exists an $f$-invariant compact set $\Lambda\subset S$ such that some power of $f_{\mid\Lambda}$ is conjugate to a full shift (see \cite[corollary 4.3]{Kat}). In particular, a  $C^{1+\alpha}$ diffeomorphism of a compact surface with positive topological entropy cannot be minimal.
\end{itemize}

Beyond the mere existence of minimal homeomorphisms of $\TT^d$ with positive topological entropy, the technique used by Rees to construct such a homeomorphism is very interesting. This technique can be seen as a very sophisticated generalisation of the one used by A. Denjoy to construct his famous counter-example (a periodic orbit free homeomorphism of $\SS^1$ which is not conjugate to a rotation, \cite{Den}). Indeed, the basic idea of Rees  construction is to start with an irrational rotation of $\TT^d$, and to ``blow-up'' some orbits, just as in Denjoy counter-example. Of course, the construction of Rees is much more complicated and delicate than the one of Denjoy; for example, to get a homeomorphism with positive topological entropy, one has to blow up a set of orbits of positive Lebesgue measure.

The aim of the present paper is to describe a general setting for what we call the \emph{Denjoy-Rees technique}. This general setting includes as particular cases the construction of various ``Denjoy counter-examples" in any dimension, and Rees construction of a minimal homeomorphism of $\TT^d$ with positive topological entropy. Moreover, we will develop a new technique which allows to control that the homeomorphisms we obtain ``do not contain too much dynamics". This yields new results such as the existence of minimal \emph{uniquely ergodic} homeomorphisms with positive topological entropy, or the possibility to realise many measurable dynamical systems as minimal homeomorphisms on manifolds.

\subsection{Strictly ergodic homeomorphisms with positive topological entropy}
\label{ss.minimal-positive-entropy}

A homeomorphism is said to be \emph{strictly ergodic} if it is minimal and uniquely ergodic. As an application of Denjoy-Rees technique, we will prove the following theorem.

\begin{theo}
\label{t.entropy}
Any compact manifold of dimension $d\geq 2$ which carries a strictly ergodic homeomorphism also carries a  strictly ergodic homeomorphism with positive topological entropy.  
\end{theo}

A. Fathi and M. Herman have proved that every compact manifold of dimension $d\geq 2$ admitting a locally free action of the circle\footnote{An action of the circle is said to be \emph{locally free} if no orbit  of this action is reduced to a point.} carries a strictly ergodic homeomorphism (see \cite{FatHer}). Putting theorem~\ref{t.entropy} together with Fathi-Herman result yields many examples. 
In particular, the torus $\TT^d$ for $d\geq 2$, the sphere $\SS^{2n+1}$ for $n\geq 1$, any Seifert manifold, any manifold obtained as a quotient of a compact connected Lie group, \emph{etc.}, carry strictly ergodic homeomorphisms with positive topological entropy (see the discussion in \cite{FatHer}).  

Theorem~\ref{t.entropy} (as well as Rees example and Katok theorem cited above) can be seen as a piece of answer to a general question of Herman asking ``whether, for diffeomorphisms,  positive topological entropy is compatible with minimality, or strict ergodicity" (see~\cite[page 141]{Kat}). 
Katok answered negatively to Herman question in the case of $C^{1+\alpha}$ diffeomorphisms of surfaces.  Then, Herman himself constructed an analytic minimal diffeomorphism with positive topological entropy on a $4$-manifold (see \cite{Her1}), and Rees constructed a minimal homeomorphism with positive topological entropy on $\TT^d$. But neither Herman, nor Rees managed to make their examples strictly ergodic (see the introductions of~\cite{Ree} and~\cite{Kat}). Theorem~\ref{t.entropy} shows that positive topological entropy is compatible with strict ergodicity for homeomorphisms (in any dimension). 
To complete the answer to Herman question, it essentially remains to determine what is the best possible regularity  for a minimal (resp. strictly ergodic) homeomorphism on $\TT^2$ with positive entropy (H\"older? $C^1$?).
We do not have any idea of the best regularity one can obtain for a homeomorphism constructed via Denjoy-Rees technique.

Note that, by Oxtoby-Ulam theorem~\cite{OxtUla}, one can assume that the unique measure preserved by the homeomorphism provided by theorem~\ref{t.entropy} is a Lebesgue measure.

\subsection{Realising measurable dynamical systems as homeomorphisms on manifolds}
\label{ss.isomorphic}

The main difference between Rees result and our theorem~\ref{t.entropy} is the fact that the homeomorphisms we construct are uniquely ergodic. More generally, we develop a technique which allows us to control the number of invariant measures of the homeomorphisms obtained by constructions \emph{\`a la Denjoy-Rees}. What is the point of controlling the invariant measures? In short:
\begin{itemize}
\item[--] the Denjoy-Rees technique by itself is a way for constructing examples of ``curious" minimal homeomorphisms, 
\item[--] the Denjoy-Rees technique combined with the possibility of controlling the invariant measures is 
not only a way for constructing examples, but also a way for realising measurable dynamical systems as homeomorphisms on manifolds.
\end{itemize}

Let us explain this. In her paper, Rees constructed a homeomorphism $f$ on $\TT^d$ which is minimal and possesses an invariant probability measure $\mu$ such that $f$ has a rich dynamics from the point of view of the measure $\mu$: in particular,  the metric entropy $h_\mu(f)$ is positive. By the variational principle, this implies that the topological  entropy $h_{\text{top}}(f)$ is also positive. Nevertheless, $f$ might possess some dynamics that is not detected by the measure $\mu$ (for example, $h_{\text{top}}(f)$ might be much bigger than $h_\mu(f)$).
So, roughly speaking, Rees constructed a minimal homeomorphism which has a rich dynamics, but without being able to control how rich this dynamics is. Now, if we can control what are the invariant measures of $f$, then we know exactly what $f$ looks like from the measurable point of view.

To make this precise, we need some  definitions. For us, a \emph{measurable dynamical system}
 $(X,\cA,S)$ is a bijective bi-measurable map $S$ on a set $X$ with a $\sigma$-algebra $\cA$.
A $S$-invariant set $X_{0} \subset X$ is \emph{full} if it has full measure for any $S$-invariant probability measure on $(X,\cA)$. Two measurable systems $(X,\cA,S)$  and $(Y,\cB,T)$  are  \emph{isomorphic} if there exist a $S$-invariant full set $X_{0} \subset X$, a $T$-invariant full set  $Y_{0} \subset Y$ and a bijective bi-measurable map $\Theta : X_{0} \ra Y_{0}$ such that $\Theta \circ S = T \circ \Theta$. With these definitions, an interesting general question is:

\begin{ques}
\label{q.realisability}
Given any measurable dynamical system $(X,\cA,S)$, does there exist a homeomorphism on a manifold which is isomophic to $(X,\cA,S)$?
\end{ques}

In this direction, using Denjoy-Rees technique together with our technique for controlling the invariant measures, we will prove the following ``realisation theorem" (which implies theorem~\ref{t.entropy}, see below).

\begin{theo}
\label{t.realisability}
Let $R$ be a uniquely ergodic aperiodic homeomorphism of a compact manifold $\cM$ of dimension $d\geq 2$. Let $h_C$ be a homeomorphism on some Cantor space $C$. Then there exists a homeomorphism $f:\cM\to\cM$ isomorphic to $R\times h_C:\cM\times C\to\cM\times C$. 

Furthermore, the homeomorphism $f$ is a topological extension of $R$: there exists a continuous map $\Phi:\cM\to\cM$ such that $\Phi\circ f=R\circ\Phi$. If $R$ is minimal (resp. transitive), then $f$ can be chosen minimal (resp. transitive).  
\end{theo}

\begin{rema*}
In~\cite{Ree}, Rees considered the case where $R$ is an irrational rotation of $\TT^d$ and $h_C$ is a full shift. She constructed a minimal homeomorphism $g$ which had a subsystem isomorphic to  $R\times h_C$, but was not isomorphic to $R\times h_C$. 
\end{rema*}

Let us comment our definitions and theorem~\ref{t.realisability}.  

The problem of realising measurable dynamical systems as topological dynamical systems admits many alternative versions. One possibility is to prescribe the invariant measure, that is, to deal with \emph{measured} systems instead of \emph{measurable} systems. In this context, the realisability problem consists in finding homeomorphisms $f$ on a manifold $\cM$ and an $f$-invariant measure $\nu$ such that $(\cM,f,\nu)$ is metrically conjugate to a given dynamical system $(X,T,\mu)$. In this direction, D. Lind and J.-P. Thouvenot proved that every finite entropy measured dynamical system is metrically conjugate  to  some shift map on a finite alphabet, and thus also to a Lebesgue measure-preserving homeomorphism of the two-torus  (see \cite{LinTho}).

Then, one can consider the same problem but with the additional requirement that the realising homeomorphism is uniquely ergodic. In this direction (but not on manifolds), one has the celebrated Jewett-Krieger theorem: any ergodic system is metrically conjugate  to a uniquely ergodic homeomorphism on a Cantor space (see e.g. \cite{DenGriSig}).  In a forthcoming paper (see \cite{BegCroLeR}), we will adress this problem for homeomorphisms of manifolds: we use theorem~\ref{t.realisability} (more precisely, the generalisation of theorem~\ref{t.realisability} stated in subsection~\ref{ss.more-statements}) to prove that any measured system whose discrete spectrum contains an irrational number (\textit{i. e.} which is a measurable extension of an irrational rotation of the circle) is metrically conjugate to a minimal uniquely ergodic homeomorphism of the two-torus $\TT^2$.

If one seeks realisations of measurable (or measured) dynamical systems by smooth maps then the Denjoy-Rees technique seems to be useless. The main technique for constructing ``curious" minimal $C^\infty$ diffeomorphisms  on manifolds was introduced by D. Anosov and Katok in \cite{AnoKat}, and developed by many authors, including Fathi, Herman, B. Fayad, Katok, A. Windsor, etc. Note that this technique is not really adapted for constructing diffeomorphisms that are isomorphic to a given measurable dynamical system,  but it allows to construct minimal diffeomorphisms that have exactly $n$ ergodic invariant probability measures, or  with such or such spectral property.  See~\cite{FayKat} for a survey. 

An interesting feature of Anosov-Katok technique is that it allows to construct examples of \emph{irrational pseudo-rotations} of the torus $\TT^2$, i.e. homeomorphisms whose rotation set is reduced to a single irrational point of $\RR^2/\ZZ^2$ (this irrational point is generally liouvillian, which is the price to pay in order to get the smoothness of the pseudo-rotation). Observe that, if one applies theorem~\ref{t.realisability} with the homeomorphism $R$ being an irrational rotation of the torus $\TT^2$ (and $h_C$ being any homeomorphism  of a Cantor set), then it is easy to see that the resulting homeomorphism $f$ is also an irrational pseudo-rotation. Varying the homeomorphism $h_C$, one gets lots of examples of  ``exotic" irrational pseudo-rotations on the torus $\TT^2$. We point out that with this method the angle of the rotation $R$ may be any irrational point of $\RR^2/\ZZ^2$ (but the pseudo-rotation we obtained are just homeomorphisms). For a general discussion on pseudo-rotations, see the introduction of~\cite{BegCroLeRPat}.

\bigskip

Let us try to give a brief idea of what the homeomorphism $f$ provided by theorem~\ref{t.realisability} looks like (see also section~\ref{ss.outline}). On the one hand, from the topological point of view, $f$ looks very much like the initial homeomorphism $R$ (which typically can be a very simple homeomorphism, like an irrational rotation of the torus $\TT^2$). Indeed, the continuous map $\Phi:\cM\to\cM$ realising the topological semi-conjugacy between $f$ and $R$ is an ``almost  conjugacy": there exists an $f$-invariant $G_{\delta}$-dense set $X$ on which $\Phi$ is one-to-one. This implies that $f$ is minimal. On the other hand, from the measurable point of view, $f$ is isomorphic to the product $R\times h_C$ (which might exhibit a very rich dynamics since $h_C$ is an arbitrary homeomorphism on a Cantor set). As often, the paradox comes from the fact that the set $X$ is big from the topological viewpoint (it is a $G_\delta$ dense set), but small from the measurable viewpoint (it has zero measure for every $f$-invariant measure).

\bigskip

We end this section by explaining how theorem~\ref{t.realisability} implies theorem~\ref{t.entropy}. 
Let $\cM$ be a manifold of dimension $d\geq 2$, and assume that there exists a strictly ergodic homeomorphism $R$ on $\cM$. We have to construct a strictly ergodic homeomorphism with positive entropy on $\cM$. For this purpose, we may assume that the topological entropy of $R$ is equal to zero, otherwise there is nothing to do. 

Let $\sigma$ be the shift map on $\{0,1\}^\ZZ$, and $\mu$ be the usual Bernoulli measure on $\{0,1\}^\ZZ$. By Jewett-Krieger  theorem (see e.g.~\cite{DenGriSig}), there exists a uniquely ergodic homeomorphism $h_C$ of a Cantor set $C$ which is metrically conjugate to $(\{0,1\}^\ZZ,\sigma,\mu)$. Since the shift map is a $K$-system, and since $R$ is uniquely ergodic and has zero topological entropy, the product map $R \times h_{C}$ is also uniquely ergodic  (see \cite[proposition 4.6.(1)]{Tho2}). Denote by $\nu$ the unique invariant measure of $R\times h_C$. Then the metric entropy $h_\nu(R\times h_{C})$ is equal to $h_\mu(\sigma)=\log 2$. Now theorem~\ref{t.realisability} provides us with a minimal homeomorphism $f$ on $\cM$, which is isomorphic to  $R\times h_{C}$. Denote by $\Theta$ the map realising the isomorphism between $R\times h_C$  and $f$. Since $R\times h_{C}$ is uniquely ergodic, $f$ is also uniquely ergodic: $\Theta_*\nu$ is the unique $f$-invariant measure). Moreover, the metric entropy $h_{\Theta_*\nu}(f)$ is equal to $h_\nu(R\times h_C)$, which is positive.  The variational principle then implies that the topological entropy of $f$ is positive. Hence, $f$ is a strictly ergodic homeomorphism with positive topological entropy.

\subsection{A more general statement}
\label{ss.more-statements}

Theorem~\ref{t.realisability} is only a particular case of a more general statement: theorem~\ref{t.main} below will allow to consider a homeomorphism $R$ that is not uniquely ergodic, and to replace the product map $R\times h_C$ by any map that fibres over $R$.
 
To be more precise, let $R$ be a  homeomorphism of a manifold $\cM$ and $A$ be any measurable subset of $\cM$. Let $C$ be a Cantor set. We consider a bijective bi-measurable map which is fibered over $R$:
\begin{eqnarray*}
h:\bigcup_{i\in\ZZ} R^i(A)\times C & \longrightarrow   &  \bigcup_{i\in\ZZ} R^i(A)\times C \\
(x,c)  & \longmapsto & (R(x), h_{x}(c))
\end{eqnarray*}
where $(h_{x})_{x \in \cup R^i(A)}$ is a family of homeomorphisms of $C$. We make the following continuity assumption: for every integer $i$, the map $h^i$ is continuous on $A\times C$. 

\begin{rema*}
In the case $A = \cM$, the continuity assumption implies that $h$ is a homeomorphism of $\cM \times C$; if $\cM$ is connected, then $h$ has to be a product as in theorem~\ref{t.realisability}.  Thus we do not restrict ourselves to the case where $A = \cM$; we rather think of $A$ as a Cantor set in $\cM$. Also note that the continuity assumption amounts to requiring that $h$ is continuous on $R^i(A) \times C$ for every $i$. Note that this does \emph{not} imply that $h$ is continuous on $\bigcup_{i\in\ZZ} R^i(A)\times C$.
\end{rema*}

We will prove the following general statement.

\begin{theo}
\label{t.main}
Let $R$ be a homeomorphism on a compact manifold $\cM$ of dimension $d\geq 2$. Let $\mu$ be an aperiodic\footnote{An invariant measure $\mu$ is \emph{aperiodic} if the set of periodic points of $R$ has measure $0$ for $\mu$.  When $\mu$ is ergodic, this just says that $\mu$ is not supported by a periodic orbit of $R$.} ergodic measure for $R$, and $A\subset\cM$ be a set which has positive measure for $\mu$ and has zero-measure for every other ergodic $R$-invariant measure. Let $h$ be a bijective map which is fibered over $R$ and satisfies the continuity assumption as above.
Then there exists a homeomorphism $f:\cM\to\cM$ such that $(\cM,f)$ is isomorphic to the disjoint union
$$\left(\bigcup_{i\in\ZZ} R^i(A)\times C,h\right) \bigsqcup
 \left(\cM\setminus\bigcup_{i\in\ZZ} R^i(A),R\right).$$
Furthermore, $f$ is a topological extension of $R$: there exists a continuous map $\Phi:\cM\to\cM$, which is one-to-one outside the set $\Phi^{-1}(\supp\mu)$, such that $\Phi\circ f=R\circ\Phi$. If $R$ is minimal (resp. transitive), then $f$ can be chosen to be minimal (resp. transitive).   
\end{theo}

The case where $R$ is not minimal on $\cM$ is considered  in the following addendum.

\begin{adde}
\label{a.main}
In any case, the dynamics $f$ is transitive on $\Phi^{-1}(\supp\mu)$. Moreover, if $R$ is minimal on $\supp\mu$, then $f$ can be chosen to be minimal on $\Phi^{-1}(\supp\mu)$. 
\end{adde}

In the case where $R$ is uniquely ergodic, theorem~\ref{t.main} asserts that  there exists a homeomorphism $f:\cM\to\cM$ which is isomorphic to the fibered map $h$. In the general case, it roughly says that  there exists a homeomorphism $f:\cM\to\cM$ which, from the measurable point of view, ``looks like" $R$ outside $\Phi^{-1}\left(\bigcup_{i\in\ZZ} R^i(A)\right)$ and like $h$ on $\Phi^{-1}\left(\bigcup_{i\in\ZZ} R^i(A)\right)$. In other words, it allows to replace the dynamics of $R$ on the iterates of $A$ by the dynamics of $h$. 

\bigskip

Our main motivation for considering maps that fibre over a homeomorphism $R$ but are not direct product is the forthcoming paper~\cite{BegCroLeR}.
Indeed, in that paper,  we will prove that  any measured extension of a uniquely ergodic homeomorphism $R$  is metrically conjugate to fibered map over $R$ satisfying a continuity assumption as above. Using theorem~\ref{t.main}, this will allow us to construct some uniquely ergodic homeomorphism metrically conjugate to any measured extension of $R$. And, choosing $R$ carefully, this will allow us to construct a strictly ergodic homeomorphism  $f$ on $\TT^2$ metrically conjugate to a given ergodic dynamical system $(X,T,\mu)$, provided that the pointwise spectrum of $(X,T,\mu)$ contains an irrational number.

To conclude, here is an example where we use theorem~\ref{t.main} with a homeomorphism $R$ that is not uniquely ergodic. Consider a euclidean rotation of irrational angle on the sphere $\SS^2$. Let $\mu$ be the ergodic measure supported by some invariant circle $A$. Let $h$ be the product of $R \mid _{A}$ with some (positive entropy) homeomorphism of the Cantor space $C$. 
Then the map $f$ provided by theorem~\ref{t.main} is an irrational pseudo-rotation of the $2$-sphere with positive topological entropy. The entropy is concentrated on the minimal invariant set $\Phi^{-1}(A)$. This set is one-dimensional (connected with empty interior), and it cuts the sphere into two  open half-spheres on which $f$ is $C^\infty$ conjugate to the initial rotation $R$. In other words, $f$ is a ``rotation with an invariant circle  replaced by some minimal set with a wilder dynamical behaviour". 
The possibility of  constructing such a homeomorphism was mentionned in~\cite{SanWal};
more details are given in appendix~\ref{ss.pseudo-rotation}.
It is interesting to compare this homeomorphisms with the examples constructed by M. Handel and Herman in~\cite{Han, Her2}.

\subsection{Outline of Denjoy-Rees technique}
\label{ss.outline} 

In this section, we give an idea of the Denjoy-Rees technique. For this purpose,
we first recall one particular construction of the famous Denjoy homeomorphism on the circle. Many features of Rees construction already appear in this presentation of Denjoy
construction, especially the use of microscopic perturbations (allowing the convergence of the construction) with macroscopic effects on the dynamics.

\paragraph{General method for constructing Denjoy counter-examples.}
To construct a Denjoy counter-example on the circle, one starts with an irrational rotation and blows up the orbit of some point to get an orbit of wandering intervals.
There are several ways to carry out the construction, let us outline the one that suits our needs. We first choose an irrational rotation $R$. The homeomorphism $f$ is obtained as a limit of conjugates of $R$,
$$
f = \lim f_{n} \mbox{ with } f_{n} = \Phi_{n}^{-1} R \Phi_{n}
$$
where the sequence of homeomorphisms $(\Phi_{n})$ will converge towards a non-invertible continuous map $\Phi$ that will provide a semi-conjugacy between $f$ and $R$.
To construct this sequence, 
we pick some interval $I_{0}$ (that will become a wandering interval).
The map $\Phi_{n}$ will map $I_{0}$  to some interval $I_{n}$
 which is getting smaller and smaller as $n$ increases, so that $I_{0}$ is ``more and more wandering''; more precisely, $I_{0}$ is disjoint from its $n$ first iterates under $f_{n}$.

In order to make the sequence $(f_{n})$ converge, the construction is done recursively.
The map $\Phi_{n+1}$ is obtained by post-composing $\Phi_{n}$ with some homeomorphism $M_{n+1}$ that maps $I_{n}$ on $I_{n+1}$, and 
whose support is the disjoint union of the $n$  backward and forward  iterates, under the rotation $R$,  of an interval $\hat I_{n}$ which is slightly bigger than $I_{n}$. Thus the uniform distance from $\Phi_{n}$ to $\Phi_{n+1}$ is roughly equal to the size of $I_{n}$.  This guarantees the convergence of the sequence $(\Phi_{n})$ if the size of $I_{n}$ tends to $0$ quickly enough. 

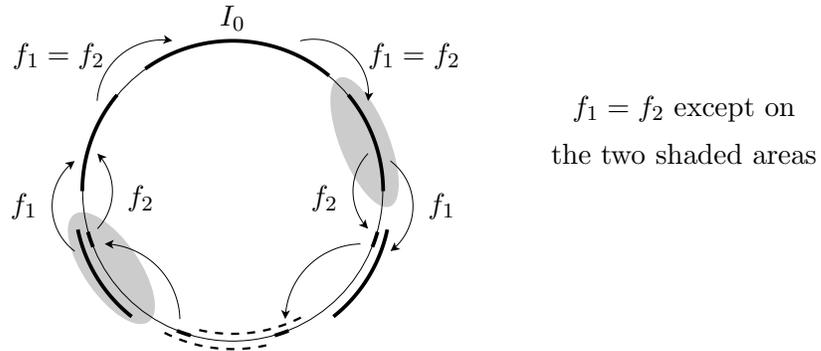
\begin{figure}[htbp]
\hspace{-2cm}
\ifx\JPicScale\undefined\def\JPicScale{1}\fi
\psset{unit=\JPicScale mm}
\psset{linewidth=0.3,dotsep=1,hatchwidth=0.3,hatchsep=1.5,shadowsize=1}
\psset{dotsize=0.7 2.5,dotscale=1 1,fillcolor=black}
\psset{arrowsize=1 2,arrowlength=1,arrowinset=0.25,tbarsize=0.7 5,bracketlength=0.15,rbracketlength=0.15}
\begin{pspicture}(0,0)(140,47.1)
\newrgbcolor{userFillColour}{0.8 0.8 0.8}
\rput{111.53}(97.55,26.53){\psellipse[linewidth=0,linestyle=none,fillcolor=userFillColour,fillstyle=solid](0,0)(9.2,-3.28)}
\newrgbcolor{userFillColour}{0.8 0.8 0.8}
\rput{124.64}(63.89,9.77){\psellipse[linewidth=0,linestyle=none,fillcolor=userFillColour,fillstyle=solid](0,0)(8.68,-3.68)}
\rput{0}(80,20){\psellipse[linewidth=0.1](0,0)(20,20)}
\rput{0}(80,20){\parametricplot[linewidth=0.5,arrows=-]{-125.54}{-50.19}{ t cos 20 mul t sin -20 mul }}
\rput{0}(80,20){\parametricplot[linewidth=0.5,arrows=-]{180}{219.81}{ t cos 20 mul t sin -20 mul }}
\rput{0}(79,19){\parametricplot[linewidth=0.5,arrows=-]{128.38}{168.19}{ t cos 20 mul t sin -20 mul }}
\rput{0}(80,19){\parametricplot[linestyle=dashed,dash=1 1,arrows=-]{77.13}{116.94}{ t cos 20 mul t sin -20 mul }}
\rput{0}(80,20){\parametricplot[linewidth=0.5,arrows=-]{-0}{39.81}{ t cos 20 mul t sin 20 mul }}
\rput{0}(81,19){\parametricplot[linewidth=0.5,arrows=-]{-51.62}{-11.81}{ t cos 20 mul t sin 20 mul }}
\rput{0}(80,21){\parametricplot[linestyle=dashed,dash=1 1,arrows=-]{-102.87}{-63.06}{ t cos 20 mul t sin 20 mul }}
\rput{0}(80,20){\parametricplot[linewidth=0.5,arrows=-]{106.19}{111.8}{ t cos 20 mul t sin -20 mul }}
\rput{0}(80,20){\parametricplot[linewidth=0.5,arrows=-]{-73.81}{-68.2}{ t cos 20 mul t sin 20 mul }}
\rput{0}(80,20){\parametricplot[linewidth=0.5,arrows=-]{158.2}{164.29}{ t cos 20 mul t sin -20 mul }}
\rput{0}(80,20){\parametricplot[linewidth=0.5,arrows=-]{-21.8}{-15.71}{ t cos 20 mul t sin 20 mul }}
\rput{0}(70.67,31.42){\parametricplot[linewidth=0.1,arrows=<-]{81.17}{176.15}{ t cos 8.69 mul t sin 8.69 mul }}
\rput{0}(63.5,18){\parametricplot[linewidth=0.1,arrows=<-]{126.87}{233.13}{ t cos 7.5 mul t sin 7.5 mul }}
\rput{0}(90.76,32.92){\parametricplot[linewidth=0.1,arrows=<-]{-7.25}{103.99}{ t cos 7.3 mul t sin 7.3 mul }}
\rput{0}(96.5,18){\parametricplot[linewidth=0.1,arrows=<-]{-53.13}{53.13}{ t cos 7.5 mul t sin 7.5 mul }}
\rput{0}(56.75,20){\parametricplot[linewidth=0.1,arrows=->]{-43.6}{43.6}{ t cos 7.25 mul t sin 7.25 mul }}
\rput{0}(103.25,20){\parametricplot[linewidth=0.1,arrows=->]{136.4}{223.6}{ t cos 7.25 mul t sin 7.25 mul }}
\rput[r](54,18){$f_1$}
\rput[l](106,18){}
\rput[r](107,16){}
\rput[l](106,18){$f_1$}
\rput[l](66,19){$f_2$}
\rput[r](94,19){$f_2$}
\rput[r](63,38){$f_1=f_2$}
\rput[l](98,38){$f_1=f_2$}
\rput{0}(97.25,2.75){\parametricplot[linewidth=0.1,arrows=->]{91.4}{178.6}{ t cos 10.25 mul t sin 10.25 mul }}
\rput{0}(62.75,2.75){\parametricplot[linewidth=0.1,arrows=->]{1.4}{88.6}{ t cos 10.25 mul t sin 10.25 mul }}
\rput(80,43){$I_0$}
\rput(140,31){$f_1=f_2$ except on}
\rput(140,25){the two shaded areas}
\end{pspicture}
\caption{Denjoy construction on the circle: the interval $I_{0}$ 
has  four consecutive disjoint iterates under $f_{1}$, and  six consecutive disjoint iterates under $f_{2}$.}
\label{f.denjoy}
\end{figure}

Clearly, this is not enough for the convergence of $(f_{n})$ (for example, if $\Phi_{n}$ was the identity outside a little neighbourhood of $I_{0}$, then $(f_{n})$ could converge to a map that crashes $I_{0}$ onto a point). We also demand that $M_{n+1}$ commutes with the rotation $R$ except on the union of two small intervals $I_{n}^{\mathrm{in}}$ and $I_{n}^{\mathrm{out}}$ (namely, $I_{n}^{\mathrm{in}} = R^{-(n+1)}(\hat I_{n})$ and $I_{n}^{\mathrm{out}} = R^{n}(\hat I_{n})$). Then an immediate computation
 shows that the map $f_{n+1}$ coincides with $f_{n}$ except on the set $\Phi_{n}^{-1}(I_{n}^{\mathrm{in}} \cup I_{n}^{\mathrm{out}})$, which happens to be equal to $\Phi_{n-1}^{-1}(I_{n}^{\mathrm{in}} \cup I_{n}^{\mathrm{out}})$ (due to the  condition on the support of $M_{n}$, see figure ~\ref{f.denjoy}). Note that the interval  $\hat I_{n}$ \emph{is chosen after the map $\Phi_{n-1}$ has been designed}, and thus we see that this set $\Phi_{n-1}^{-1}(I_{n}^{\mathrm{in}} \cup I_{n}^{\mathrm{out}})$ can have been made arbitrarily small (by choosing $\hat I_n$ small enough), so that $f_{n+1}$ is arbitrarily close to $f_{n}$.

\paragraph{The Denjoy-Rees technique.}
We now turn to the generalisation of the Denjoy construction developed by Rees. We have in mind the easiest setting: the map $R$ is an irrational rotation of the two-torus ${\cM}= \TT^2$, we are given some homeomorphism $h_{C}$ on some abstract Cantor space $C$, and we want to construct a minimal homeomorphism $f$ of $\TT^2$ which is isomorphic to $R\times h_C$. In some sense, we aim to blow up the dynamics of $R$ and to ``embed the dynamics of $h_{C}$'' into the blown-up homeomorphism $f$. The main difference with the Denjoy construction is that we  have to blow up the orbits of all the points of a positive measure Cantor set $K$. Whereas for the previous construction the point to be blown up was disjoint from all its iterates under the rotation $R$, obviously \emph{$K$ will  meet some of its iterates}: one has to deal with the recurrence of $K$, which adds considerable difficulty.

At step $k$, we will know an approximation $E_{k}^0$ of $K$ by a union of small rectangles. The key property of these rectangles, that enables the construction in spite of the recurrence, is that they are \emph{dynamically coherent}: if $X_1, X_2$ are two connected components of $E_k^0$, if $-k \leq l,l'  \leq k$, then the rectangles $R^l(X_{1})$ and $R^{l'}(X_{2})$ are either disjoint or equal.

\begin{figure}[htbp]
\def\JPicScale{1.3}
\hspace{-1.5cm}
\ifx\JPicScale\undefined\def\JPicScale{1}\fi
\psset{unit=\JPicScale mm}
\psset{linewidth=0.3,dotsep=1,hatchwidth=0.3,hatchsep=1.5,shadowsize=1}
\psset{dotsize=0.7 2.5,dotscale=1 1,fillcolor=black}
\psset{arrowsize=1 2,arrowlength=1,arrowinset=0.25,tbarsize=0.7 5,bracketlength=0.15,rbracketlength=0.15}
\begin{pspicture}(0,0)(124,90)
\rput(15,20){}
\rput{0}(85,70){\psellipse[linewidth=0.1](0,0)(35,-20)}
\rput{0}(85,20){\psellipse[linewidth=0.1](0,0)(35,-20)}
\psline[linewidth=0.4](70,69)(70,67.89)
\psline[linewidth=0.4](70,66.78)(70,65.67)
\psline[linewidth=0.4](70,62.33)(70,61.22)
\psline[linewidth=0.4](70,60.11)(70,59)
\psline[linewidth=0.4](70,85)(70,83.89)
\psline[linewidth=0.4](70,82.78)(70,81.67)
\psline[linewidth=0.4](70,78.33)(70,77.22)
\psline[linewidth=0.4](70,76.11)(70,75)
\psline[linewidth=0.4,fillstyle=solid](80,69)(80,67.89)
\psline[linewidth=0.4,fillstyle=solid](80,66.78)(80,65.67)
\psline[linewidth=0.4,fillstyle=solid](80,62.33)(80,61.22)
\psline[linewidth=0.4,fillstyle=solid](80,60.11)(80,59)
\psline[linewidth=0.4,fillstyle=solid](80,85)(80,83.89)
\psline[linewidth=0.4,fillstyle=solid](80,82.78)(80,81.67)
\psline[linewidth=0.4,fillstyle=solid](80,78.33)(80,77.22)
\psline[linewidth=0.4,fillstyle=solid](80,76.11)(80,75)
\psline[linewidth=0.1,fillstyle=solid](84,69)(84,67.89)
\psline[linewidth=0.1,fillstyle=solid](84,66.78)(84,65.67)
\psline[linewidth=0.1,fillstyle=solid](84,62.33)(84,61.22)
\psline[linewidth=0.1,fillstyle=solid](84,60.11)(84,59)
\psline[linewidth=0.1,fillstyle=solid](84,85)(84,83.89)
\psline[linewidth=0.1,fillstyle=solid](84,82.78)(84,81.67)
\psline[linewidth=0.1,fillstyle=solid](84,78.33)(84,77.22)
\psline[linewidth=0.1,fillstyle=solid](84,76.11)(84,75)
\psline[linewidth=0.1,fillstyle=solid](81,69)(81,67.89)
\psline[linewidth=0.1,fillstyle=solid](81,66.78)(81,65.67)
\psline[linewidth=0.1,fillstyle=solid](81,62.33)(81,61.22)
\psline[linewidth=0.1,fillstyle=solid](81,60.11)(81,59)
\psline[linewidth=0.1,fillstyle=solid](81,85)(81,83.89)
\psline[linewidth=0.1,fillstyle=solid](81,82.78)(81,81.67)
\psline[linewidth=0.1,fillstyle=solid](81,78.33)(81,77.22)
\psline[linewidth=0.1,fillstyle=solid](81,76.11)(81,75)
\psline[linewidth=0.1,fillstyle=solid](85,69)(85,67.89)
\psline[linewidth=0.1,fillstyle=solid](85,66.78)(85,65.67)
\psline[linewidth=0.1,fillstyle=solid](85,62.33)(85,61.22)
\psline[linewidth=0.1,fillstyle=solid](85,60.11)(85,59)
\psline[linewidth=0.1,fillstyle=solid](85,85)(85,83.89)
\psline[linewidth=0.1,fillstyle=solid](85,82.78)(85,81.67)
\psline[linewidth=0.1,fillstyle=solid](85,78.33)(85,77.22)
\psline[linewidth=0.1,fillstyle=solid](85,76.11)(85,75)
\psline[linewidth=0.15]{->}(70,58)(70,24)
\psline[linewidth=0.15]{->}(80,57.5)(80,20.5)
\psline[linewidth=0.05]{->}(71,24.5)(83.5,32.5)
(83.5,32.5)(94.5,28)
(94.5,28)(105,30.5)
(105,30.5)(114,24)
(114,24)(113.5,16.5)
(113.5,16.5)(101,16)
(101,16)(92.5,12.5)
(92.5,12.5)(82.5,18)
\psline[linewidth=0.05]{->}(70.5,86)(85,89)
(85,89)(95,84)
(95,84)(105,85)
(105,85)(114,73)
(114,73)(114,62)
(114,62)(102,61)
(102,61)(89,53)
(89,53)(81,58)
\rput[l](104.5,58.5){$f^p$}
\rput[l](104,13){$R^p$}
\rput[l](106,45){$\Phi$}
\rput[r](69.5,23){$\scriptstyle x$}
\rput[r](79.5,20){$\scriptstyle y$}
\rput[r](69.5,71){\rotatebox{90}{$\scriptstyle\{x\}\times C$}}
\rput(78,71){\rotatebox{90}{$\scriptstyle \{y\}\times C$}}
\psline[linewidth=0.1,fillstyle=solid](74,69)(74,67.89)
\psline[linewidth=0.1,fillstyle=solid](74,66.78)(74,65.67)
\psline[linewidth=0.1,fillstyle=solid](74,62.33)(74,61.22)
\psline[linewidth=0.1,fillstyle=solid](74,60.11)(74,59)
\psline[linewidth=0.1,fillstyle=solid](74,85)(74,83.89)
\psline[linewidth=0.1,fillstyle=solid](74,82.78)(74,81.67)
\psline[linewidth=0.1,fillstyle=solid](74,78.33)(74,77.22)
\psline[linewidth=0.1,fillstyle=solid](74,76.11)(74,75)
\psline[linewidth=0.1,fillstyle=solid](71,69)(71,67.89)
\psline[linewidth=0.1,fillstyle=solid](71,66.78)(71,65.67)
\psline[linewidth=0.1,fillstyle=solid](71,62.33)(71,61.22)
\psline[linewidth=0.1,fillstyle=solid](71,60.11)(71,59)
\psline[linewidth=0.1,fillstyle=solid](71,85)(71,83.89)
\psline[linewidth=0.1,fillstyle=solid](71,82.78)(71,81.67)
\psline[linewidth=0.1,fillstyle=solid](71,78.33)(71,77.22)
\psline[linewidth=0.1,fillstyle=solid](71,76.11)(71,75)
\psline[linewidth=0.1,fillstyle=solid](75,69)(75,67.89)
\psline[linewidth=0.1,fillstyle=solid](75,66.78)(75,65.67)
\psline[linewidth=0.1,fillstyle=solid](75,62.33)(75,61.22)
\psline[linewidth=0.1,fillstyle=solid](75,60.11)(75,59)
\psline[linewidth=0.1,fillstyle=solid](75,85)(75,83.89)
\psline[linewidth=0.1,fillstyle=solid](75,82.78)(75,81.67)
\psline[linewidth=0.1,fillstyle=solid](75,78.33)(75,77.22)
\psline[linewidth=0.1,fillstyle=solid](75,76.11)(75,75)
\psline[linewidth=1,fillstyle=solid](24,69)(24,67.89)
\psline[linewidth=1,fillstyle=solid](24,66.78)(24,65.67)
\psline[linewidth=1,fillstyle=solid](24,62.33)(24,61.22)
\psline[linewidth=1,fillstyle=solid](24,60.11)(24,59)
\psline[linewidth=1,fillstyle=solid](24,85)(24,83.89)
\psline[linewidth=1,fillstyle=solid](24,82.78)(24,81.67)
\psline[linewidth=1,fillstyle=solid](24,78.33)(24,77.22)
\psline[linewidth=1,fillstyle=solid](24,76.11)(24,75)
\psline[linewidth=0.5]{->}(65,73)(27,73)
\rput[r](48.5,75.5){$C \leftarrow K \times C$}
\rput[r](15.5,70){$h_{C}^p$}
\rput[r](5,60){}
\rput(24,54){$C$}
\rput(75,56){$K \times C$}
\rput[r](14,41){}
\rput{0}(102,72.5){\parametricplot[linewidth=0.1,arrows=-]{-153.43}{-26.57}{ t cos 5.59 mul t sin 5.59 mul }}
\rput{0}(102,66){\parametricplot[linewidth=0.1,arrows=-]{36.87}{143.13}{ t cos 5 mul t sin 5 mul }}
\rput{0}(102.41,24.5){\parametricplot[linewidth=0.1,arrows=-]{-153.43}{-26.57}{ t cos 5.59 mul t sin 5.59 mul }}
\rput{0}(102.41,18){\parametricplot[linewidth=0.1,arrows=-]{36.87}{143.13}{ t cos 5 mul t sin 5 mul }}
\rput(77.5,15){$K$}
\rput[r](18,52){}
\psline[linewidth=0.05]{->}(22,84)(17,80)
(17,80)(15,75)
(15,75)(18,71)
(18,71)(16,68)
(16,68)(21,62)
\psline[linewidth=0.5]{->}(104,51)(104,38)
\rput{0}(117.08,81.23){\parametricplot[linewidth=0.5,arrows=<-]{-53.13}{118.07}{ t cos 6.54 mul t sin 6.54 mul }}
\rput(124,88){$f$}
\rput{0}(117.08,31.23){\parametricplot[linewidth=0.5,arrows=<-]{-53.13}{118.07}{ t cos 6.54 mul t sin 6.54 mul }}
\rput(124,38){$R$}
\rput[l](116,7){$\cal M$}
\rput[l](115.5,56.5){$\wt{\cal M} = {\cal M}$}
\rput(122,15){}
\pspolygon[linewidth=0.1,fillstyle=solid](70,23.5)(70.56,23.5)(70.56,22.98)(70,22.98)
\pspolygon[linewidth=0.1,fillstyle=solid](70,22.46)(70.56,22.46)(70.56,21.94)(70,21.94)
\pspolygon[linewidth=0.1,fillstyle=solid](71.11,22.46)(71.67,22.46)(71.67,21.94)(71.11,21.94)
\pspolygon[linewidth=0.1,fillstyle=solid](71.11,23.5)(71.67,23.5)(71.67,22.98)(71.11,22.98)
\pspolygon[linewidth=0.1,fillstyle=solid](73.33,20.39)(73.89,20.39)(73.89,19.87)(73.33,19.87)
\pspolygon[linewidth=0.1,fillstyle=solid](73.33,19.35)(73.89,19.35)(73.89,18.83)(73.33,18.83)
\pspolygon[linewidth=0.1,fillstyle=solid](74.44,19.35)(75,19.35)(75,18.83)(74.44,18.83)
\pspolygon[linewidth=0.1,fillstyle=solid](74.44,20.39)(75,20.39)(75,19.87)(74.44,19.87)
\pspolygon[linewidth=0.1,fillstyle=solid](70,20.39)(70.56,20.39)(70.56,19.87)(70,19.87)
\pspolygon[linewidth=0.1,fillstyle=solid](70,19.35)(70.56,19.35)(70.56,18.83)(70,18.83)
\pspolygon[linewidth=0.1,fillstyle=solid](71.11,19.35)(71.67,19.35)(71.67,18.83)(71.11,18.83)
\pspolygon[linewidth=0.1,fillstyle=solid](71.11,20.39)(71.67,20.39)(71.67,19.87)(71.11,19.87)
\pspolygon[linewidth=0.1,fillstyle=solid](73.33,23.5)(73.89,23.5)(73.89,22.98)(73.33,22.98)
\pspolygon[linewidth=0.1,fillstyle=solid](73.33,22.46)(73.89,22.46)(73.89,21.94)(73.33,21.94)
\pspolygon[linewidth=0.1,fillstyle=solid](74.44,22.46)(75,22.46)(75,21.94)(74.44,21.94)
\pspolygon[linewidth=0.1,fillstyle=solid](74.44,23.5)(75,23.5)(75,22.98)(74.44,22.98)
\pspolygon[linewidth=0.1,fillstyle=solid](80,23.5)(80.56,23.5)(80.56,22.98)(80,22.98)
\pspolygon[linewidth=0.1,fillstyle=solid](80,22.46)(80.56,22.46)(80.56,21.94)(80,21.94)
\pspolygon[linewidth=0.1,fillstyle=solid](81.11,22.46)(81.67,22.46)(81.67,21.94)(81.11,21.94)
\pspolygon[linewidth=0.1,fillstyle=solid](81.11,23.5)(81.67,23.5)(81.67,22.98)(81.11,22.98)
\pspolygon[linewidth=0.1,fillstyle=solid](83.33,20.39)(83.89,20.39)(83.89,19.87)(83.33,19.87)
\pspolygon[linewidth=0.1,fillstyle=solid](83.33,19.35)(83.89,19.35)(83.89,18.83)(83.33,18.83)
\pspolygon[linewidth=0.1,fillstyle=solid](84.44,19.35)(85,19.35)(85,18.83)(84.44,18.83)
\pspolygon[linewidth=0.1,fillstyle=solid](84.44,20.39)(85,20.39)(85,19.87)(84.44,19.87)
\pspolygon[linewidth=0.1,fillstyle=solid](80,20.39)(80.56,20.39)(80.56,19.87)(80,19.87)
\pspolygon[linewidth=0.1,fillstyle=solid](80,19.35)(80.56,19.35)(80.56,18.83)(80,18.83)
\pspolygon[linewidth=0.1,fillstyle=solid](81.11,19.35)(81.67,19.35)(81.67,18.83)(81.11,18.83)
\pspolygon[linewidth=0.1,fillstyle=solid](81.11,20.39)(81.67,20.39)(81.67,19.87)(81.11,19.87)
\pspolygon[linewidth=0.1,fillstyle=solid](83.33,23.5)(83.89,23.5)(83.89,22.98)(83.33,22.98)
\pspolygon[linewidth=0.1,fillstyle=solid](83.33,22.46)(83.89,22.46)(83.89,21.94)(83.33,21.94)
\pspolygon[linewidth=0.1,fillstyle=solid](84.44,22.46)(85,22.46)(85,21.94)(84.44,21.94)
\pspolygon[linewidth=0.1,fillstyle=solid](84.44,23.5)(85,23.5)(85,22.98)(84.44,22.98)
\rput(42.5,41){$f^p(\{x\}\times C)=\{y\}\times C$}
\rput(42.5,31.5){$R^p(x)=y$}
\rput(2.5,1){}
\end{pspicture}
\caption{Rees construction on the torus: description of the isomorphism between $f$ and $R \times h_{C}$. The topological flavour is given by the vertical map $\Phi$, the measurable flavour is given by the horizontal map $K \times C \to C$.}
\label{f.isomorphism}
\end{figure}
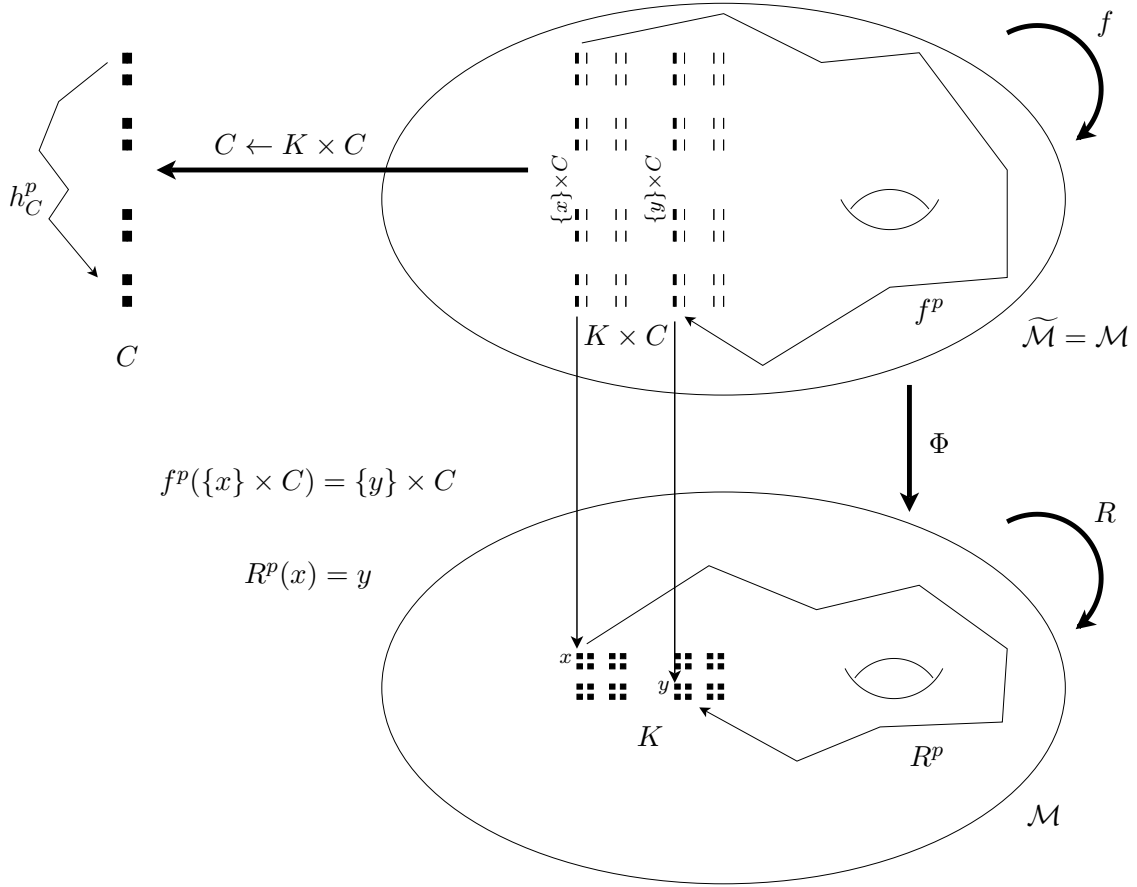

The final map $f$ will again be a limit of conjugates of $R$,
$$
f = \lim f_{k} \mbox{ with } f_{k} = \Phi_{k}^{-1} R \Phi_{k} \mbox{ and } \Phi = \lim \Phi_{k}.
$$
We fix an embedding of the abstract product Cantor set $K \times C$ in the manifold $\cM$ (see figure~\ref{f.isomorphism}); in some sense, this set $K\times C$ will play for the Rees map the same role as the interval $I_{0}$ for the Denjoy map. 
Indeed, each point $x$ of $K$ will be blown-up by $\Phi^{-1}$, so that the fibre $\Phi^{-1}(x)$ contains the ``vertical'' Cantor set $\{x\} \times C$ embedded in $\cM$. Furthermore, the first return map of $f$ in $\wt K :=\Phi^{-1}(K)$ will leave the embedded product Cantor space $K \times C$ invariant: in other words, for each point $x \in K$ and each (return time) $p$ such that $R^p(x)$ belongs to $K$, the homeomorphism $f^p$ will map the vertical Cantor set $\{x\} \times C$ onto the vertical Cantor set $\{R^p(x)\} \times C$. When both vertical Cantor sets are identified to $C$ by way of the second coordinate on $K\times C$, the map $f^p$ induces a homeomorphism of $C$, which will be equal to $h_{C}^p$. This is the way one embeds the dynamics of $h_{C}$ into the dynamics of $f$, and gets an isomorphism between the product $R \times h_{C}$ and the restriction of the map $f$ to the set $\bigcup_{i} f^i(K\times C)$. Note that theorem~\ref{t.realisability} requires more, namely an isomorphism between $R \times h_{C}$ and  the map $f$ \emph{on the whole manifold $\cM$}. We will explain in the last paragraph  how one can further obtain that $f$ is isomorphic to its restriction to the set $\bigcup_{i} f^i(K\times C)$.

Once again the construction will be carried out recursively. At step $k$  we will take care of all return times less than $2k+1$, that is,  the map $\Phi_{k}$ will be constructed so that the approximation $f_{k}$ of $f$ will satisfy the description of figure~\ref{f.isomorphism} for $\mid p \mid \leq 2k+1$. This property will be transmitted to $f_{k+1}$ (and so gradually to $f$) because for any point $x$ of $K$ whose return time in $K$ is less than or equal to $2k+1$, $f_{k+1}^p(x) = f_{k}^p(x)$. Actually the equality $f_{k+1}=f_{k}$ will hold except on a very set which becomes smaller and smaller as $k$ increases (just as in the Denjoy construction). On the other hand  $f_{k+1}$  will take care of return times equal to $2k+2$ and $2k+3$. The convergence of the sequences $(\Phi_{k})$ and $(f_{k})$ will be obtained using essentially the same argument as in the Denjoy construction.

Another feature of the construction is that we want $f$ to inherit from the minimality of $R$. This will be an easy  consequence of the two following properties. Firstly the fibre $\Phi^{-1}(x)$ above a point $x$ that does not belong to an iterate of $K$ will be reduced to a point. Secondly the other fibres will have empty interior.

\paragraph{Control of invariant measures}
Until here, we have been dealing with the control of the dynamics on the iterates of the product Cantor set $K\times C$. This is enough for $f$ to admit the product $R \times h_{C}$ as a subsystem, and to get an example with positive topological entropy (Rees initial result). If we want the much stronger property that $f$ is uniquely ergodic (in theorem~\ref{t.entropy}) or isomorphic to $R \times h_{C}$ (in theorem~\ref{t.realisability}), we need to gain some  control of  the dynamics outside the iterates of $K\times C$, on the whole manifold $\cM$. With this in view, we first note that, since $K$ has positive measure,  the (unique) invariant measure for $R$ gives full measure to $\bigcup_{i}R^i(K)$. The automatic consequence for $f$ is that any invariant measure for $f$ gives full measure to $\bigcup_{i}f^i(\wt K)$ (with $\wt K = \Phi^{-1}(K)$). It now remains to put further constraints on the construction to ensure that 
 any invariant measure for $f$ will give measure $0$ to $\wt K \setminus K \times C$. This will be done by considering the first return map of $f$ in $\wt K$, and by forcing the $\omega$-limit set of any point $x \in \wt K$, with respect to this first return map, to be included in $K\times C$.

\subsection{Structure of the paper}
\label{ss.structure-paper}

Recall that our goal is to prove theorem~\ref{t.main} (which implies theorem~\ref{t.realisability} and theorem~\ref{t.entropy}, see the end of section~\ref{ss.isomorphic}). So we are given, in particular, a homeomorphism $R$ on a manifold $\cM$ and a map $h$ which fibres over $R$, and we aim to construct a homeomorphism $f$ on $\cM$ which is  isomorphic to $\left(\bigcup_{i\in\ZZ} R^i(A)\times C,h\right) \bigsqcup \left(\cM\setminus\bigcup_{i\in\ZZ} R^i(A),R\right).$ 

\bigskip

The paper is divided into three parts.
\begin{itemize}
\item[--] In part~A, we construct a Cantor set $K$, obtained as a decreasing intersections of a sequence sets $(E_n^0)_{n\in\NN}$, where $E_n^0$ is a finite collection of pairwise disjoint rectangles for every $n$. 
\item[--] In part~B,  we explain how to blow-up the orbits of the points of $K$: we construct a sequence of homeomorphisms $(M_{n})$ whose infinite composition is a map $\Psi:\cM \to \cM$ such that $\Psi^{-1}(K)$ contains a copy of the product Cantor set $K \times C$.
\item[--] In part~C, we explain how to insert the dynamics of $h$ in the blowing-up of the orbits of the points of $K$. In order to improve the convergence, one will define an extracted sequence $(\bM_{k})$ of $(M_{n})$. One also needs to ``twist" the dynamics by constructing a sequence of homeomorphisms $(H_k)$ and by replacing each $\bM_k$ by the homeomorphism $H_k\circ \bM_k$. The infinite composition of the homeomorphisms $H_k\circ \bM_k$ is a map $\Phi$ and the desired homeomorphism $f:\cM\to\cM$ is a topological extension of $R$ by $\Phi$.\\
The construction of part C is further divided into two main parts: in section~\ref{s.insert-dynamics}
we obtain Rees theorem (on any manifold);  then we explain in section~\ref{s.control} how to get unique ergodicity.
\end{itemize}
The proof of our main theorem~\ref{t.main} is given at the end of part~C (subsections~\ref{ss.bilan} and~\ref{ss.proof}). The proof of addendum~\ref{a.main} is more technical,  we postpone it to appendix~\ref{a.recurrence}.

In each of the three parts of the paper, we will proceed as follows. First, we introduce some new objects (for example, a sequence of homeomorphisms $(M_n)_{n\geq 1}$) and some hypotheses on these objects (for example, the diameter of every connected component of the support of $M_n$ is less than $2^{-n}$). Second, we prove some consequences of the hypotheses (for example, if the above hypothesis is satisfied, and if we set $\Psi_n=M_n\circ\dots\circ M_1$, then the sequence of maps $(\Psi_n)_{n\geq 1}$ converges). And third, we prove that the hypotheses are realisable (for example, we  construct  a sequence of homeomorphisms $(M_n)$ satisfying the required hypotheses). Most of the time, the main difficulties lie in the ``realisability" results and it might be a good idea to skip the proofs of these results for the first reading (these are mainly subsections~\ref{ss.construction-cantor},~\ref{ss.realis-B_12456},~\ref{ss.construc-n},~\ref{ss.realis-C_1256}, and~\ref{ss.control-existence}).

We have already mentioned the main novelties with respect to Rees paper: we can enrich the dynamics of any homeomorphism (not only rotations of the torus) and we can control the measurable dynamics that have been inserted. Apart from this, the general scheme of our construction is quite similar to Rees original one. However, let us point out a couple of  differences.
\begin{itemize}
\item[--] We have tried, as much as possible, to divide the construction into  independent steps: the Cantor set $K$ is defined once and for all in part A, the sequence $(M_{n})$ is constructed once and for all in part B. Then (in part C) an extraction process will be used to get the convergence.
As already explained above, each step is structured according to a fixed pattern (objects / hypotheses / consequences / construction). 
\item[--] In order to insert the dynamics of the model $h$ in the blowing-up of the orbits of the points of $K$, we notice that it is enough to consider the first-return map of $f$ to the blowing-up of $K$. This avoids many technical complications.
\end{itemize}

\newpage
\part{Construction of a Cantor set $K$}
\label{p.construction-rectangles}

In this part, we assume that we are given a compact topological manifold $\cM$ and a homeomorphism $R:\cM\to\cM$. For technical purposes, we choose a metric on $\cM$. We will explain how to construct a Cantor set $K\subset\cM$ which has a nice behaviour with respect to the action of the homeomorphism $R$. In part B, this nice behaviour will allow us to ``blow-up" the orbits of the points of $K$. The construction described below is a generalisation of Rees construction which works in the case where $R$ is an irrational rotation of the torus $\TT^d$ (see \cite{Ree}). Since rotations of $\TT^d$ are products of rotations of $\SS^1$, Rees original construction essentially takes place in dimension $1$, and is much easier than the present one. The reader mainly  interested in the rotation case can refer to~\cite{Ree} instead of section~\ref{ss.construction-cantor} below.

\section{Dynamically coherent Cantor sets}
\label{s.rectangles}

\subsection{Definitions}

Remember that a \emph{Cantor set} is a metrisable totally discontinuous compact topological space without isolated points; any two such spaces are homeomorphic. In what follows, a \emph{rectangle of $\cM$}  is a subset $X$ of $\cM$ homeomorphic to the closed unit ball in $\RR^d$, where $d$ is the dimension of $\cM$. The phrase ``collection of rectangles'' will always refer to a finite family of pairwise disjoint rectangles of $\cM$. Such a collection will be denoted by a calligraphic letter (like ``$\cE$")~; the corresponding straight letter (like ``$E$") will represent the reunion of all the rectangles of the collection. We denote by $R^k({\cE})$ the collection of rectangles whose elements are the images under the map $R^k$ of the elements of $\cE$.

\begin{defi}
Let $\cE$ and $\cF$ be two collections of rectangles.
We say that $\cal F$ \emph{refines} $\cal E$ if
\begin{itemize}
\item[--] every element of $\cE$ contains at least one element of $\cF$;
\item[--] for every elements $X$ of $\cE$ and $Y$ of $\cF$, either $X$ and $Y$ are disjoint or $Y$ is included in the interior of $X$.
\end{itemize}
\end{defi}

The second property is equivalent to  $F \cap \partial E = \emptyset$ (where $\partial E$ is the boundary of $E$).

\begin{defi}
Let $p$ be a positive integer. A collection of rectangles ${\cal E}^0$ is \emph{$p$
times iterable} if for every  element $X,X'$ of ${\cal E}^0$, for every integers $k,l$ such that $|k|, |l| \leq
p$, the rectangles $R^k(X)$ et $R^l(X')$ are disjoint or equal.
\end{defi}

Equivalently, the collection ${\cal E}^0$ is $p$ times iterable if, for every element $X$ of ${\cal E}^0$ and every integer $m$ such that $|m| \leq 2p$, if $R^m(X)$ meets $E^0$, then $R^m(X)$ also belongs to ${\cal E}^0$.

Let $\cE^0$ be a collection of  rectangles which is $p$ times iterable. For any positive integer $n\leq p$, we define the collection of rectangles
$$
{\cal E}^n:=\bigcup_{|k| \leq n}R^k({\cal E}^0)
$$
and we denote by $E^n$ the union of all the elements of $\cE^n$. For any positive integer $n<p$, we consider the oriented graph $\cG(\cE^n)$, whose vertices are the elements of $\cE^n$, and whose edges represent the dynamics of $R$: there is an edge from $X\in\cE^n$ to $X'\in\cE^n$ if and only if $X'=R(X)$. 

\begin{defi}
Let $\cE^0$ be a  $p$ times iterable collection of rectangles, and $n<p$.
The collection ${\cal E}^n$  is said to be \emph{without cycle}
if the graph $\cG(\cE^n)$ has no cycle.
\end{defi}

The most important definition is the following.

\begin{defi}
\label{d.compatibility}
Let ${\cal E}^0$ and $\cF^0$ be two collections of rectangles, such that $\cE^0$ is $p$ times iterable, and ${\cal F}^0$ is $p+1$ times iterable.  Assume that  $\cF^{p+1}$ refines $\cE^p$. We say that ${\cal F}^0$ is \emph{compatible with ${\cal E}^0$ for $p$ iterates} if $F^0\subset E^0$ and
$$
\mbox{for every $k$ such that $|k|\leq 2p+1$, we have $R^k(F^{0})\cap E^{0}\subset F^{0}.$}
$$
\end{defi}

This definition can be reformulated in several different ways. For example, one can check that $\cF^0$ is compatible with $\cE^0$ for $p$ iterates if and only if $F^0 \subset E^0$ and 
\begin{equation}
 \mbox{for every $k$ such that $|k|\leq p$, we have }  F^{p+1}\cap R^k(E^0)=R^k(F^0).
\end{equation}
Indeed, assume that $F^0\subset E^0$ and consider an integer $k$ such that $|k|\leq p$. 
 Then the reverse inclusion of (1), $F^{p+1}\cap R^k(E^0)\supset R^k(F^0)$, is always satisfied. Furthermore, one has
$$
 \bigcup_{|k'| \leq 2p+1} R^{k'}(F^0) =  \bigcup_{|k| \leq p} R^{-k}(F^{p+1}).
$$
On the one hand, the compatibility says that the intersection of the left-hand side with $E_{0}$ is included in $F_{0}$. On the other hand, property (1) says that the  intersection of the right-hand side  with $E_{0}$ is included in $F_{0}$. One deduces that the compatibility is equivalent to property (1).

\subsection{Hypotheses $\mathbf{A_{1,2,3}}$}

Let $({\cal E}^0_n)_{n\in \NN}$ be a sequence of collections of rectangles. We introduce the following hypotheses.

\begin{itemize}
\item[]
\begin{itemize}
\item[$\mathbf{A_{1}}$] \emph{(Combinatorics of  rectangles)}
\begin{itemize}
\item[$\mathbf{a}$] For every $n\in \NN$, the collection ${\cal E}^0_n$ is $n+1$ times iterable and the collection $\cE_n^{n}$ has no cycle;
\item[$\mathbf{b}$] for every $0 \leq m<n$, the collection ${\cal E}^{n+1}_{n}$
refines the collection ${\cal E}^{m+1}_{m}$;
\item[$\mathbf{c}$] for every $n\in \NN$, the collection ${\cal E}^0_{n+1}$ is
compatible with ${\cal E}^0_n$ for $n+1$ iterates.
\end{itemize}
\item[$\mathbf{A_{2}}$] \emph{(No isolated point)}\\
For every $n\in\NN$ and every rectangle $X\in{\cal E}^0_{n}$,  there are at least two elements of ${\cal E}^0_{n+1}$ contained in~$X$. 
\item[$\mathbf{A_{3}}$] \emph{(Decay of the collections of  rectangles)}\\
 The supremum of the diameters of the elements of $\cE^n_{n}$ tends to $0$ when $n\to+\infty$.
\end{itemize}
\end{itemize}
For sake of simplicity, we will often assume that $\cE_{0}^0$ contains a single rectangle $X_{0}$.

\begin{figure}[htbp]
\hspace{-1.cm}
\ifx\JPicScale\undefined\def\JPicScale{1}\fi
\psset{unit=\JPicScale mm}
\psset{linewidth=0.3,dotsep=1,hatchwidth=0.3,hatchsep=1.5,shadowsize=1}
\psset{dotsize=0.7 2.5,dotscale=1 1,fillcolor=black}
\psset{arrowsize=1 2,arrowlength=1,arrowinset=0.25,tbarsize=0.7 5,bracketlength=0.15,rbracketlength=0.15}
\begin{pspicture}(0,0)(155,42)
\pspolygon[linewidth=0.4](65,42)(105,42)(105,2)(65,2)
\pspolygon[linewidth=0.2](90,40)(100,40)(100,30)(90,30)
\pspolygon[linewidth=0.2](70,30)(80,30)(80,20)(70,20)
\pspolygon[linewidth=0.2](90,20)(100,20)(100,10)(90,10)
\pspolygon[linewidth=0.2](130,40)(140,40)(140,30)(130,30)
\pspolygon[linewidth=0.2](110,30)(120,30)(120,20)(110,20)
\pspolygon[linewidth=0.2](130,20)(140,20)(140,10)(130,10)
\pspolygon[linewidth=0.2](50,40)(60,40)(60,30)(50,30)
\pspolygon[linewidth=0.2](50,20)(60,20)(60,10)(50,10)
\pspolygon[linewidth=0.2](30,30)(40,30)(40,20)(30,20)
\pspolygon[linewidth=0.1](72,28)(74,28)(74,26)(72,26)
\pspolygon[linewidth=0.1](76,24)(78,24)(78,22)(76,22)
\pspolygon[linewidth=0.1](32,28)(34,28)(34,26)(32,26)
\pspolygon[linewidth=0.1](36,24)(38,24)(38,22)(36,22)
\pspolygon[linewidth=0.1](52,38)(54,38)(54,36)(52,36)
\pspolygon[linewidth=0.1](56,34)(58,34)(58,32)(56,32)
\pspolygon[linewidth=0.1](112,28)(114,28)(114,26)(112,26)
\pspolygon[linewidth=0.1](116,24)(118,24)(118,22)(116,22)
\pspolygon[linewidth=0.1](132,38)(134,38)(134,36)(132,36)
\pspolygon[linewidth=0.1](136,34)(138,34)(138,32)(136,32)
\pspolygon[linewidth=0.1](96,18)(98,18)(98,16)(96,16)
\pspolygon[linewidth=0.1](92,14)(94,14)(94,12)(92,12)
\pspolygon[linewidth=0.1](56,18)(58,18)(58,16)(56,16)
\pspolygon[linewidth=0.1](52,14)(54,14)(54,12)(52,12)
\pspolygon[linewidth=0.1](12,28)(14,28)(14,26)(12,26)
\pspolygon[linewidth=0.1](16,24)(18,24)(18,22)(16,22)
\pspolygon[linewidth=0.1](26,18)(28,18)(28,16)(26,16)
\pspolygon[linewidth=0.1](22,14)(24,14)(24,12)(22,12)
\pspolygon[linewidth=0.1](22,38)(24,38)(24,36)(22,36)
\pspolygon[linewidth=0.1](26,34)(28,34)(28,32)(26,32)
\psline[linewidth=0.4]{->}(60,35)(90,35)
\psline[linewidth=0.4]{->}(100,35)(130,35)
\psline[linewidth=0.4]{->}(110,25)(80,25)
\psline[linewidth=0.4]{->}(70,25)(40,25)
\psline[linewidth=0.4]{->}(60,15)(90,15)
\psline[linewidth=0.4]{->}(100,15)(130,15)
\psline[linewidth=0.05]{->}(24,37)(52,37)
\psline[linewidth=0.05]{->}(28,33)(56,33)
\psline[linewidth=0.05]{->}(58,33)(96,33)
\psline[linewidth=0.05]{->}(98,33)(136,33)
\psline[linewidth=0.05]{->}(116,23)(78,23)
\psline[linewidth=0.05]{->}(76,23)(38,23)
\psline[linewidth=0.05]{->}(36,23)(18,23)
\psline[linewidth=0.05]{->}(28,17)(56,17)
\psline[linewidth=0.05]{->}(58,17)(96,17)
\psline[linewidth=0.05]{->}(98,17)(136,17)
\pscustom[linewidth=0.05]{\psline(16,23)(12,23)
\psbezier(12,23)(10,23)(10,22)(10,20)
\psbezier(10,20)(10,18)(10,17)(12,17)
\psbezier(12,17)(14,17)(26,17)(26,17)
}
\psline[linewidth=0.05]{->}(54,37)(92,37)
\psline[linewidth=0.05]{->}(94,37)(132,37)
\psline[linewidth=0.05]{->}(112,27)(74,27)
\psline[linewidth=0.05]{->}(72,27)(34,27)
\psline[linewidth=0.05]{->}(32,27)(14,27)
\psline[linewidth=0.05]{->}(24,13)(52,13)
\psline[linewidth=0.05]{->}(54,13)(92,13)
\psline[linewidth=0.05]{->}(94,13)(132,13)
\psline[linewidth=0.05]{->}(17.5,17)(26,17)
\psline[linewidth=0.05]{->}(87,5)(79,5)
\pscustom[linewidth=0.4]{\psline(140,35)(150,35)
\psbezier(150,35)(155,35)(155,35)(155,30)
\psbezier(155,30)(155,25)(155,25)(150,25)
\psbezier(150,25)(145,25)(120,25)(120,25)
}
\psline[linewidth=0.4]{->}(129,25)(120,25)
\pscustom[linewidth=0.05]{\psline(138,33)(148,33)
\psbezier(148,33)(153,33)(153,33)(153,28)
\psbezier(153,28)(153,23)(153,23)(148,23)
\psbezier(148,23)(143,23)(118,23)(118,23)
}
\psline[linewidth=0.05]{->}(127,23)(118,23)
\pspolygon[linewidth=0.1](92,38)(94,38)(94,36)(92,36)
\pspolygon[linewidth=0.1](96,34)(98,34)(98,32)(96,32)
\pspolygon[linewidth=0.1](136,18)(138,18)(138,16)(136,16)
\pspolygon[linewidth=0.1](132,14)(134,14)(134,12)(132,12)
\pspolygon[linewidth=0.1](81,10)(83,10)(83,8)(81,8)
\pspolygon[linewidth=0.1](77,6)(79,6)(79,4)(77,4)
\pscustom[linewidth=0.05]{\psbezier(138,17)(138,17)(142,17)(145,17)
\psbezier(145,17)(148,17)(148,16)(148,13)
\psbezier(148,13)(148,10)(148,9)(145,9)
\psbezier(145,9)(142,9)(131,9)(127,9)
\psbezier(127,9)(123,9)(83,9)(83,9)
}
\pscustom[linewidth=0.05]{\psline(134,37)(144,37)
\psbezier(144,37)(149,37)(149,37)(149,32)
\psbezier(149,32)(149,27)(149,27)(144,27)
\psbezier(144,27)(139,27)(114,27)(114,27)
}
\psline[linewidth=0.05]{->}(123,27)(114,27)
\pscustom[linewidth=0.05]{\psbezier(134,13)(134,13)(138,13)(141,13)
\psbezier(141,13)(144,13)(144,12)(144,9)
\psbezier(144,9)(144,6)(144,5)(141,5)
\psbezier(141,5)(138,5)(127,5)(123,5)
\psbezier(123,5)(119,5)(79,5)(79,5)
}
\psline[linewidth=0.05]{->}(91,9)(83,9)
\end{pspicture}
\caption{The three successive collections of rectangles $\cE_{0}^0, \cE_{1}^1, \cE_{2}^2$. The graph $\cG_{0}^0$ is reduced to a single rectangle; the graph $\cG_{1}^1$ has two connected components, the graph $\cG_{2}^2$ has three connected components. Note that $E_{2}^2 \cap E_{0}^0 \neq E_{2}^0$.}
\label{f.combinatorics}
\end{figure}
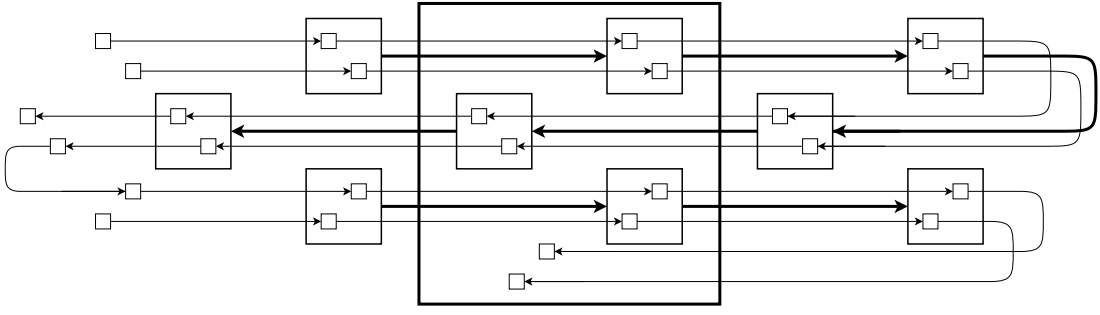

\subsection{Consequences of hypotheses $\mathbf{A_{1,2,3}}$}
\label{ss.consequences-A}
Let $({\cal E}^0_n)_{n\in \NN}$ be a sequence of collections of rectangles. 

\paragraph{Main consequence of hypothesis $\mathbf{A_{1.c}}$ (compatibility).} 
One immediately checks that, under hypothesis $\mathbf{A_{1.c}}$, we have, for every $n\geq 0$,
$$E_{n+1}^{n+1} \cap E_{n}^0 = E^{0}_{n+1}.$$ 
Moreover, using an easy induction,  one can check that, under hypothesis $\mathbf{A_{1.c}}$, 
we have, for every  $n,k\geq 0$
$$E_{n+k}^{n+1} \cap E_{n}^0 = E^{0}_{n+k}.$$ 
Be careful: in general, it is not true\footnote{Actually, if $R$ is minimal, then for every $n\geq 0$ and every rectangle $X\in\cE_n^0$, there exists an integer $k$ such that $E_{n+k}^{n+k} \cap X\neq E^{0}_{n+k}$.}  that $E_{n+k}^{n+p} \cap E_{n}^0=E^{0}_{n+k}$
for $p\geq 2$ (see figure~\ref{f.combinatorics}).

\paragraph{The graphs $\cG_n^p$} 
Under hypothesis $\bm{A_{1.a}}$, for every $n\in \NN$ and every $p\in \{0\dots,n\}$ we will consider the oriented graph $\cG^p_n=\cG(\cE^p_n)$. By hypothesis $\bm{A_{1.a}}$, the graphs $\cG^n_n$ have no cycle. More precisely, each connected component of $\cG_n^n$ is a trivial linear graph 
$$X \ra R(X) \ra \cdots \ra R^l(X)$$
(for some $l \geq 2n+1$).
Hypothesis $\mathbf{A_{1.c}}$ implies that each connected component of $\cG_{n+1}^n$ is isomorphic to a connected component $\cG_n^n$ (the isomorphism
is given by the inclusion of the rectangles). 

By definition of the graphs $\cG_n^p$, the first and the last rectangles of any
connected component of the graph $\cG_{n+1}^{n+1}$ do not belong to the graph $\cG_{n+1}^n$.
Moreover, two consecutive connected components of $\cG_{n+1}^n$ within a connected component of $\cG_{n+1}^{n+1}$ are separated by one or two rectangles (see figure~\ref{f.combinatorics}).

\paragraph{The Cantor set $K$}
We denote by $K$ the decreasing intersection of all the compact sets $E^0_{n}$.
Hypothesis $\mathbf{A_{2}}$ implies that $K$ has no isolated point and
hypothesis $\mathbf{A_{3}}$ implies that $K$ is totally disconnected; hence, if both hypotheses $\mathbf{A_2}$ and $\mathbf{A_3}$ hold, then $K$ is a Cantor set, homeomorphic to the usual triadic Cantor set in $[0,1]$. 

\begin{lemma}
\label{l.KX}
Let $n \geq 0$. Let $X \ra \cdots \ra X'=R^p(X)$ be a path in the graph $\cG_{n}^n$ with $X, X' \in \cE_{n}^0$. Then  $R^p(K \cap X) =K \cap X'$.
\end{lemma}

\begin{proof}
First note that if $|p| > 2n+1$ then the path $X \ra \cdots \ra X'$ must cross $\cE_{n}^0$; by decomposing it into shorter paths we see that it suffices to deal with the case $|p| \leq 2n+1$. Now we have
$$
R^p(E_{n+1}^0 \cap X) = R^p(E_{n+1}^0) \cap X' \subset R^p(E_{n+1}^0) \cap E_{n}^0 \subset E_{n+1}^0
$$
where the last inclusion follows from hypothesis $\mathbf{A_{1.c}}$ (compatibility of $\cE_{n+1}^0$ with $\cE_n^0$ for $n+1$ iterates). Thus 
$$
R^p(E_{n+1}^0 \cap X) \subset E_{n+1}^0 \cap X'.
$$ 
Exchanging the role of $X$ and $X'$, we get the reverse inclusion. Applying this argument recursively, we get that for any positive $m$
$$
R^p(E_{n+m}^0 \cap X)  = E_{n+m}^0 \cap X'.
$$
We now take the intersection on all positive $m$ to get the desired equality.
\end{proof}

\subsection{Construction of the Cantor set $K$: realisation of hypotheses $\mathbf{A_{1,2,3}}$}
\label{ss.construction-cantor}

In this subsection, we first give a characterisation of the Cantor sets $K$ that appear as the intersection of collections of rectangles satisfying hypotheses~$\bm{A_{1,2,3}}$. Then, using Rokhlin lemma,  we construct such a Cantor set inside a given set $A$. 

\begin{defi}
\label{d.tame-Cantor}
A Cantor set (or, more generally, a totally discontinuous set) $K\subset \cM$ is \emph{tamely embedded} if there exist arbitrarily small neighbourhoods of $K$ that are a finite union of pairwise disjoint rectangles of $\cM$.
\end{defi}

The geometry of tamely embedded  Cantor sets is discussed in appendix~\ref{a.cantor}.

\begin{defi}
A Cantor set $K\subset \cM$ is said to be \emph{dynamically coherent} for $R$ if
\begin{itemize}
\item[--] for each integer $k$, the intersection $R^k(K)\cap K$ is open in $K$;
\item[--] for any integer $p\geq 0$ and any point $x\in K$, there exists $p+1$ consecutive positive iterates $R^k(x),\dots,R^{k+p}(x)$ and $p+1$ consecutive negative iterates $R^{-\ell}(x),\dots,R^{-\ell-p}(x)$ outside $K$.
\end{itemize}
\end{defi}

\begin{defi}
\label{d.meagre}
A Cantor set $K\subset\cM$ is said to be \emph{dynamically meagre} for $R$ if $K$ has empty interior in
$\displaystyle\Lambda=\adhe\left(\bigcup_{i\in \ZZ}R^i(K)\right)$.
\end{defi}

\begin{prop}
\label{p.exist-rectangles-1}
A tamely embedded Cantor set $K\subset \cM$ is dynamically coherent if and only if there exists a sequence of collections of rectangles $(\cE^0_n)_{n\in\NN}$ satisfying hypotheses
$\mathbf{A_{1,2,3}}$ and such that  $K=\bigcap_{n \geq 0} E_{n}^0$. 
\end{prop}

\begin{prop}
\label{p.exist-rectangles-2}
Assume that we are given an aperiodic $R$-invariant probability measure $\mu$, and a measurable set $A\subset\cM$ such that $\mu(A)>0$. Then, there exists a Cantor set $K\subset A$ which is dynamically coherent. Furthermore, $K$ can be chosen such that:
\begin{itemize}
\item[--]  it is disjoint from its image $R(K)$;
\item[--]  it is tamely embedded;
\item[--] it is dynamically meagre;
\item[--] $\mu\left(\bigcup_{k\in\ZZ} R^k(K)\right)$ is arbitrarily close to $\mu\left(\bigcup_{k\in\ZZ} R^k(A)\right)$ (in particular, $\mu(K)>0$).
\end{itemize}
\end{prop}

\begin{rema}
\label{r.single-rectangle}
Let $K$ be a Cantor set provided by proposition~\ref{p.exist-rectangles-2}, and $(\cE_n^0)_{n\in\NN}$ be a sequence of collections of rectangles provided by proposition~\ref{p.exist-rectangles-1}.  
\begin{enumerate}
\item One can choose the collection of rectangles $(\cE^0_n)_{n\in\NN}$ such that $E^0_0$ is an arbitrarily small neighbourhood of $K$. Since  $R(K) \cap K = \emptyset$ by proposition~\ref{p.exist-rectangles-2},  we may assume that $R(E_{0}^0)\cap E_{0}^0= \emptyset$, so that  the graph $\cG^0_0$ has no edge.
\item We can always assume that $A$ is included in the support of $\mu$ (by replacing $A$ with $A \cap \supp(\mu)$).
\item We are mainly interested in the case when the measure $\mu$ is supposed to be ergodic (see hypotheses of theorem~\ref{t.main}).
In this case the set $\bigcup_{k\in\ZZ} R^k(K)$ has full measure and its closure is equal to the support of $\mu$. Furthermore, one may assume that the collection $\cE_{0}^0$ is reduced to a single rectangle $X_{0}$. This property is obtained by picking one rectangle $X_{0} \in \cE_{0}^0$ such that the set $K' = K \cap X_{0}$ has positive measure, and replacing the sequence $(\cE_{n}^0)_{n\in\NN}$ by the sequence  $({\cE'}_{n}^0)_{n\in\NN}$ where ${\cE'}_{n}^0$ is the collection of the rectangles of $\cE_{n}^0$ that are included in $X_{0}$. One easily checks that the new sequence $({\cE'}_{n}^0)_{n\in\NN}$ still satisfies hypotheses $\bm{A_{1,2,3}}$.
\end{enumerate}
\end{rema}

\begin{proof}[Proof of proposition~\ref{p.exist-rectangles-1}]
Let us first assume that $K$ is defined as an intersection $K=\bigcap_{n\in\NN} E_n^0$ where $(\cE^0_n)_{n\in\NN}$ is a sequence of collections of rectangles satisfying hypotheses $\mathbf{A_{1,2,3}}$. Clearly, $K$ is tamely  embedded in $\cM$. Fix an integer $n\geq 1$. For every rectangle $X\in\cE_n^0$, let $K_X:=K\cap X$. Then, $\{K_X\}_{X\in\cE_n^0}$ is a partition of $K$ into clopen (closed and open) subsets. Moreover, for any rectangles $X,X'\in \cE^0_n$, if $R^n(X)$ intersects $X'$ then $R^n(X) = X'$ (by hypothesis $\mathbf{A_{1.a}}$), and  the sets $R^n(K_X)$ and $K_{X'}$ coincide (by compatibility, see lemma~\ref{l.KX}). This implies that the intersection $R^n(K)\cap K$ is the union of some of the $K_X$'s, and is thus open in $K$. Finally, by hypothesis $\mathbf{A_{1.a}}$ (no cycle), any point $x\in K$ has a positive iterate $f^k(x)$ and a negative iterate $f^{-\ell}(x)$ which are in $E^{n}_n\setminus E^{n-1}_n$:
the $n-1$ first backward and forward iterates of $f^k(x)$ and $f^{-\ell}(x)$ are disjoint from $E^0_n$ and thus from $K$. So we have proven that $K$ is dynamically coherent for $R$.\\

Conversely we consider a tamely embedded Cantor set $K$ which is dynamically coherent for $R$.
We have to build a sequence of collections of rectangles $(\cE^0_n)_{n\in\NN}$ satisfying hypotheses $\mathbf{A_{1,2,3}}$ such that $K=\bigcap_{n\in\NN} E_n^0$. The construction is made by induction. 

Let us assume that we have already constructed some collections of rectangles $\cE^0_{0},\dots\cE^0_{n-1}$, and that the following induction hypothesis is satisfied for $m\leq n-1$:  
\begin{itemize}
\item[]
\begin{itemize}
\item[$\mathbf{(P_{n})}$]  the boundary $\partial X$ of each rectangle $X\in\cE_m^0$ is disjoint from the set $\bigcup_{i\in\ZZ} R^i(K)$.
\end{itemize}
\end{itemize}
We will now explain how to construct a collection $\cE_n^0$. 

Hypotheses $\mathbf{(P_{1})},\dots,\mathbf{(P_{n-1})}$ imply that, for every $m\leq n-1$, every $X\in\cE^0_m$ and every $k$ such that $|k|<2(n+1)$, the set $K\cap R^{-k}(X)$ is a clopen subset of $K$.
Hence, we can find a partition $\cP_n$ of $K$ into clopen subsets such that, for each $P\in\cP_n$, if $P$ intersects a rectangle $R^{-k}(X)$ with $X\in\cE^0_m$ and $m\leq n-1$, and if $|k|<2(n+1)$, then $P$ is contained in $R^{-k}(X)$. 
Now, we introduce the set $O_{n}$ made of the points $x\in K$ such that  $R^k(x)$ is outside $K$ for $0<k<2(n+1)$. Since $K$ is dynamically coherent, $O_{n}$ is a clopen subset of $K$ and any point of $K$ belongs to or has a positive iterate in $O_{n}$. Similarly, we define the set $I_{n}$ made of the points $x\in K$  such that $R^{-\ell}(x)$ is outside $K$ for $0<\ell<2(n+1)$. Then, $I_n$ is a clopen subset of $K$, and any point in $K$ belongs to or has a negative iterate in $I_{n}$.
(Note that the first return of a point of $O_{n}$ to $K$ occurs in $I_{n}$.) One now defines a new partition $\cD_n$ of $K$ into clopen subsets of $K$ which is finer than $\cP_n$: two points $x$ and $x'$ of $K$ belong to a same element of $\cD_n$ if
\begin{itemize}
\item[--] their first entry time $s\geq 0$ to $O_{n}$
and their last exit time $r\geq 0$ from $I_n$ coincide: in other words,
we have $R^{s}(x)\in O_{n}$, $R^{-r}(x)\in I_{n}$, $R^j(x)\notin I_{n}\cup O_{n}$
for $-r<j<s$ and the same properties hold for $x'$;
\item[--] for any  $-r\leq j\leq s$, the two iterates $f^j(x)$ and $f^j(x')$ are both outside $K$ or  belong to the same element of $\cP_n$.
\end{itemize}
By construction, the collection $\cD_n$ satisfies an equivariance
property: for any $D,D'\in \cD_n$ and any integers $k,\ell$
such that $|k|,|\ell| \leq n+1$, the sets $R^k(D)$ and $R^\ell(D')$
are disjoint or coincide. This shows that the collection $\cD_n$
is organised as the vertices of an oriented graph $\Gamma_n$: one puts
an edge from $D$ to $D'$ if $R^k(D)=D'$ for some positive $k<2(n+1)$
and if $R^j(D)$ is disjoint from $K$ for all the positive $j<k$.
Since any element of $\cD_n$ has a positive iterate included in $O_{n}$
and a negative iterate included in $I_{n}$, the graph $\Gamma_n$
has no cycle and its connected components are linear graphs $D_{0} \subset I_{n} \ra R^{k_{0}}(D_{0}) \ra \cdots \ra R^{k_{p}}(D_{0}) \subset O_{n}$.

The collection $\cE_n^0$ is obtained by ``thickening" the elements of $\cD_n$. Indeed, since $K$ is tamely embedded in $\cM$, every element $D\in \cD_n$ can be thickened as a union $\hat D$ of finitely many disjoint rectangles; this yields a finite collection of rectangles $\cE_n^0$.  

Hypothesis $\mathbf{A_{1.a}}$ is obtained as follows. Since the graph $\Gamma_n$ has no cycle, for each connected component $\Sigma$ of this graph, one can first thicken a single element $D\in \Sigma$ as a set $\hat D$; then any other element $D'\in \Sigma$ can be uniquely written as an image $R^k(D)$ and we set $\hat D'=R^k(\hat D)$. If moreover the connected components of each set $\hat D$ are small enough, then the collection $\cE_n^0$ is $n+1$ times iterable and the associated graph $\cG_n^{n}$ has no cycle (hypothesis $\mathbf{A_{1.a}}$). 

Hypothesis $\mathbf{A_{1.b}}$ is obtained by choosing the elements of $\cE_n^0$ small enough. Indeed, by definition of the partition $\cP_n$, and since  $\cD_n$ is finer than $\cP_n$,  one has the following property: for every $m<n$, every $X\in\cE^0_m$ and $k$ such that $|k|<2(n+1)$, if $D$ is an element of the partition $\cD_n$ with intersects $R^{-k}(X)$, then $D\subset R^{-k}(X)$. Hence, choosing the elements of $\cE_n^0$ small enough, we obtain that:  for every $m<n$, every $X\in\cE^0_m$ and $k$ such that $|k|<2(n+1)$, if $X'$ is an element of the partition $\cE_n^0$ with intersects $R^{-k}(X)$, then $X'\subset R^{-k}(X)$. In particular, $\cE_n^{n+1}$ refines $\cE_m^{m+1}$ for $m< n$ (hypothesis $\mathbf{A_{1.b}}$).  

Hypothesis $\mathbf{A_{1.c}}$ is obtained as follows. First, choosing the elements of $\cE_n^0$ small enough, we get $E^0_{n}\subset E^0_{n-1}$. Second, since $K$ is dynamically coherent, the set $$\bigcup_{|k|\leq 2(n+1)+1}R^k(K)\setminus K$$ is a closed set. Hence, we can choose the set $E_n^0$ in such a way that the following equality holds
$$
E_n^0\cap \bigcup_{|i|\leq 2(n+1)+1}R^i(K)=K\quad\quad\quad (\star)
$$
Let us now consider a rectangle $X\in \cE^0_{n}$ and an iterate $R^{k}(X)$, with \mbox{$|k| \leq 2((n-1))+1$}, that intersects $E^0_{n-1}$. We have to show that $R^k(X)$ also belongs to $\cE^0_n$. We therefore introduce the rectangle $\hat X \in\cE^0_{n-1}$ that contains $X$ and a rectangle $\hat X'\in\cE^0_{n-1}$ that intersects $R^k(X)$. Since the collection $\cE^{n+1}_n$ refines the collection $\cE^{(n-1)+1}_{n-1}$, we have $R^k(X)\subset \hat X'$. We obtain $R^k(X\cap K)\subset \hat X'$ so that by ($\star$) we get  $R^k(X\cap K)\subset K$. Hence $R^k(X)$ intersects a rectangle $X'\in \cE^0_n$ and since $\cE^0_n$ is $n+1$ iterable, we have $X'=R^k(X)$. This proves that $E_n^0$ is compatible with $E_{n-1}^0$ for $n$ iterates (hypothesis $\mathbf{A_{1.c}}$).

Hypothesis $\mathbf{A_2}$ and $\mathbf{A_3}$ are easily obtained by choosing the elements of $\cE_n^0$ small enough. 

Finally, the induction property $\mathbf{(P_n)}$ can be obtained as follows. The space $\cH$ of homeomorphisms $\phi$ from the closed unit ball $\DD$ to a subset of $\cM$, endowed with the uniform convergence topology is a Baire space. Since $K$ is tamely embedded, the homeomorphisms $\phi$ such that $\phi(\partial \DD)$ is disjoint from $K$ is an open and dense subset of $\cH$. Hence, for any homeomorphism $\phi$ in a dense $G_\delta$ subset of $\cH$, the boundary $\phi(\partial \DD)$
is disjoint from all the iterates of $K$. Since each ball $X\in\cE^0_n$ is the image of a homeomorphism $\varphi\in \cH$, one modify $X$ in order to get hypothesis $(P_n)$ by considering a homeomorphism close to $\phi$ in this $G_\delta$ set.
\end{proof}

\begin{proof}[Proof of proposition~\ref{p.exist-rectangles-2}]
We will build a decreasing sequence of compact sets $(K_n)_{n\in\NN}$ contained in $A$ with the following main properties (see figure~\ref{f.tower}).
\begin{enumerate}
\item For each $n$, there exists a compact set $F_n$ such that $K_n=F_n\cup T_n(F_n)\cup\dots\cup T_n^{\ell_n}(F_n)$, where $T_n$ is the return map to $K_n$ for $R$.
\item The sets $F_n, T_n(F_n),\dots,T_n^{\ell_n}(F_n)$ are  pairwise disjoint.
\item For every finite sequence $(t_i)=(t_{1}, \dots , t_{l_{n}})$, define the set  $F_{n}^{(t_i)}$ of those $x \in F_{n}$ whose sequence of successive return times in $K_{n}$ is $(t_i)=(t_{1}, \dots , t_{l_{n}})$ (that is, for every $0 \leq i \leq \ell_{n}-1$, the return time of $T_{n}^i(x)$ in $K_n$ is $t_{i+1}$). Then all but a finite number of those sets $F_{n}^{(t_{i})}$ are empty, and the other ones are clopen subsets of $F_{n}$.
\item The sets $R^{-k}(F_n)$ and $R^{k}(T_n^{\ell_n}(F_n))$, with $1\leq k\leq n$ are disjoint from $K_n$.
\item For $n\leq p$, we have $F_{p} \subset F_{n}$.  Let $F_p(n)$ be the union of all the sets $T_p^k(F_p)$ with $0\leq k\leq \ell_p$ that are contained in $F_n$. Then, $K_p=F_p(n)\cup T_n(F_p(n))\cup\dots\cup T_n^{\ell_n}(F_p(n))$. In particular, we have $F_{p}(n) = K_{p} \cap F_{n}$.
\end{enumerate}
Then we will prove that these properties imply that the set $K:=\bigcap_{n\in\NN} K_n$ is dynamically coherent. Finally, we will explain how to get the other desired properties for the set $K$. 

\begin{figure}[htbp]
\label{f.tower}
\def\JPicScale{0.8}
\hspace{-1.5cm}
\ifx\JPicScale\undefined\def\JPicScale{1}\fi
\psset{unit=\JPicScale mm}
\psset{linewidth=0.3,dotsep=1,hatchwidth=0.3,hatchsep=1.5,shadowsize=1}
\psset{dotsize=0.7 2.5,dotscale=1 1,fillcolor=black}
\psset{arrowsize=1 2,arrowlength=1,arrowinset=0.25,tbarsize=0.7 5,bracketlength=0.15,rbracketlength=0.15}
\begin{pspicture}(0,0)(204,43)
\newrgbcolor{userFillColour}{0.9 0.9 0.9}
\rput{0}(182,23.5){\psellipse[linewidth=0.1,linecolor=white,fillcolor=userFillColour,fillstyle=solid](0,0)(10,10)}
\newrgbcolor{userFillColour}{0.9 0.9 0.9}
\rput{0}(77,23.5){\psellipse[linewidth=0.2,linestyle=dashed,dash=1 1,fillcolor=userFillColour,fillstyle=solid](0,0)(10,10)}
\newrgbcolor{userFillColour}{0.9 0.9 0.9}
\rput{0}(147,23.5){\psellipse[linewidth=0.1,linecolor=white,fillcolor=userFillColour,fillstyle=solid](0,0)(10,10)}
\newrgbcolor{userFillColour}{0.9 0.9 0.9}
\rput{0}(112,23.5){\psellipse[linewidth=0.1,linecolor=white,fillcolor=userFillColour,fillstyle=solid](0,0)(10,10)}
\rput{0}(179.5,28.5){\psellipse[linewidth=0,linecolor=white,fillcolor=lightgray,fillstyle=solid](0,0)(2.5,2.5)}
\rput{0}(179.5,18.5){\psellipse[linewidth=0,linecolor=white,fillcolor=lightgray,fillstyle=solid](0,0)(2.5,2.5)}
\rput{0}(187,23.5){\psellipse[linewidth=0,linecolor=white,fillcolor=lightgray,fillstyle=solid](0,0)(2.5,2.5)}
\rput{0}(144.5,28.5){\psellipse[linewidth=0,linecolor=white,fillcolor=lightgray,fillstyle=solid](0,0)(2.5,2.5)}
\rput{0}(144.5,18.5){\psellipse[linewidth=0,linecolor=white,fillcolor=lightgray,fillstyle=solid](0,0)(2.5,2.5)}
\rput{0}(152,23.5){\psellipse[linewidth=0,linecolor=white,fillcolor=lightgray,fillstyle=solid](0,0)(2.5,2.5)}
\rput{0}(109.5,28.5){\psellipse[linewidth=0,linecolor=white,fillcolor=lightgray,fillstyle=solid](0,0)(2.5,2.5)}
\rput{0}(109.5,18.5){\psellipse[linewidth=0,linecolor=white,fillcolor=lightgray,fillstyle=solid](0,0)(2.5,2.5)}
\rput{0}(117,23.5){\psellipse[linewidth=0,linecolor=white,fillcolor=lightgray,fillstyle=solid](0,0)(2.5,2.5)}
\rput{0}(94.5,3.5){\parametricplot[linewidth=0.2,arrows=<-]{67.38}{112.62}{ t cos 32.5 mul t sin 32.5 mul }}
\rput(94.5,41){$T_{n-1}$}
\rput{0}(129.5,3.5){\parametricplot[linewidth=0.2,arrows=<-]{67.38}{112.62}{ t cos 32.5 mul t sin 32.5 mul }}
\rput(129.5,41){$T_{n-1}$}
\rput{0}(164.5,3.5){\parametricplot[linewidth=0.2,arrows=<-]{67.38}{112.62}{ t cos 32.5 mul t sin 32.5 mul }}
\rput(164.5,41){$T_{n-1}$}
\psline[linewidth=0.2]{->}(78,28.5)(106,28.5)
\psline[linewidth=0.2]{->}(148,18.5)(176,18.5)
\psline[linewidth=0.2]{->}(113,18.5)(141,18.5)
\psline[linewidth=0.2]{->}(156,23.5)(184,23.5)
\psline[linewidth=0.2]{->}(121,23.5)(149,23.5)
\psline[linewidth=0.2]{->}(86,23.5)(114,23.5)
\psline[linewidth=0.2]{->}(148,28.5)(176,28.5)
\psline[linewidth=0.2]{->}(113,28.5)(141,28.5)
\psline[linewidth=0.2]{->}(78,18.5)(106,18.5)
\pscustom[]{\psbezier{-}(183,28.5)(192.8,28.5)(197,25.5)(197,18.5)
\psbezier(197,18.5)(197,11.5)(193.85,8.5)(186.5,8.5)
\psbezier(186.5,8.5)(179.15,8.5)(161.45,8.5)(127.5,8.5)
\psbezier(127.5,8.5)(93.55,8.5)(76.45,8.5)(70.5,8.5)
\psbezier(70.5,8.5)(64.55,8.5)(62,10.75)(62,16)
\psbezier{->}(62,16)(62,21.25)(67.1,23.5)(79,23.5)
}
\pscustom[]{\psbezier{-}(190,23.5)(199.8,23.5)(204,20.5)(204,13.5)
\psbezier(204,13.5)(204,6.5)(200.85,3.5)(193.5,3.5)
\psbezier(193.5,3.5)(186.15,3.5)(166.2,3.5)(127,3.5)
\psbezier(127,3.5)(87.8,3.5)(68.9,3.5)(64,3.5)
\psbezier(64,3.5)(59.1,3.5)(57,5.75)(57,11)
\psbezier{->}(57,11)(57,16.25)(61.2,18.5)(71,18.5)
}
\rput[l](201,23.5){$T_n$}
\rput{0}(74.5,28.5){\psellipse[linewidth=0.2,linestyle=dashed,dash=1 1,fillcolor=lightgray,fillstyle=solid](0,0)(2.5,2.5)}
\rput{0}(74.5,18.5){\psellipse[linewidth=0.2,fillcolor=lightgray,fillstyle=solid](0,0)(2.5,2.5)}
\rput{0}(82,23.5){\psellipse[linewidth=0.2,fillcolor=lightgray,fillstyle=solid](0,0)(2.5,2.5)}
\rput[l](82,21.5){}
\rput[l](24,41){$K_{n-1}$}
\rput[l](24,27){$K_n$}
\rput[l](24,34){$F_{n-1}$}
\rput[l](24,20){$F_n$}
\rput[l](24,13){$F_n(n-1)$}
\rput{0}(19,13){\psellipse[linewidth=0.2,fillcolor=lightgray,fillstyle=solid](0,0)(2,2)}
\rput{0}(19,20){\psellipse[linewidth=0.2,linestyle=dashed,dash=1 1,fillcolor=lightgray,fillstyle=solid](0,0)(2,2)}
\newrgbcolor{userFillColour}{0.9 0.9 0.9}
\rput{0}(19,34){\psellipse[linewidth=0.2,linestyle=dashed,dash=1 1,fillcolor=userFillColour,fillstyle=solid](0,0)(2,2)}
\newrgbcolor{userFillColour}{0.9 0.9 0.9}
\rput{0}(19,41){\psellipse[linewidth=0,linecolor=white,fillcolor=userFillColour,fillstyle=solid](0,0)(2,2)}
\rput{0}(19,27){\psellipse[linewidth=0,linecolor=white,fillcolor=lightgray,fillstyle=solid](0,0)(2,2)}
\rput{0}(12,13){\psellipse[linewidth=0.2,linestyle=dashed,dash=1 1,fillcolor=lightgray,fillstyle=solid](0,0)(2,2)}
\rput(15,11){$,$}
\rput(14,11){}
\end{pspicture}
\caption{Two successive towers in the construction of $K$}
\end{figure}
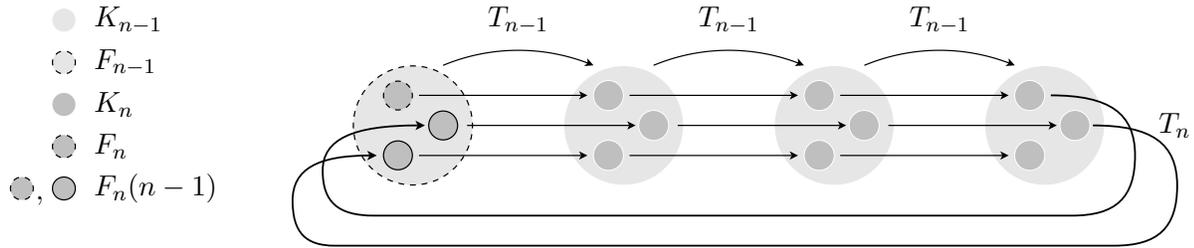

Let us first explain how to build a sequence $(K_n)_{n\in\NN}$ satisfying properties~$1,...,5$. One chooses a compact set $K_0=F_0$ in $A$ with positive $\mu$-measure. Let us assume that $K_{n-1}$ has been defined. On $F_{n-1}$, one considers the dynamics of the return map $R_{F_{n-1}}$ associated to $R$. By Rokhlin lemma, there exists a measurable set $B\subset F_{n-1}$ and an (arbitrarily large) integer $s$ such that the subsets $B,R_{F_{n-1}}(B),\dots,R_{F_{n-1}}^{s}(B)$ of $F_{n-1}$ are pairwise disjoint, and such that the measure $\mu\left(F_{n-1}\setminus \left(B\cup R_{F_{n-1}}(B)\cup\dots\cup R_{F_{n-1}}^{s}(B)\right)\right)$ is arbitrarily small. One sets $F_n=B$ and:
\begin{itemize}
\item[a.] $F_n(n-1)= F_{n} \cup R_{F_{n-1}}(F_{n})\cup\dots\cup R_{F_{n-1}}^{s}( F_{n} )$
\item[b.] $K_n=F_n(n-1)\cup T_{n-1}(F_n(n-1))\cup\dots\cup T_{n-1}^{\ell_{n-1}}(F_n(n-1))$.
\end{itemize}
Denoting by $T_n$ be the first return map of $R$ to $K_n$, we observe that we have
\begin{itemize}
\item[c.] $K_n=F_n\cup T_{n}(F_n)\cup\dots\cup T_{n}^{\ell_{n}}(F_n)$, where $\ell_n=(s+1)(\ell_{n-1}+1)-1$.
\end{itemize}
We define the partition of $F_{n}$ by the sets $F_{n}^{(t_{i})}$ according to the $\ell_{n}$ first return times in $K_{n}$ (as in property~3 above).  Shrinking a little the set $F_{n}$, one can assume that this is a finite partition by compact (and thus clopen) subsets of $F_{n}$ (this makes use of the regularity of the measure $\mu$); in this process we keep properties a, b and c.

One checks inductively that the equality $K_n=F_n(m)\cup T_m(F_n(m))\cup\dots\cup T_m^{\ell_m}(F_n(m))$ also holds for any $m<n-1$. In this construction, one can also consider a first Rokhlin tower $C,R_{F_{n-1}}(C),\dots,R_{F_{n-1}}^{t}(C)$: since the measure of $C\cup R_{F_{n-1}}(C)\cup\dots\cup R_{F_{n-1}}^{n-1}(C)$ and $R_{F_{n-1}}^{t-n+1}(C)\cup\dots\cup R_{F_{n-1}}^{t}(C)$ is arbitrarily small if $t$ is large, one can choose in the previous proof $B=R_{F_{n-1}}^n(C)$ and $s=t-2n$. By this choice, the sets $R^{-k}(F_n)$ and $R^{k}(T_n^{\ell_n}(F_n))$, with $1\leq k\leq n$ are disjoint from $K_n$. Hence, all the required properties on $K_n$ are satisfied.\\

Now, assume that we have such a sequence of sets $(K_n)_{n\in\NN}$ satisfying properties~$1,\dots,5$. Let $K:=\bigcap_{n\in\NN} K_n$. 
Clearly, $K$ is a compact subset of $A$. Let us prove that $K$ is dynamically coherent. 
For every $n\in\NN$, let $F(n) := K \cap F_{n}$. Then, we have  $K=F(n)\cup T_n(F(n))\cup\dots\cup T_n^{\ell_n}(F(n))$, and the sets $R^{-k}(F(n))$ and $R^{k}(T_n^{\ell_n}(F(n)))$ are disjoint from $K$ for every $k$ such that $1\leq k\leq n$ . This implies that the compact set $K$ satisfies the second property of the definition of dynamical coherence. Let us explain how to prove the first one.
The partition of the set $F_{n}$ into clopen subsets $F_n^{(t_i)}$ having the same sequence $(t_i)=(t_1,\dots,t_{l_n})$ of $\ell_{n}$ first return times in $K_{n}$ induces a partition of the set $F(n)$ by clopen sets having the same sequence $(t_i)$ of $\ell_{n}$ first return times in $K$  (the return times of $R$ in $K$ and $K_n$ are the same).
The set $F(n) \cap R^{-n} (K)$ can be written as a union of some elements of this partition, thus it is also a clopen subset\footnote{Here we use that on $T_{n}^{\ell_{n}}(F(n))$ the return time is larger that $n$, thus we do not need to know the return times on that set.}  of $F(n)$. Similarly, the set $T_{n}(F(n))$ admits a finite partition into compact (hence clopen) subsets $R^{t_{1}}(F_{n}^{(t_{i})})$; thus it is a compact set whose subset $T_{n}(F(n)) \cap R^{-n} (K)$ is a clopen subset obtained as a union of some sets of this partition. The same argument works for all the   sets $T_{n}^i(F(n))$ for $i=0, \dots , \ell_{n}$.
Consequently, we see firstly that these sets constitute a (finite) partition of $K$ by clopen subsets, and secondly that  the set $K \cap R^{-n}(K)$ is a finite union of clopen subsets of $K$ and thus it is also a clopen subset of $K$. Hence $K$ is dynamically coherent. \\

Finally, we have to explain how to choose the sequence $(K_n)_{n\in\NN}$ in order to get the other desired properties for the compact set $K$: it is a Cantor set,  it is disjoint from $R(K)$, tamely embedded , dynamically meagre, and the measure of $\bigcup_{k\in \ZZ}R^k(K)$ is arbitrarily close to the measure of $\bigcup_{k\in \ZZ}R^k(A)$.

Clearly, the sequence $(K_n)_{n\in\NN}$ can be chosen such that  the diameter of the connected components of $K_n$ tends to $0$ when $n\to\infty$. This implies that $K=\bigcap_n  K_n$ is totally disconnected. Moreover, at each step, one may choose $K_n$ such that the Hausdorff distance between $K_{n}$ and $K_{n-1}$ is arbitrarily small. One can also assume (the measure $\mu$ has no atom) that each $K_{n}$ has no isolated point. These two facts imply that $K$ has no isolated point. Hence, the sequence $(K_n)_{n\in\NN}$ can be chosen such that $K$ is a Cantor set.

In order to get a set $K$ which is disjoint from $R(K)$, it is enough to choose carefully the set $K_0$: by Rokhlin lemma, there exists a measurable set $A'\subset A$ such that $A'$ and $R(A')$ are disjoint and such that the measures of $\bigcup_{k\in \ZZ}R^k(A)$ and of $\bigcup_{k\in \ZZ}R^k(A')$ are arbitrarily close.
One then chooses the compact set $K_0$ in $A'$, such that the measure $\mu(A'\setminus K_0)$ is arbitrarily small. We get $K\subset A'$ so that $K$ and $R(K)$ are disjoint. 

Let us explain how to get a Cantor set $K$ which is tamely embedded. Let $(T_i)$ be a sequence of triangulations of $\cM$ whose simplexes that are not of maximal dimension have zero $\mu$-measure,
and whose diameters  decrease to $0$ when $i$ tends to infinity. By removing from $K_0$ a set with arbitrarily small measure, one can assume that $K_0$ does not intersect the boundary of the simplexes of the $T_i$. Hence by removing from $K_0$ a set with arbitrarily small measure, we may assume that $K_0$ has a neighbourhood which is a disjoint union of arbitrarily small rectangles\footnote{It may happen that the topological manifold $\cM$ does not admit any triangulation, but this is not a problem for our purpose. Indeed, one can always cover $\cM$ by a finite number of balls $B_1,\dots,B_n$ such that the boundary of these balls has zero measure. Then, one can choose an open neighbourhood $U$ of the union of the boundaries of the $B_i$'s such that $\mu(U)$ is small, and replace $K_0$ by $K_0':=K_0\setminus U$. By construction, the compact set $K_0'$ is included in the interiors of the balls $B_1,\dots,B_n$; so, for $K_0'$, one can use triangulations of the balls $B_1,\dots,B_n$.}. Repeating inductively the argument for $K_1,K_2,\dots$, we obtain a set $K$ which has a basis of neighbourhoods made of disjoint unions of rectangles, i.e. a set $K$ which is tamely embedded.

In order to get a set $K$ which is dynamically meagre, we proceed as follows. Let us consider an integer $n\geq 0$. One can require that all the non-empty open sets of $K_n$ have positive $\mu$-measure: in particular, $K_n$ contains a dense set of recurrent points. By removing a set with arbitrarily small $\mu$-measure, this implies that $K_n$ is contained in the $\frac 1 n$-neighbourhood of $\left(\bigcup_{k\in \ZZ}R^k(K_n)\right)\setminus K_n$, and hence of $\left(\bigcup_{|k|\leq r_n}R^k(K_n)\right)\setminus K_n$ for some integer $r_n\geq 1$. By choosing the sets $K_{n+i}$  (with $i\geq 1$) to be $\frac 1 n$-dense in $K_n$, one gets that $K_n$ is contained in the $\frac 2 n$-neighbourhood of $\left(\bigcup_{|k|\leq r_n}R^k(K)\right)\setminus K$. This implies that $K$ is dynamically meagre.

Finally, observe that, in the construction of the sequence $(K_n)_{n\in\NN}$, the measure $\mu(A\setminus K_0)$ is arbitrarily small, and the measure $\mu(K_n\setminus K_{n+1})$ is also arbitrarily small for every $n$. This implies that the measures of $\bigcup_{k\in \ZZ}R^k(A)$ and of $\bigcup_{k\in \ZZ}R^k(K)$ are arbitrarily close, as required.
\end{proof}

\clearpage

\part{Blowing-up of the orbit of $K$}
\label{part.M_n}

All along part~\ref{part.M_n}, we assume that we are given a sequence of collections of rectangles $(\cE^0_n)_{n\in\NN}$ such that hypotheses $\mathbf{A_1}$ and $\mathbf{A_3}$ are satisfied\footnote{In section~\ref{s.Cantor-sets-in-the-fibres}, we will assume moreover that hypothesis $\mathbf{A_2}$ is satisfied, so that $K$ will be a Cantor set, but we do not need this hypothesis in sections~\ref{s.general-scheme} and~\ref{s.transitivity-minimality}.}. We use the notations $\cE_n^p$, $E_n^p$,  $\cG_n^p$ and $K$ defined in section~\ref{s.rectangles}. We also assume that the graph $\cG_{0}^0$ has no edge (see remark~\ref{r.single-rectangle}).
\begin{itemize}
\item[--] In section~\ref{s.general-scheme}, we will introduce a sequence of homeomorphisms $(M_n)_{n\geq1}$, and some hypotheses on this sequence ($\mathbf{B_1}$, $\mathbf{B_2}$ and $\mathbf{B_3}$). Under these hypotheses, we prove that there exist a continuous onto map  $\Psi = \cdots \circ M_{n} \circ \cdots \circ M_1$ and a homeomorphism $g:\cM\rightarrow\cM$ such that $\Psi\circ g=R\circ \Psi$.
\item[--] In section~\ref{s.transitivity-minimality}, we will formulate another hypothesis, denoted by $\mathbf{B_4}$. If this hypothesis is satisfied, and if $R$ is minimal (resp. transitive), then $g$ is also  minimal (resp. transitive).
\item[--] In section~\ref{s.Cantor-sets-in-the-fibres}, we will state two additional hypotheses ($\mathbf{B_5}$ and $\mathbf{B_6}$) in order to embed a Cantor set $K\times C$ in $\Psi^{-1}(K)$, such that the dynamics of $g$ on the orbit of $K\times C$ is isomorphic to the trivial dynamics $R\times\mbox{Id}$  (for the notion of isomorphism defined in the introduction). This section also contains the construction of a sequence $(M_{n})_{n \geq 1}$ satisfying hypotheses $\mathbf{B_{1,2,4,5,6}}$ (but maybe not $\mathbf{B_{3}}$).
\item[--] Finally, in section~\ref{s.extraction}, we will explain the extraction process, that allows in particular to get hypothesis $\mathbf{B_{3}}$.
\end{itemize}
The purpose of part~\ref{part.H_k} will be to modify the construction in order to get a non-trivial dynamics on the orbit of the Cantor set $K\times C$.

\section{General scheme}
\label{s.general-scheme}

\subsection{The sequences of homeomorphisms $(M_n)_{n \geq 1}$, $(\Psi_n)_{n \in \NN}$ and $(g_n)_{n \in \NN}$}
\label{ss.def-M}

We consider a sequence $(M_n)_{n\geq 1}$ of homeomorphisms of the manifold $\cM$.  Given this sequence, we consider, for every $n \geq 1$, the homeomorphism $\Psi_n$ defined by 
$$
\Psi_n=M_n\circ \dots \circ  M_2\circ M_1,
$$
and the homeomorphism $g_n$ defined by
$$
g_n=\Psi_n^{-1}\circ R \circ \Psi_n.
$$
We also set $\Psi_{0}= \mathrm{Id}$ and $g_{0}= R$.

\subsection{Hypotheses $\mathbf{B_{1,2,3}}$}
\label{ss.hyp-B_123}

We consider the following hypotheses on the homeomorphisms $(M_n)_{n \geq 1}$.
\begin{itemize}
\item[]
\begin{itemize}
\item[$\mathbf{B_1}(n)$] \emph{(Support)} ~\\
 The support of the homeomorphism $M_n$ is contained in the set $E_{n-1}^{n-1}$. 
 \item[$\mathbf{B_2}(n)$] \emph{(Commutation)} ~\\
The maps $M_n$ and $R$ commute along the edges of the graph $\cG_{n-1}^{n-1}$.
\item[$\mathbf{B_3}$] \emph{(Convergence)} ~\\
Let $\cA_{n} = \cE_{n}^{n+1} \setminus  \cE_{n}^{n-1}$. Then the supremum of the diameters of the rectangles $\Psi^{-1}_{n-1}(X))$ with $X \in \cA_{n}$ tends to $0$
when $n$ tends to $+\infty$.
\end{itemize}
\end{itemize}
In $\mathbf{B_{1,2}}(n)$ we assume $n \geq 1$.
The precise meaning of $\mathbf{B_2}$ is: for every $X \in \cE_{n-1}^{n-1}$
such that $R(X) \in \cE_{n-1}^{n-1}$, the equality $M_{n} \circ R = R \circ M_{n}$ holds on $X$.

\subsection{Main consequences of hypotheses $\mathbf{B_{1,2,3}}$}
\label{ss.cons-B_123}

The fundamental properties concerning the convergence of the sequences $(\Psi_{n})$ and $(g_{n})$ are stated in proposition~\ref{p.extension}. This proposition is one  step in the proof of theorem~\ref{t.main}. But it is also an interesting result by itself: in appendix~\ref{s.Denjoy}, we will show how proposition~\ref{p.extension} can be used by to produce various kinds of ``Denjoy counter-examples".

\begin{prop}[Existence of $\Psi$ and $g$]
\label{p.extension}
Assume that hypotheses $\mathbf{B_{1,2,3}}$  are satisfied for every $n\geq1$. Then:
\begin{enumerate}
\item The sequence of homeomorphisms $(\Psi_n)$ converges uniformly towards a continuous map \hbox{$\Psi \colon \cM \to \cM$} (which is not invertible in general).
\item The sequence of homeomorphisms $(g_{n})$ converges uniformly towards a homeomorphism $g$ of $\cM$, and the sequence $(g_{n}^{-1})$ converges uniformly towards $g^{-1}$. 
\item The  homeomorphism $g$ is a topological extension of $R$:  one has $R\circ \Psi=\Psi\circ g$.
\end{enumerate}
\end{prop}

Note that the convergence of $(\Psi_n)$ (item 1 above) only uses hypothesis~$\mathbf{B_{1}}$ (but not hypotheses $\mathbf{B_2}$ and $\mathbf{B_3}$).

\begin{proof}[Proof of proposition~\ref{p.extension}]~
\paragraph{1. Convergence of the sequence $(\Psi_{n})$}
Let $n,p$ be two positive integers. The map $\Psi_{n+p}$ is obtained by post-composing the map $\Psi_{n}$ by the homeomorphism $M_{n+p}\circ\dots\circ M_{n+1}$. Hypotheses $\mathbf{B_1}$ and $\mathbf{A_{1.b}}$ imply that each connected component of the support of this homeomorphism is included in a rectangle of the collection $\cE_m^m$ for some $m \in \{n,\dots ,n+p-1\}$. We will use this property several times thereafter, so we label it ``property $(\star)$".

Property $(\star)$ implies that
the uniform distance from $\Psi_{n}$ to $\Psi_{n+p}$ is smaller than the supremum of the diameters of the rectangles of the collection $\bigcup_{m=n}^{n+p}\cE_{m}^{m}$. On the other hand, the supremum of the diameters of the rectangles of the collection $\cE_{m}^m$ is assumed to tends to $0$ when $m\to 0$ (hypothesis $\mathbf{A_{3}}$). Hence, the sequence of maps $(\Psi_{n})$ is a Cauchy sequence. This proves the first assertion of proposition~\ref{p.extension}.

\paragraph{2. Convergence of the sequence  $(g_{n})$}

\begin{lemma}
\label{l.rectangles-pieges}
Consider a point $x\in \cM$, an integer $n \geq 0$, and a rectangle $X \in \cE_{n}^{n+1}$. Then
$$(\Psi_{n}(x) \in X)\Longleftrightarrow (\Psi_{n+p}(x) \in X\mbox{ for every }p\geq 0)\Longleftrightarrow 
(\Psi(x) \in X).$$
\end{lemma}

\begin{proof}[Proof of lemma~\ref{l.rectangles-pieges}]
Let $p$ be a positive integer. By property $(\star)$,  the map $\Psi_{n+p}$ is obtained by post-composing $\Psi_{n}$ by  the homeomorphism $M_{n+p}\circ\dots\circ M_{n+1}$, and each connected component of the support of this homeomorphism $M_{n+p}\circ\dots\circ M_{n+1}$ is included in a rectangle of the collection $\cE_m^m$ for some $m\in\{n,\dots,n+p-1\}$. In particular, each connected component of the support of the homeomorphism $M_{n+p}\circ\dots\circ M_{n+1}$ is either contained in the rectangle $X$, or disjoint from $X$ (see hypothesis $\mathbf{A_{1.b}}$). Hence, the map $\Psi_{n+p}$ is obtained by post-composing $\Psi_{n}$ by a homeomorphism which preserves the rectangle $X$. 

This shows that $\left(\Psi_{n}(x) \in X\right)\Longleftrightarrow \left(\Psi_{n+p}(x) \in X \mbox{ for every }p\geq 0\right).$ Since $X$ is compact, letting $p\to\infty$, we also get that $\left(\Psi_{n}(x) \in X\right)\Longrightarrow \left(\Psi(x) \in X\right).$ So, the only implication we are left to be proved is $\left(\Psi_{n}(x) \notin X\right)\Longrightarrow \left(\Psi(x) \notin X\right).$

Suppose $\Psi_{n}(x) \not \in X$. Consider the sequence of points $(\Psi_{n+p} (x))_{p\geq 0}$. Because of property $(\star)$, two successive points $\Psi_{n+p} (x),\Psi_{n+p+1} (x)$ in this sequence are either equal, or both belong to some rectangle $X' \in \cE_{n+p}^{n+p}$. Hence, either the sequence $(\Psi_{n+p} (x))_{p\geq 0}$ is constant (in which case the property we want to prove is obvious), or there is a $p_0$ such that $\Psi_{n+p_0}(x)$ belongs to some rectangle $X' \in \cE_{n+p_0}^{n+p_0}$, and then the whole  sequence  of points $(\Psi_{n+p}(x))_{p\geq p_0}$ is trapped in $X'$.  Since $X'$ is compact, this implies that the point $\Psi(x)$ is also in $X'$. And since $X'$ is disjoint from $X$ (hypothesis $\mathbf{A_{1}.c}$), we get that $\Psi(x)\notin X$.
\end{proof}

\begin{lemma}
\label{l.An}
Let $n$ be a positive integer. Consider the set $A_n:=E_n^{n+1}\setminus E_n^{n-1}$. Then:
\begin{enumerate}
\item $\Psi_{n-1}^{-1}(A_n)=\Psi_n^{-1}(A_n)= \Psi_{n+p} ^{-1}(A_n) $ for every positive integer $p$,
\item $g_{n+1}=g_n$ and $g_{n}^{-1}= g_{n+1}^{-1}$ on the set  $\cM\setminus\Psi_n^{-1}(A_n)$. 
\end{enumerate}
\end{lemma}

\begin{proof}[Proof of lemma~\ref{l.An}]
Hypothesis $\mathbf{A_{1.c}}$ implies that  $E_{n}^{n+1} \cap E_{n-1}^{n-1} = E_{n}^{n-1}$. Hence, the set $A_n$ is disjoint from the set $E_{n-1}^{n-1}$ which contains the support of the homeomorphism $M_n^{-1}$ (hypothesis~$\mathbf{B_1}$). Hence, $M_n^{-1}(A_n)=A_n$. Hence, $\Psi_n^{-1}(A_n)=\Psi_{n-1}^{-1} \circ M_n^{-1}(A_n)= \Psi_{n-1}^{-1}(A_n)$. Moreover, lemma~\ref{l.rectangles-pieges} implies that $\Psi_n^{-1}(A_n)= \Psi_{n+p} ^{-1}(A_n)$ for every positive integer $p$. This completes the proof of item~1.

\medskip

Concerning item 2, we only prove the equality $g_{n+1}=g_{n}$; the proof the equality concerning  the inverse maps is completely similar. Recall that 
$$
g_n=\Psi_n^{-1}\circ R\circ\Psi_n\quad\mbox{ and }\quad g_{n+1}=\Psi_n^{-1}\circ M_{n+1}^{-1}\circ R\circ M_{n+1}\circ\Psi_n.
$$
Therefore, proving that $g_{n+1}(x)=g_n(x)$ for every $x\in\cM\setminus\Psi_n^{-1}(A_n)$ amounts to proving that $M_{n+1}^{-1}\circ R\circ M_{n+1}(x)=R(x)$ for every $x\in\cM\setminus A_n$. Let $x$ be a point in $\cM\setminus A_n$. 
\begin{itemize}
\item[--] Either $x \not \in E_{n}^{n+1}$.  Then $R(x)\notin E_n^n$. 
In this case the equality holds because  both the points $x$ and $R(x)$ are outside the support of the homeomorphism $M_{n+1}$ (hypothesis $\mathbf{B_1}$). 
\item[--] Or $x$ is in a rectangle $X\in\cE_{n}^{n-1}$. Then the rectangle $R(X)$ is in $\cE_{n}^{n}$, and thus, both rectangles $X$ and $R(X)$ are vertices of the graph  $\cG^{n}_{n}$. In this case the equality holds because $M_{n+1}$ commutes with $R$ along this graph (hypothesis $\mathbf{B_2}$).
\end{itemize}
\end{proof}

We now turn to the proof of the convergence of the sequence homeomorphisms $(g_{n})$. We fix two positive integers $n,p$. We will check that the uniform distance between the two maps $g_{n}$ and $g_{n+p}$ and between the two maps $g_{n}^{-1}$ and $g_{n+p}^{-1}$ is less than some quantity $\varepsilon_{n}$ which tends to $0$ when $n$ goes to infinity. Applying item 2 of lemma~\ref{l.An}, we get that both equalities $g_{n} = g_{n+p}$ and $g_{n}^{-1} = g_{n+p}^{-1}$ hold outside the set 
$$
F_{n,p} = \bigcup_{k=0}^{p-1} \Psi_{n+k}^{-1} (A_{n+k})
$$
where $A_{n+k}=E_{n+k}^{n+k+1}\setminus E_{n+k}^{n+k-1}$. 
Consequently, we have $d(g_{n},g_{n+p})<\varepsilon_n$  and $d (g_{n}^{-1} , g_{n+p}^{-1}) \leq \varepsilon_{n}$ where $\varepsilon_{n}$ is the supremum of the diameters of the connected components of $F_n$. So, we are left to prove that this supremum tends to  $0$ when $n\to \infty$.
Let $X$ be a  rectangle in $\cE_{n+k}^{n+k+1} \setminus \cE_{n+k}^{n+k-1}$ for some $k\in\{0,\dots,p-1\}$. 
By item 1 of lemma~\ref{l.An}, we have $\Psi_{n+k}^{-1} ( X ) =  \Psi_{n+p}^{-1} ( X ) $. Hence, we can rewrite $F_{n,p}$ as
$$
F_{n,p} =  \Psi_{n+p}^{-1} \left( \bigcup_{k=0}^{p-1}  A_{n+k} \right).
$$
Hence, hypothesis $\mathbf{A_{1.b}}$ implies that every connected component of $F_{n,p}$ is a set $\Psi_{n+p}^{-1} (X)$ with $X \in \cE_{n+k}^{n+k+1} \setminus \cE_{n+k}^{n+k-1}$ and $0\leq k\leq p-1$. Using once again item 1 of lemma~\ref{l.An}, we get  that every connected component of $F_{n,p}$ is  some $\Psi_{n+k-1}^{-1} (X)$ with $X \in \cE_{n+k}^{n+k+1} \setminus \cE_{n+k}^{n+k-1}$ and $0\leq k\leq p-1$. Hence hypothesis $\mathbf{B_{3}}$ implies that the supremum of the diameters of the connected components of $F_{n,p}$ tends to $0$ when $n\to\infty$. 

Thus we have proved that the sequences $(g_{n})$ and $(g_{n}^{-1})$ are Cauchy sequences. Hence, they converge respectively towards a homeomorphism $g$ and its inverse $g^{-1}$. 

\paragraph{3. Relation between  $g$ and $R$.}
The semi-conjugacy $R\circ \Psi=\Psi\circ g$ is obtained by taking the limit in the equality $R\circ \Psi_n=\Psi_n\circ g_n$.
\end{proof}

\subsection{Some other consequences of hypothesis $\mathbf{B_{1,2,3}}$}

Let us call \emph{orbit of $K$} the set $\displaystyle\bigcup_{i\in\ZZ} R^i(K)$.

\begin{prop}
\label{p.extension2}
Assume that hypotheses $\mathbf{B_{1,2,3}}$ are satisfied for every $n \geq 1$. Then:
\begin{enumerate}
\item For every point $x \in R^{n_{0}}(K)$ for some integer $n_{0}$, let $(X_{n})_{n \geq \mid n_{0} \mid}$, $X_{n} \in \cE_{n}^{n_{0}}$ be the decreasing sequence of rectangles containing $x$. Then
$$
\Psi^{-1}(x) = \bigcap_{n \geq n_{0}} \Psi_{n}^{-1}(X_{n}).
$$
\item  For every point $y\in\cM$ which does not belong to the orbit of $K$, the fibre $\Psi^{-1}(y)$ is a single point.
\end{enumerate}
\end{prop}

\begin{proof}[Proof of proposition~\ref{p.extension2}]
The first assertion is a direct consequence of lemma~\ref{l.rectangles-pieges}. The proof of the second assertion begins by a lemma. 

\begin{lemma}
\label{l.piege2}
Let $x$ be a point of $M$. The following properties are equivalent:
\begin{itemize}
\item[1.]  the point $\Psi(x)$ belongs to the orbit of $K$;
\item[2.]  for every $n$ large enough, the point $\Psi_{n}(x)$ belongs to the set $E_{n}^n$.
\end{itemize}
\end{lemma}

\begin{proof}[Proof of lemma~\ref{l.piege2}]
Recall that 
$$
K = \bigcap_{n \geq 0} E_{n}^0\quad\mbox{ and }\quad E_{n}^n = \bigcup_{ -n \leq i \leq n} R^i(E_{n}^0).
$$
Implication $1\Longrightarrow 2$ is a consequence of these equalities and lemma~\ref{l.rectangles-pieges}. 

In order to prove implication $2\Longrightarrow 1$, we consider a point $x$ and an integer $n_0$ such that $\Psi_{n}(x)\in E_{n}^n$ for every $n\geq n_0$. In particular, the point $\Psi_{n_0}(x)$ is in the set $E_{n_0}^{n_0}$. Hence, lemma~\ref{l.rectangles-pieges} implies that the point $\Psi_{n_0+1}(x)$ is also in the set $E_{n_0}^{n_0}$. Hence, the point $ \Psi_{n_0+1}(x)$ is in the set 
$$
E_{n_0}^{n_0} \cap E_{n_0+1}^{n_0+1} = E_{n_0+1}^{n_0}
$$
(the equality follows from hypothesis $\mathbf{A_{1.c}}$). By induction, we obtain that the point $\Psi_{n_0+p}(x)$ is in the set $E_{n_0+p}^{n_0}$ for every $p\geq 0$. Since $(E_{n_0+p}^{n_0})_{p\geq 0}$ is a decreasing sequence, this implies that the point $\Psi(x)=\displaystyle\mathop{\lim}_{p\to\infty}\Psi_{n_0+p}(x)$ is in the set 
$\bigcap_{p\geq 0}E_{n_0+p}^{n_0}=\bigcup_{|i| \leq n_0} R^i(K)$.
In particular, the point $\Psi(x)$ belongs to the orbit of $K$.
\end{proof}

Let us turn to the proof of the second assertion of proposition~\ref{p.extension2}.
Pick two points $x,x'\in\cM$ such that $\Psi(x) = \Psi(x') = y$ and such that $y$ is not in the orbit of $K$. We have to prove that $x=x'$. On the one hand, lemma~\ref{l.rectangles-pieges} implies that, for every $n\geq 0$, 
\begin{itemize}
\item[--] either $\Psi_{n}(x) \not \in E_{n}^n$, and then   $\Psi_{n}(x') \not \in E_{n}^n$,
\item[--] or $\Psi_{n}(x) $ is in some rectangle $X \in  \cE_{n}^n$, and then $\Psi_{n}(x')$ is in the same rectangle $X$.
\end{itemize}
On the other hand, lemma~\ref{l.piege2} shows that there exists arbitrarily large integers $n$ such that $\Psi_{n}(x) \not \in E_{n}^n$. The easy case is when this happens for \emph{every} $n$ large enough, say every $n$ bigger than some $n_{0}$: in this case, one has $\Psi(x) = \Psi_{n_{0}}(x)$ and $\Psi(x') = \Psi_{n_{0}} (x')$ (because of hypothesis $\mathbf{B_1}$ on the support of the $M_n$'s), so that $\Psi_{n_0}(x)=\Psi_{n_0}(x')$, and $x=x'$ since the map $\Psi_{n_{0}}$ is invertible. 
If we are not in the easy case, then there exists an increasing sequence of integers $(n_k)_{k\in\NN}$ such that, for every $k$,
\begin{enumerate}
\item $\Psi_{n_k-1}(x), \Psi_{n_k-1}(x')   \not \in E_{n_k-1}^{n_k-1}$
\item $\Psi_{n_k}(x), \Psi_{n_k}(x')  \in X_{k}$ where $X_{k}$ is some rectangle in  $\cE_{n_k}^{n_k}$.
\end{enumerate}
From item 1 and lemma~\ref{l.rectangles-pieges}, we see that, for every $k$, the rectangle $X_{k}$ is not included in (and thus is disjoint from) the set $E_{n_k-1}^{n_k-1}$. Now remember that $\Psi_{n_k} = M_{n_k} \circ \Psi_{n_k-1}$ and that the support of $M_{n_k}$ is included in the set $E_{n_k-1}^{n_k-1}$. From this we get $\Psi_{n_k}^{-1}(X_{k}) = \Psi_{n_k-1}^{-1}(X_{k})$ for every $k$. Item  2 says that the points $x$ and $x'$ both belong to the set $\Psi_{n_k}^{-1}(X_{k})= \Psi_{n_k-1}^{-1}(X_{k})$. We now apply hypothesis $\mathbf{B_{3}}$, which says that the diameter of the set $\Psi_{n_k-1}^{-1}(X_{n_k})$ tends to $0$ when $k\to\infty$. This proves that $x=x'$, as wanted.
\end{proof}

\begin{rema}[Comparison between the maps $g_n$ and $g$]
\label{r.gngnp1}~
\begin{enumerate}
\item Item 2 of  lemma~\ref{l.An} is not optimal: the equality $g_{n}=g_{n+1}$ also holds on every rectangle $\Psi_{n}^{-1}(X)$  such that both $X$ and $R(X)$ belong to $\cE_{n}^n$.
\item For every $n\geq 1$,  we have $g=g_{n}$ on the set   $$\Psi^{-1}\left(R^{-n}(K) \cup \cdots \cup R^{n-1}(K)\right).$$
\end{enumerate}
\end{rema}

\begin{proof}
The first claim follows easily from the commutation of $M_{n+1}$ and $R$ along the edges of the graph $\cG_n^n$ (hypothesis $\mathbf{B_{2}}(n+1)$).

In order to prove the second claim, we consider an integer $n \in \NN$ and a point $\wt x\in\cM$ such that $x:=\Psi(\wt x)$ belongs to $R^{-n}(K) \cup \cdots \cup R^{n-1}(K)$. Then for any $p \geq 0$, both $x$ and $R(x)$ belongs to $E_{n+p}^{n+p}$. We apply the first claim: if $x \in X \in \cE_{n+p}^{n+p}$ then $g_{n+p}= g_{n+p+1}$ on $\Psi_{n+p}^{-1}(X)$. The first assertion of proposition~\ref{p.extension2} says that the point $\wt x$ belongs to this last set. Thus $g_{n}(\wt x) = g_{n+p} (\wt x)$ for any positive $p$. The second claim follows.
\end{proof}

\section{Transitivity, minimality}
\label{s.transitivity-minimality}

In this section, we assume that we are given a sequence $(M_n)_{n\geq 1}$ of homeomorphisms of~$\cM$. We use the notations $\Psi_n$ and $g_n$ defined in section~\ref{s.general-scheme}. We recall that the sequence $(\cE_{n}^0)$ is supposed to satisfy hypotheses $\mathbf{A_{1,3}}$. 
When hypotheses $\mathbf{B_{1,2,3}}$ are satisfied, proposition~\ref{p.extension} provides us with a homeomorphism $g=\lim g_n$ and a continuous map $\Psi=\lim\Psi_n$ such that $\Psi\circ g=R\circ\Psi$. The purpose of the section is to present  hypothesis $\mathbf{B_{4}}$ on the $M_n$'s: when $R$ is supposed to be transitive (or minimal), this hypothesis implies that $g$ is transitive (or  minimal). The more general situation when the dynamics of $R$ is only supposed to be transitive or minimal on a subset of $\cM$ (see addendum~\ref{a.main}) will be treated in appendix~\ref{a.recurrence}.

\subsection{Hypotheses $\mathbf{B_4}$}
\label{ss.hyp-B_4}

We call  \emph{internal radius} of a set $F\subset\cM$ the supremum of the radii of the balls included in~$F$. We consider the following hypothesis.
\begin{itemize}
\item[]
\begin{itemize}
\item[$\mathbf{B_4}$] \emph{(Fibres are thin)}\\
The internal radius of the set $\Psi_n^{-1}(E_n^0)$ goes to $0$ when $n\to\infty$. 
\end{itemize}
\end{itemize}

\subsection{Consequences of hypothesis $\mathbf{B_4}$}
\label{ss.cons-B_4}
We now assume that the sequence $(M_n)_{n\geq 1}$ satisfies hypotheses $\mathbf{B_{1,2,3}}$.

\begin{prop}
\label{p.reformulation}~
 Hypothesis $\mathbf{B_4}$ is satisfied if and only if $\Psi^{-1}(K)$ has empty interior.
\end{prop}

\begin{proof}
Hypothesis $\mathbf{B_{1}}$ on the support easily implies that for any $n$ one has $\Psi^{-1}(E_{n}^0) = \Psi_{n}^{-1}(E_{n}^0)$ (see for example lemma~\ref{l.rectangles-pieges}). Then the  definition  of $K$ gives
$$
\Psi^{-1}(K) = \bigcap_{n \in \NN} \Psi^{-1}(E_{n}^0) = \bigcap_{n \in \NN} \Psi_{n}^{-1}(E_{n}^0)
$$
where  the intersections are decreasing. The proposition follows from these equalities and an easy compactness argument.
\end{proof}

\begin{prop}
\label{p.recurrence}
If the dynamics $R$ on $\cM$ is transitive (resp. minimal) and hypothesis ${\mathbf{B_4}}$ is satisfied, then the dynamics $g$ on $\cM$ is also transitive (resp. minimal).
\end{prop}

\begin{lemma}
\label{l.recurrence}
If hypothesis ${\mathbf{B_4}}$ is satisfied and $F$ is a compact  $g$-invariant subset  of  $\cM$ such that $\Psi(F)=\cM$, then $F= \cM$.
\end{lemma}

\begin{proof}
Under the hypotheses of the lemma, assume $F \neq \cM$.  The second assertion of proposition~\ref{p.extension2} implies that the complement $O$ of $F$ in $\cM$ is an open set contained in $\bigcup_{i\in \ZZ}g^i(\Psi^{-1}(K))$. By Baire theorem, there exists $i\in\ZZ$ such that the set $g^i(\Psi^{-1}(K))$ has non-empty interior in $\cM$. Hence $\Psi^{-1}(K)$ has non-empty interior in $\cM$.  Hence proposition~\ref{p.reformulation} implies that $\mathbf{B_4}$ is not satisfied.
\end{proof}

\begin{proof}[Proof of proposition~\ref{p.recurrence}]
Assume that hypothesis ${\mathbf{B_4}}$ is satisfied, and that the dynamics of $R$ on $\cM$ is minimal. Consider any non-empty $g$-invariant compact set $F\subset \cM$. Then  $\Psi(F)$ is an $R$-invariant compact set, so  $\Psi(F)=\cM$. Lemma \ref{l.recurrence} implies $F=\cM$. Hence the dynamics of $g$ on $\cM$ is minimal. 

Now, assume that hypothesis ${\mathbf{B_4}}$ is satisfied, and that the dynamics of $R$ on $\cM$ is transitive. Consider a point $x\in \cM$  whose $R$-orbit is dense in $\cM$ and choose any lift $\tilde x\in \Psi^{-1}(x)$ of $x$. The   closure of the $g$-orbit of $\tilde x$ projects by $\Psi$ on the closure of the $R$-orbit of $x$, that is on $\cM$. Hence, by lemma \ref{l.recurrence}, the closure of the $g$-orbit of $\tilde x$ is $\cM$. Hence, the dynamics of $g$ on $\cM$ is transitive. 
\end{proof}

\begin{rema*}
Conversely, if the dynamics of $g$ on $\cM$ is minimal, then the dynamics of $R$ on $\cM$ is minimal and hypothesis ${\mathbf{B_4}}$ is satisfied. We will not use this fact,
and the (easy) proof is left to the reader.
\end{rema*}

\section{Cantor sets in the fibres of $\Psi$}
\label{s.Cantor-sets-in-the-fibres}

We recall that we are considering a sequence of collection of rectangles $(\cE_{n}^0)$ such that hypotheses $\mathbf{A_{1,3}}$ are satisfied. Moreover, we will now assume that hypothesis $\mathbf{A_2}$ is also satisfied (so that $K$ is a Cantor set). 

In this section, we will again consider a sequence $(M_n)_{n\geq 1}$ of homeomorphisms of~$\cM$.
We will state two additional hypotheses $\mathbf{B_{5,6}}$  concerning the homeomorphisms $M_n$. Roughly speaking, we want to  ensure that, for any $x$ in $K$, the fibre $\Psi^{-1}(x)$ will contain a Cantor set identified to $\{x\} \times C$ (hypothesis $\mathbf{B_5}$), and that the homeomorphism $g$ will induce a ``trivial" dynamics on $K\times C=\bigcup_{x\in K}\{x\}\times C$ (hypothesis $\mathbf{B_6}$). Then we will explain how to construct a sequence of homeomorphisms $(M_n)_{n\geq 1}$ satisfying hypotheses $\mathbf{B_{1,2,4,5,6}}$. 

\paragraph{Additional assumption.} For sake of simplicity, we will assume that the collection $\cE_0^0$
is made of a single rectangle $X_0$ (see remark~\ref{r.single-rectangle}).

\subsection{The Cantor set $K\times C$}
\label{ss.objects-cantor}

We set $\t X_{0}= X_{0}$.  Even if both sets are equal, this notation trick allows us to deal with $\t X_{0}$ for objects that have a dynamical meaning for the map $g$, and to deal with $X_{0}$ for those which have a dynamical meaning for $R$.  In particular, the semi-conjugacy  $\Psi$ will be constructed so that $\Psi(\t X_{0}) = X_{0}$.

We consider an (abstract) Cantor set $C$. Thus $K \times C$ is again a Cantor set. 
We embed this Cantor set $K\times C$ in the interior of the rectangle $\t X_0$; from now on, we will  see $K\times C$ as a Cantor set in $\t X_0$. For technical reasons  we choose a \emph{tame embedding} of $K\times C$ in $\t X_0$ (see definition~\ref{d.tame-Cantor}).

For every $n\in\NN$ and every rectangle $X\in\cE_n^0$,  we denote by $K_X$ the Cantor set $K\cap X$. Then $K_{X} \times C$  is a  sub-Cantor set of $K\times C$ embedded in the rectangle $\t X_0$. 

\begin{rema*} 
We will often consider some path $X \ra \cdots \ra X'=R^p(X)$ in the graph $\cG_{n}^n$. Remember that this means that the rectangles $X, R(X), \dots , R^p(X)=X'$ are in $\cE_n^n$, and observe that the integer $p$ is unique (since the graph has no cycle).  Moreover, we will always consider the case where the rectangles $X,X'$ are in $\cE_n^0$; in particular $X$ and $X'$ are included in the rectangle $X_0$. 
\end{rema*}

\subsection{Hypotheses $\mathbf{B_{5,6}}$}
\label{ss.hyp-B_56}

\begin{figure}[htbp]
\def\JPicScale{1.3} \hspace{-2cm}
\input{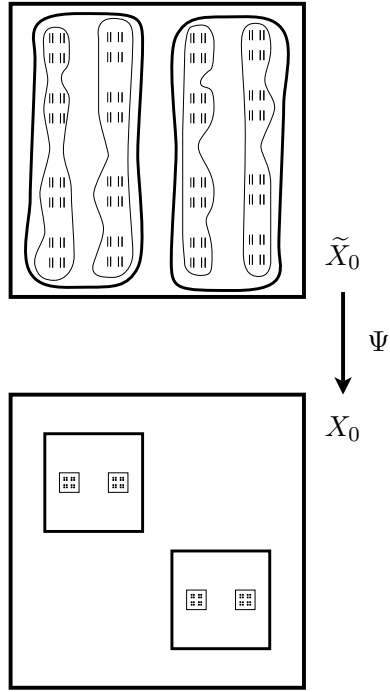}
\caption{Cantor sets in the fibres (hypothesis~$\mathbf{B_5}$).
Below, the Cantor set $K$ inside the rectangles $X$ of $\cE_{n}^0$; above, the Cantor set $K \times C$ inside the rectangle $\Psi^{-1}(X)$.}
\label{f.cantor-fibres}
\end{figure}

We consider the following hypotheses.
\begin{itemize}
\item[]
\begin{itemize}
\item[$\mathbf{B_5}(n)$] \emph{(Cantor sets in the fibres of $\Psi$, see figure~\ref{f.cantor-fibres})} ~\\
For every rectangle $X\in \cE_{n}^0$,  the open set $\Psi_n^{-1}(\inte(X))$ contains the Cantor set $K_{X} \times C.$
\item[$\mathbf{B_6}(n)$] \emph{(Embedding a trivial dynamics)} ~\\
If $X \ra \cdots \ra X'=R^p(X)$ is a path in the graph $\cG_{n}^n$ with $X,X'\in \cE_n^0$, then
$$
g_{n}^p(x,c)= ( R^p(x)  ,c)\quad\mbox{ for every }\quad(x,c) \in K_{X}  \times C.
$$
\end{itemize}
\end{itemize}

Note that in hypothesis $\mathbf{B_6}$ the point $R^p(x)$ belongs to $K_{X'}$
(this is due to the compatibility hypothesis $\mathbf{A_{1.c}}$, see lemma~\ref{l.KX}),
so that $(R^p(x),c)$ is a point of $K_{X'}\times C \subset \cM$, and  the statement is meaningful.

\subsection{Consequences of hypotheses $\mathbf{B_{5,6}}$}
\label{ss.cons-B_56}

Proposition~\ref{p.conjugacy} below states some consequences of hypotheses  $\mathbf{B_{6}}$ on the dynamics of~$g$. We will not use formally this proposition;  we only state it for ``pedagogical reasons": in part~C, there will be a similar proposition, where the map $R \times \mbox{Id}$ will be replaced by some non-trivial fibered dynamics $h$, and whose proof will be entirely similar to the proof of proposition~\ref{p.conjugacy}. We use  the notion of isomorphism defined in the introduction (subsection~\ref{ss.isomorphic}).
\begin{prop}
\label{p.conjugacy}
Assume that hypotheses $\mathbf{B_{1,2,3}}$ are satisfied. Then the dynamics of $g$ on the set $\displaystyle\bigcup_{n\in \ZZ} g^n(K\times C)$ is isomorphic to the dynamics of $R\times\mbox{Id}$ on the set $\displaystyle \bigcup_{n\in \ZZ} R^n(K)\times C$.
\end{prop}

\begin{lemma}
\label{l.conjugacy}
Assume hypotheses $\mathbf{B_{1,2,3}}$.
 Under hypothesis $\mathbf{B_5}$, for every $x \in K$,
$$
\{x\} \times C \subset \Psi^{-1}(x),
$$
that is, for every $(x,c)\in K\times C$, the formula $\Psi(x,c) = x$ holds.
\end{lemma}

\begin{proof}
This is an immediate consequence of the equality $\displaystyle\Psi^{-1}(x) = \bigcap_{n \in\NN} \Psi_{n}^{-1}(X_{n})$ of proposition~\ref{p.extension2}.
\end{proof}

\begin{proof}[Proof of proposition~\ref{p.conjugacy}]
We first note that the return times of $g$ on $K \times C$ and of $R$ on $K$ are coherent: for every $(x,c)$ in $K\times C$ and every integer $p$,
$$
g^p(x,c) \in K \times C \Leftrightarrow R^p(x) \in K.
$$
Indeed, the direct implication follows easily from lemma~\ref{l.conjugacy} and the semi-conjugacy $\Psi g = R \Psi$. The other implication follows not so easily from hypothesis $\mathbf{B_{6}}$.
More precisely, consider any integer $n \geq \mid p\mid$ and the rectangle $X$ such that $x \in X \in \cE_{n}^0$. Since $R^p(X)$ contains $R^p(x)$ which is supposed to belong to $K$, the rectangle $R^p(X)$ also belongs to $\cE_{n}^0$ (as a consequence of the iterability, hypothesis $\mathbf{A_{1.a}}$). Then we can apply hypothesis $\mathbf{B_{6}}$ to the path $X \ra \cdots \ra R^p(X)$ in $\cG_{n}^n$, which shows in particular that  $g_{n}^p(x,c) \in K \times C$. It remains to note that applying
$p$ times remark~\ref{r.gngnp1}  yields $g_{n}^p(x,c)=g^p(x,c)$.

Now consider the map $\Theta$ from the embedded Cantor set $K \times C \subset \t X_{0}$ to
the  set $K \times C \subset X_{0} \times C$ which is just the inverse map of our embedding $K \times C \hookrightarrow \t X_{0}$. Let $(x,c)$ be a point  of $K \times C \subset \t X_{0}$ which returns in $K \times C$ under iteration of $g$, and let $p \geq 1$ be the return time. Then $p$ is also the return time of $\Psi(x,c) = x$ in $K$ under iteration of $R$. Let $n \geq p$. Then hypothesis $\mathbf{B_6}$ implies that 
$$\Theta (g_{n}^p(x,c)) = (R \times \mbox{Id})^p(\Theta(x,c)).$$ 
Since $g_{n}^p(x,c) = g^p(x,c)$ (see remark~\ref{r.gngnp1}), we get that $\Theta$ realises an isomorphism between the first return map of $g$ on $K \times C \subset \cM$ and the first return map of $R \times \mbox{Id}$ on $K \times C \subset \cM \times C$, with compatible return times. Hence $\Theta$ extends (in a unique way) to an isomorphism as required by the proposition.
\end{proof}

\subsection{Construction of the $M_{n}$'s: realisation of hypotheses $\mathbf{B_{1,2,4,5,6}}$}
\label{ss.realis-B_12456}

 The following proposition is one of the main steps of the proof of our main theorem.

\begin{prop}[Existence of $(M_{n})$]
\label{p.existence-M}
Assume that hypotheses $\mathbf{A_{1,2,3}}$ hold and that the graph $\cG_{0}^0$ has no edge.
Then there exists a sequence $(M_{n})_{n\geq 1}$ of homeomorphisms of $\cM$ such that hypotheses $\mathbf{B_{1,2,4,5,6}}$ are satisfied.   
\end{prop}

\begin{rema*}
Note that proposition~\ref{p.existence-M} does not ensure hypothesis $\mathbf{B_3}$. Indeed, it is not possible to get this hypothesis for \emph{any} sequence of collections of rectangles $(\cE_n^0)_{n\in\NN}$. In section~\ref{s.extraction}, we will explain how to obtain hypothesis~$\mathbf{B_3}$ using an extraction process.
\end{rema*}

\begin{proof}
We will proceed by induction.
 For these purpose, we need to state a quantitative version of $\mathbf{B_4}$:
\begin{itemize}
\item[]
\begin{itemize}
\item[$\mathbf{B_4}(n)$] For every $X$ in $\cE_n^0$, the internal radius of  the rectangle $\Psi_{n}^{-1}(X)$ is less than $\frac{1}{n}$.
\end{itemize}
\end{itemize}
Of course, one may replace the sequence $\left(\frac{1}{n}\right)_{n\geq 1}$ by any sequence of positive numbers $(\epsilon_n)_{n\geq 1}$ such that $\epsilon_n\to 0$ when $n\to 0$. Clearly, if hypothesis $\mathbf{B_4}(n)$ is satisfied for every $n$, then hypothesis $\mathbf{B_4}$ is satisfied. 

Now let $n \geq 1$. We assume that some homeomorphisms $M_1,\dots,M_{n-1}$ have been constructed and that hypotheses $\mathbf{B_{1,2,4,5,6}}(m)$ are satisfied for every $m\in\{1,\dots,n-1\}$.
We will explain how to construct a homeomorphism $M_n$ such that hypotheses $\mathbf{B_{1,2,4,5,6}}(n)$ are satisfied.

\paragraph{Step 1. Choice of a tame embedding   $\bm{i_X}$  of $\bm{K_X\times C}$ in $\bm{\inte(X)}$.}
We first need to choose such an embedding for every rectangle $X\in\cE_n^0$. Furthermore, we need these embeddings to satisfy the following equivariance property: if $X \ra \cdots \ra X'=R^p(X)$ is a path in the graph $\cG_{n}^n$ with $X,X' \in \cE_{n}^0$, then   
\begin{equation}
\label{e.equivariance-i_X}
R^p\circ i_X (x,c)= i_{X'}(R^p(x),c)
\end{equation}
for every $(x,c)\in K_X\times C$. Remember that the point $(R^p(x),c)$ is known to belong to $K_{X'} \times C$ (lemma~\ref{l.KX}), so this equality is compatible with the the requirement that ${i_{X'}}(K_{X'}\times C) \subset {\inte(X')}$. The construction of such a family of homeomorphisms is straightforward, since the graph $\cG_{n}^n$ has no cycle (hypothesis $\mathbf{A_{1.a}}$).

\bigskip

Recall that the support of $M_n$ has to be included $E_{n-1}^{n-1}$ (hypothesis $\mathbf{B_1}(n)$); so we need to define $M_n$ on each rectangle of $\cE_{n-1}^{n-1}$, i.e. on each vertex of the graph  $\cG_{n-1}^{n-1}$. From now on we will treat the different connected components of the graph $\cG_{n-1}^{n-1}$  independently\footnote{However note that the homeomorphisms $(i_X)_{X\in\cE_n^n}$ defined at step~1 are equivariant not only along the graph $\cG_n^{n-1}$, but also along the graph $\cG_n^n$. This will automatically establish a link between the restrictions of $M_{n}$ to any two rectangles belonging to the same connected component of $\cG_{n}^{n}$.}. Moreover, since $M_n$ has to commute with $R$ along the edges of the graph $\cG_{n-1}^{n-1}$ (hypothesis $\mathbf{B_2}(n)$), we will first define $M_n$ on one vertex of each connected component of the graph  $\cG_{n-1}^{n-1}$ (step~2) and then extend $M_{n}$ along the component by commutation (step~3). Also note that for $n=1$ the graph $\cG_{n-1}^{n-1}$ has no edge (by assumption), so step 3 is useless.

\paragraph{Step 2. Definition of $\bm{M_n}$ on one vertex of each connected component of the graph $\bm{\cG_{n-1}^{n-1}}$.} 
  For each connected component $\Gamma$ of the graph $\cG_{n-1}^{n-1}$, we choose one vertex $\widehat X$ of $\Gamma$ which is in $\cE_{n-1}^0$. We define $M_n$ on this rectangle $\widehat X$ as follows (see figure~\ref{f.Mn}).   

We denote by $X_1,\dots,X_\ell$ the immediate sub-rectangles of $\widehat X$, \emph{i.e.} the rectangles of $\cE_n^0$ contained in $\widehat X$. For first we observe that the set $\Psi_{n-1}^{-1}(\widehat X)$ contains the  Cantor set $K_{\widehat X}\times C$ (by the induction\footnote{When $n = 1$ the observation rather follows from the convention $\Psi_{0}=\mathrm{Id}$.} hypothesis $\mathbf{B_{5}}(n-1)$),  and that the Cantor set $K_{\widehat X}\times C$ is the  disjoint union of the Cantor sets $K_{X_1}\times C,\dots,K_{X_\ell}\times C$. As a consequence, the Cantor set $\Psi_{n-1}(K_{\widehat X}\times C)$ in $\inte(\widehat X)$ is the disjoint union of the Cantor sets  $\Psi_{n-1}(K_{X_1}\times C),\dots,\Psi_{n-1}(K_{X_\ell}\times C)$.  Thus we can construct  the homeomorphism $M_{n}$ on  $\Psi_{n-1}(K_{\widehat X}\times C)$ the following way: on each subset $\Psi_{n-1}(K_{X_i}\times C)$ we define $M_{n}$  by the formula
\begin{equation}
\label{e.def-M_n}
M_{n}(\Psi_{n-1}(x,c))= i_{X_i}(x,c).
\end{equation}

\begin{figure}[htbp]
\label{f.Mn}
\begin{center}
\input{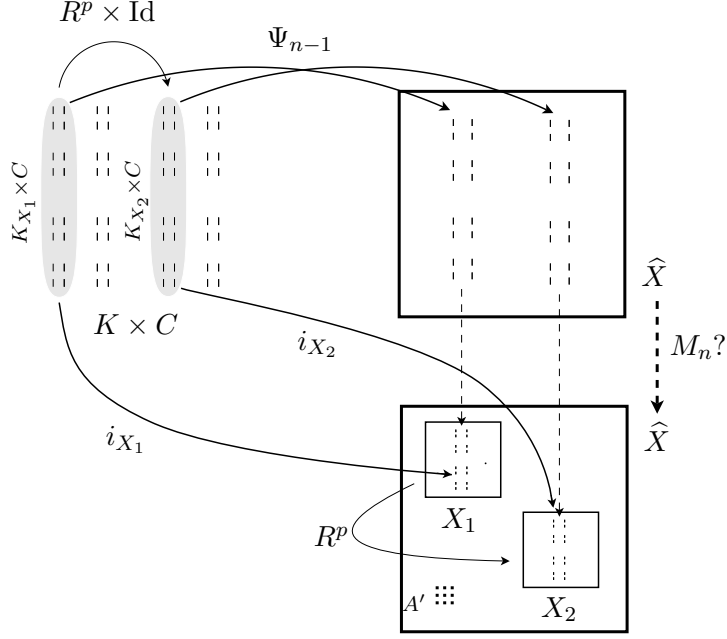}
\caption{Construction of $M_{n}$.}
\end{center}
\end{figure}

In a moment, we will extend the homeomorphism $M_{n}$ to the whole rectangle $\hat X$; before doing this, we first make an extension to a finite set that will ensure that the preimages of the sub-rectangles are ``thin''\footnote{There certainly are other ways to get this property, but this trick will also be useful in appendix~\ref{a.recurrence}.}. More precisely, we fix a small $\delta_{n} >0$ and we consider two finite subsets $A,A'$ in the interior of the rectangle $\hat X$, such that:
\begin{itemize}
\item[--] $A = \Psi_{n-1}(\wt A)$ where $\wt A$ is  disjoint from  the Cantor set $K_{\widehat X}\times C$ and $\delta _{n}$-dense in $\Psi_{n-1}^{-1}(\hat X)$;
\item[--] $A'$ is disjoint from all the sub-rectangles $X_{i}$;
\item[--] $A$ and $A'$ have the same number of elements.
\end{itemize}
We now extend $M_{n}$ to the set $A$ so that $M_{n}(A)=A'$.

Since the  set $A \cup \Psi_{n-1}(K_{\widehat X}\times C)$ is totally disconnected and tamely embedded in $\inte(\widehat X)$, using proposition~\ref{p.extension-cantor} of appendix~\ref{a.cantor} we can further extend $M_{n}$ to a homeomorphism of $\widehat X$ which is the identity on the boundary  of the rectangle.
 
Note that $M_n$ was constructed in such a way that, for every $i$, the rectangle $M_{n}^{-1}(X_{i})$ are disjoint from the set $A$. Since $\Psi_{n-1}^{-1}(A)$ is $\delta_{n}$-dense in $\Psi_{n-1}^{-1}(\hat X)$, this implies that:
\begin{equation}
\label{e.internal-diameter}
\mbox{for every }i,\mbox{ the internal radius of  }\Psi_{n-1}^{-1}(M_{n}^{-1}(X_{i}))\mbox{ is less than } \delta_{n}.
\end{equation}

\paragraph{Step 3. Definition of $\bm{M_n}$ on $\bm{\cM}$.}   
Let $\widehat X'$ be any rectangle of $\cE_{n-1}^{n-1}$. We define $M_n$ on $\widehat X'$ as follows. In step 2, we have defined $M_n$ on one (and only one) vertex $\widehat X$ which is in the same connected component of $\cG_{n-1}^{n-1}$ as $\widehat X'$. We consider the path $\widehat X \ra \cdots \ra \widehat X'=R^p(\widehat X)$ in $\cG_{n-1}^{n-1}$. Then we define $M_n$ on $\widehat X'$ by   
\begin{equation}
\label{e.def-M_n-2}
M_{n|\widehat X'}=R^p\circ M_{n|\widehat X}\circ R^{-p}.
\end{equation}
This defines the homeomorphism $M_n$ on $E_{n-1}^{n-1}$, ensures that $M_n$ preserves $E_{n-1}^{n-1}$ and is equal to the identity on the boundary of each connected component of $E_{n-1}^{n-1}$. As a consequence, we can extend $M_n$ on the whole manifold $\cM$ in such a way that $M_n$ is the identity outside $E_{n-1}^{n-1}$ (thus hypothesis $\mathbf{B_1}(n)$ is satisfied). Moreover, the definition of $M_n$ given above clearly implies that $M_n$ commutes with $R$ along the edges of $\cG_{n-1}^{n-1}$, i.e. satisfies hypothesis $\mathbf{B_2}(n)$. 

We now fix the value of the number $\delta _{n}$ of step 2 to get hypothesis $\mathbf{B_4}(n)$ (fibres are thin). Let $\ell$ be the maximal  length of the connected components of the graph $\cG_{n-1}^{n-1}$. Then $\delta _{n}$ may have been chosen  so small that for any $|p| \leq \ell$,  couple of points distant from less than $\delta _{n}$ are mapped by the homeomorphism $g_{n-1}^{p}$ to points distant from less than $\frac{1}{n}$. Now if $X'_{i}$ is an immediate sub-rectangle of $\widehat X'$, then $X_{i} = R^{-p}(X'_{i})$ is an immediate sub-rectangle of $\widehat X$. According to equality~\ref{e.internal-diameter} and the choice of $\delta _{n}$, the internal diameter of $g_{n-1}^p (\Psi_{n-1}^{-1}(M_{n}^{-1}(X_{i}))) = \Psi_{n-1}^{-1}(M_{n}^{-1}(X'_{i}))$ is less than $1/n$. This gives hypothesis $\mathbf{B_4}(n)$.

\paragraph{Step 4. Relation between $\bm{\Psi_n}$ and the $\bm{i_X}$'s.}
Note that the map $M_{n}$ is now defined, so that we can also deal with the map $\Psi_{n} = M_{n} \circ \cdots \circ M_1$. We claim that, for every rectangle $X\in\cE_{n}^0$, one has (see figure~\ref{f.Mn})  
\begin{equation}
\label{e.Psi-g}
\Psi_{n|K_{X}\times C}=i_{X|K_{X}\times C}.
\end{equation}

Observe that equality~\eqref{e.Psi-g} follows immediately from equality~(\ref{e.def-M_n})  when $X$ is included in one of the rectangles of $\cE_{n-1}^0$ on which $M_n$ was defined in step 2 (in particular, there nothing else to check for $n=1$; in what follows we assume $n \geq 2$). Now let $X'$ be any rectangle of $\cE_{n}^0$, and  denote by $\widehat X'$ the rectangle of $\cE_{n-1}^0$ containing $X'$. 
Let $\widehat X \ra \cdots \ra \widehat X'=R^p(\widehat X)$ be the (unique) path in the graph $\cG_{n-1}^{n-1}$ such that $\widehat X$ is a rectangle on which $M_{n}$ has been defined at step 2. Consider the rectangle $X = R^{-p}(X') \in \cE_{n}^0$. As explained above, equality~\eqref{e.Psi-g} holds for $X$.
Consider the following cube-shaped diagram. We aim at proving commutation of the right-hand side.
For this it suffices to prove  the commutation of the five other sides.
$$
\xymatrix{ 
&   K_{X} \times C  \ar@{>}[rr]  ^{R^p \times \mathrm{Id}}    \ar@{>}'[d] [ddd] ^(.3){i_{X}} & & K_{X'} \times C  \ar@{>}[ddd]  ^{i_{X'}} \\ 
K_{X} \times C \ar@{>}[ur]    ^{\mathrm{Id}}   \ar@{>}[rr]   ^(.65){R^p \times \mathrm{Id}}    \ar@{>}[ddd]   ^(.4){\Psi_{n}}  & & K_{X'} \times C \ar@{>}[ur]    ^{\mathrm{Id}}   \ar@{>}[ddd]   ^(.4){\Psi_{n}} \\
\\
& i_{X}(K_{X} \times C) \ar@{>}'[r]^(.8){R^p}[rr]    & & i_{X'}(K_{X'} \times C)\\
\Psi_{n}(K_{X}\times C) \ar@{>}[rr]   ^{R^p}  \ar@{>}[ur]  ^{\mathrm{Id}}  & & \Psi_{n}(K_{X'}\times C) \ar@{>}[ur]  ^{\mathrm{Id}} 
}$$
The top and bottom side of the diagram commutes trivially;
the commutation of the back side is due to the equivariance of the homeomorphisms $i_{X}$ (equality~\eqref{e.equivariance-i_X});
the commutation of the left-hand side is equality~\eqref{e.Psi-g} for $X$.
The commutation of the front side is obtained by putting together the induction 
hypothesis~$\mathbf{B_6}(n-1)$ and relation~\eqref{e.def-M_n-2}, yielding yet another  commutative diagram (see below).
This finishes the proof of equality~\eqref{e.Psi-g} for $X'$.
$$
\xymatrix{ 
K_{X}\times C \ar@{>}[d]   ^{\Psi_{n-1}} \ar@{>}[rr] ^{R^p \times \mathrm{Id}}
& &  K_{X'}\times C \ar@{>}[d]   ^{\Psi_{n-1}} \\
\Psi_{n-1}(K_{X}\times C) \ar@{>}[rr]   ^{R^p} \ar@{>}[d]   ^{M_{n}} & & \Psi_{n-1}(K_{X'}\times C)\ar@{>}[d]   ^{M_{n}} \\
\Psi_{n}(K_{X}\times C) \ar@{>}[rr]   ^{R^p}  & & \Psi_{n}(K_{X'}\times C)
}
$$

\paragraph{Step 5. Hypotheses $\mathbf{B_{1,2,4,5,6}}(n)$ are satisfied.}  It follows from the very
construction of $M_n$  that hypothesis $\mathbf{B_{1}}(n)$,  $\mathbf{B_2}(n)$ and   $\mathbf{B_4}(n)$   are satisfied (see the end of step 3). Equality~\eqref{e.Psi-g} implies hypothesis $\mathbf{B_5}(n)$ (Cantor sets in the fibres). Equality~\eqref{e.Psi-g} and the equivariance property~\eqref{e.equivariance-i_X} imply  hypothesis $\mathbf{B_6}(n)$ (embedding a trivial dynamics).  This completes the proof of the proposition. 
\end{proof}

\section{Extraction process}
\label{s.extraction}

Recall that we are considering a sequence of collections of rectangles $(\cE_n^0)_{n\in\NN}$ such that hypotheses $\mathbf{A_1}$ and $\mathbf{A_3}$ are satisfied.
In this section, we assume moreover that we are given a sequence $(M_n)_{n\geq 1}$ of homeomorphisms of $\cM$ satisfying hypothesis $\mathbf{B_{1}}$. The purpose of this section is to explain an extraction process which allows to replace the sequence of collections of rectangles $(\cE_n^0)_{n\in\NN}$ by a sub-sequence $(\cE_{n_k}^0)_{k\in\NN}$, and the  sequence of homeomorphisms $(M_n)_{n\geq 1}$ by a sequence $(\bM_k)_{k\in\NN}$. The first important point is that this extraction process will respect hypotheses $\mathbf{B_{1,2,4,5,6}}$. The second important point is that this extraction process can be used to get some convergence properties; this is the reason why it will play an important role in part~\ref{part.H_k}. To illustrate this role, we will explain below how to use the extraction process to get hypothesis $\mathbf{B_3}$.

\subsection{Definition  of the extraction process}
\label{ss.def-extraction}

Assume we are given an increasing sequence of integers $(n_k)_{k\in\NN}$, with $n_0=0$. Note that $n_{k} \geq k$ for every $k$.
Then we consider 
\begin{itemize}
\item[--] for every $k\in\NN$ and every $p\leq k+1$, the collection of rectangles  $\bcE_k^p$, the set $\bE_k^p$, and the graph $\bcG_k^p$ defined by 
$$
\bcE_k^p:=\cE_{n_k}^p\quad\quad \bE_k^p:=E_{n_k}^p\quad\quad \bcG_k^p:=\cG_{n_k}^p;
$$
\item[--] for every $k\geq 1$, the homeomorphism $\bM_k$ defined by 
$$
\bM_k:=M_{n_k|E_{n_k-1}^{k-1}}\circ\dots\circ M_{i|E_{i-1}^{k-1}} \circ\dots\circ M_{n_{k-1}+1|E_{n_{k-1}}^{k-1}}
$$
where $M_{i|E_{i-1}^{k-1}}$ denotes the homeomorphism\footnote{Hypothesis $\mathbf{B_1}$ implies that, for every $i\geq k$, the set $E_{i-1}^{k-1}$  is the union of some of the connected components of the support of the homeomorphism $M_i$. Hence, the map  $M_{i|E_{i-1}^{k-1}}$ is a homeomorphism.} which is equal to $M_i$ on $E_{i-1}^{k-1}$ and equal to the identity on $\cM\setminus E_{i-1}^{k-1}$; 
\item[--] for every $k\geq 1$, the maps $\bPsi_k$ and $\mathbf{g}_k$ defined by 
$$\bPsi_k:=\bM_{k}\circ\dots\circ\bM_1\quad\mbox{ and }\quad
\mathbf{g}_k:=\bPsi_k^{-1}\circ R\circ\bPsi_k.$$
\end{itemize}

\begin{rema}
\label{r.extraction-fibres}
It follows from these definitions that  $\bPsi_{k}^{-1}= \Psi_{n_k}^{-1}$ on the set $\bE_{k-1}^0=E_{n_{k-1}}^0$. This can be proved recursively, using the equality $\bM_{k}(\bE_{k-1}^0) = \bE_{k-1}^0$.
As a consequence,   $\bPsi_{k}^{-1} (\bE_{k-1}^0)= \Psi_{n_k}^{-1}( E_{n_{k-1}}^0)$,  and on this set we have   $\bPsi_{k}= \Psi_{n_{k}}$.
\end{rema}

\begin{rema*}
Applying two successive extractions with some sequences of integers $(n_k)_{k\geq 0}$ and $(k_l)_{l\geq 0}$ amounts to applying a single extraction with the sequence $(n_{k_l})_{l\geq 0}$. 
\end{rema*}

\subsection{Hypotheses $\mathbf{B_{1,2,4,5,6}}$ are preserved}
\label{ss.extraction-preserve-hyp}

Assume we are given an increasing sequence of integers $(n_k)_{k\in\NN}$. We consider the sequence of collection of rectangles $(\bcE_k^0)_{k\in\NN}$ and the sequence of homeomorphisms $(\bM_k)_{k\in\NN}$ defined above. Then, we can consider the conditions obtained by replacing the sequences $(\cE_n^0)_{n\in\NN}$ and $(M_n)_{n\geq 1}$ by the sequences $(\bcE_k^0)_{k\in\NN}$ and $(\bM_k)_{k\geq 1}$ in hypotheses $\mathbf{A_{1,2,3}}$ and $\mathbf{B_{1,2,3,4,5,6}}$. We obtain some new conditions that we still call ``hypotheses $\mathbf{A_{1,2,3}}$ and  $\mathbf{B_{1,2,3,4,5,6}}$". For example, the new hypotheses $\mathbf{A_{1.a}}$ and $\mathbf{B_6}$ are:
\begin{itemize}
\item[]
\begin{itemize}
\item[$\mathbf{A_{1.a}}$] For every $k\in \NN$, the collection ${\bcE}^0_k$ is $k+1$ times iterable and the collection $\bcE_k^{k}$ has no cycle;
\item[$\mathbf{B_6}(k)$] If $X \ra \cdots \ra X'=R^p(X)$ is a path in the graph $\bcG_{k}^k$, with $X,X'\in \bcE_k^0$, then 
$$\bPsi_k^{-1}\circ R^p\circ\bPsi_k(x,c)=(R^p(x),c)\quad\mbox{for every}\quad(x,c)\in K_X\times C.$$
\end{itemize}
\end{itemize}

A very important feature of the extraction process defined above is that it does not destroy the hypotheses.

\begin{prop}
\label{p.preserved-by-extraction}
\begin{itemize}
\item[--] The extracted sequence of collections of rectangles $(\bcE_k^0)_{k\in\NN}$ satisfies (the new) hypotheses $\mathbf{A_{1,3}}$. Moreover, if the original sequence $(\cE_n^0)_{n\in\NN}$ satisfies hypothesis $\mathbf{A_2}$, then the extracted sequence $(\bcE_k^0)_{k\in\NN}$ also satisfies $\mathbf{A_2}$.
\item[--] The extracted sequence of homeomorphisms $(\bM_k)_{k\geq 1}$ satisfies hypothesis $\mathbf{B_{1}}$. Moreover, if $\mathbf{B_*}$ is one of the hypotheses $\mathbf{B_2}$, $\mathbf{B_4}$, $\mathbf{B_5}$, $\mathbf{B_6}$, and if the original sequence of homeomorphisms $(M_n)_{n\geq 1}$ satisfies $\mathbf{B_*}$, then the extracted sequence $(\bM_k)_{k\geq 1}$ also satisfies $\mathbf{B_*}$.
\end{itemize}
\end{prop}

\begin{proof}
The proof is immediate for hypotheses $\mathbf{A_{1.a}}$, $\mathbf{A_{1.b}}$ and $\mathbf{A_3}$. For hypotheses $\mathbf{A_{1.c}}$ and $\mathbf{A_2}$, it follows from a very easy induction.
For $\mathbf{B_1}$, it follows from the very definition of the sequence of homeomorphisms $(\mathbf{M_k})$.  For $\mathbf{B_2}$, one just needs to observe that $\bM_{k}=M_{n_k}\circ\dots\circ M_{n_{k-1}+1}$ on every rectangle $X$ which is a vertex of the graph $\bcG_{k-1}^{k-1}=\cG_{n_{k-1}}^{k-1}$, and that all the homeomorphisms $M_{n_k},\dots,M_{n_{k-1}+1}$ commute with $R$ along the edges of the graph $\bcG_{k-1}^{k-1}$. Finally, for $\mathbf{B_4}, \mathbf{B_5}, \mathbf{B_6}$, it follows from remark~\ref{r.extraction-fibres}.
\end{proof}

\subsection{Realisation of hypothesis $\mathbf{B_3}$}
\label{ss.realis-B_3}

Proposition~\ref{p.get-B_3} shows that the extraction process can be used to get hypothesis $\mathbf{B_3}$.  Formally, we will not use it (although we will use a very similar statement in part~\ref{part.H_k}). We only state it for ``pedagogical reasons" to illustrate the use of the extraction process.

\begin{prop}
\label{p.get-B_3}
There exists an increasing sequence of integers $(n_k)_{k\geq 1}$ such that hypothesis $\mathbf{B_3}$ is satisfied for the extracted sequences $(\bcE_k^0)_{k\in\NN}$ and $(\bM_k)_{k\geq 1}$.
\end{prop}

\begin{proof}
Choose a sequence of positive real numbers $(\epsilon_k)_{k\in\NN}$ such that $\epsilon_k\to 0$, and assume that, for some $k \geq 1$, some integers $n_1<\dots<n_{k-1}$ have been constructed in such a way that 
$\max\{\diam(\bPsi_{l-1}^{-1}(X))\mid  X\in\bE_l^{l+1}\setminus \bE_{l}^{l-1}\}\leq\epsilon_l$ for $1\leq l\leq k-1$.

On the one hand, we know the integer $n_{k-1}$. Hence we know the homeomorphism $\bPsi_{k-1}$. Hence we can find $\eta>0$ such that, if $X$ is a rectangle of diameter less than $\eta$, then $\bPsi_{k-1}^{-1}(X)$ is a set of diameter less than $\epsilon_k$.  On the other hand, hypothesis $\mathbf{A_3}$ gives us an integer $n_k>n_{k-1}$ such that the diameter of every rectangle of $E_{n_k}^{k+1}\setminus E_{n_k}^{k-1}$ is less than $\eta$. Hence, we have $\max\{\diam(\bPsi_{k-1}^{-1}(X))\mid  X\in\bE_k^{k+1}\setminus \bE_{k}^{k-1}\}\leq\epsilon_k$, as wanted.
\end{proof}

\newpage
\part{Insertion of the desired dynamics in the fibres of $\Psi$}
\label{part.H_k}

In this part, we assume that we are given a sequence of collections of rectangles $(\cE^0_n)_{n\in\NN}$ such that hypotheses $\mathbf{A_{1,3}}$ are satisfied and a sequence of homeomorphisms $(M_n)_{n\geq 1}$ such that hypotheses $\mathbf{B_{1,2}}$ are satisfied\footnote{In sections~\ref{s.insert-dynamics} and~\ref{s.control}, we will also assume that hypotheses $\mathbf{B_{5,6}}$ are satisfied, but we do not need theses hypotheses for section~\ref{s.general-H}.}.
\begin{itemize}
\item[--] In section~\ref{s.general-H}, we introduce a sequence of homeomorphisms $(H_{k})_{k \geq 1}$ and hypotheses $\mathbf{C_{1,2,3,4}}$ that are the analogue of hypotheses $\mathbf{B_{1,2,3,4}}$ for the maps $M_{n}$. In particular, under these hypotheses, the sequence $(H_{k} \circ \bM_{k})$ will still satisfy hypotheses $\mathbf{B_{1,2,3,4}}(k)$, so the convergence results of part~\ref{part.M_n} will also hold for this sequence: thus we get a continuous onto map  $\Phi = \cdots \circ H_{k} \circ \bM_{k} \circ \cdots \circ H_{1} \circ \bM_{1}$ and a homeomorphism $f:\cM\rightarrow\cM$ such that $\Phi\circ f=R\circ \Phi$.
\item[--]  In section~\ref{s.insert-dynamics}, we will state two additional hypotheses (denoted by $\mathbf{C_5}$ and $\mathbf{C_6}$) in order to ensure that the dynamics of $f$ on the orbit of $K\times C$ is isomorphic to a given fibered map $h$. This section also explains how to construct  a sequence $(H_{k})$ satisfying all hypotheses $\mathbf{C_{1, \dots , 6}}$. It ends with the first part of the proof of our main theorem~\ref{t.main}, which contains Mary Rees original result.
\item[--] Finally, in section~\ref{s.control}, we will formulate hypotheses $\mathbf{C_{7,8}}$
under which  the dynamics of $f$ on the whole of $\cM$  is isomorphic to $h$.  Then we will tighten the construction of section~\ref{s.insert-dynamics} in order to get a sequence  $(H_{k})$ satisfying the additional hypotheses $\mathbf{C_{7,8}}$. This enables us to complete the proof of theorem~\ref{t.main}.
\end{itemize}

\section{General scheme}
\label{s.general-H}

\subsection{The sequences $(n_k)_{k\geq 1}$ and $(H_k)_{k\geq 1}$}
\label{ss.def-H}

We consider an extraction sequence $(n_k)_{k\geq 1}$ and a sequence $(H_k)_{k\geq 1}$ of homeomorphisms of $\cM$. Given these sequences, for every $k\geq 1$ we define  the homeomorphism $\Phi_k$ by:
$$
\Phi_k:=H_k\circ \bM_{k}\circ\dots\circ H_2\circ\bM_2\circ H_1\circ\bM_1
$$
and the homeomorphism
$$
f_k:=\Phi_k^{-1}\circ R\circ\Phi_k.
$$
One also sets $\Phi_0=\id$ and $f_0=R$.

\subsection{Hypotheses $\mathbf{C_{1,2,3,4}}$}
\label{ss.hyp-C_1234}

We consider the following hypotheses.
\begin{itemize}
\item[]
\begin{itemize}
\item[$\mathbf{C_1}(k)$] \emph{(Support)}\\
 The support of the homeomorphism $H_k$ is contained in the set $\bE_{k}^{k-1}$. 
\item[$\mathbf{C_2}(k)$] \emph{(Commutation)}\\
The maps $H_k$ and $R$ commute along the edges of the graph $\bcG_{k}^{k-1}$.
\item[$\mathbf{C_3}$] \emph{(Convergence)} 
Let $\bcA_{k} = \bcE_{k}^{k+1} \setminus  \bcE_{k}^{k-1}$. Then the supremum of the diameters of the rectangles $\Phi^{-1}_{k-1}(X))$ with $X \in \bcA_{k}$ tends to $0$
when $k$ tends to $+\infty$.
\item[$\mathbf{C_4}$] \emph{(Fibres are thin)}\\
The internal radius of the set $\Phi_k^{-1}(\bE_k^0)$ goes to $0$ when $k\to\infty$.
\end{itemize}
\end{itemize}
The last two hypotheses have also quantitative versions.
\begin{itemize}
\item[]
\begin{itemize}
\item[$\mathbf{C_3}(k)$]\emph{(Convergence, quantitative version)} \\
$$
\displaystyle\max \left\{ \diam(\Phi^{-1}_{k-1}(X)) \mid    \; X\in \bcE_k^{k+1}\setminus \bcE_k^{k-1} \right\}  \leq \frac{1}{k}.
$$
\item[$\mathbf{C_4}(k)$]\emph{(Fibres are thin, quantitative version)}\\
The internal radius of the set $\Phi_k^{-1}(\bE_k^0)$ is less than $\frac{1}{k}$.
\end{itemize}
\end{itemize}

\begin{rema*}
There is an important difference between the assumptions on $(M_n)_{n\geq 1}$ and their analogues for $(H_k)_{k\geq 1}$: hypothesis $\mathbf{C_1}$ requires the support of the homeomorphism $H_k$ to be included in the set $\bE_{k}^{k-1}=E_{n_k}^{k-1}$, whereas hypothesis $\mathbf{B_1}$ implies that the support of the homeomorphism $\bM_{k}$ is contained in the set $\bE_{k-1}^{k-1}=E_{n_{k-1}}^{k-1}$. Consequently, if all the hypotheses are satisfied and if the integer $n_k$ is much bigger than $n_{k-1}$, then the size of the support of the homeomorphism $H_k$ will be much smaller than the size of the support of the homeomorphism $\bM_k$. This will play a crucial role when we will try to get hypotheses $\mathbf{C_{3}}$ and $\mathbf{C_4}$ (see in particular proposition~\ref{p.existence-H-1}). 

Note also that there is no real difference between $\mathbf{C_2}$ and $\mathbf{B_2}$: the compatibility hypothesis $\mathbf{A_{1.c}}$ implies that, given hypothesis $\mathbf{C_1}(k)$ on the support,  nothing changes if one replaces in $\mathbf{C_2}(k)$ the graph  $\bcG_{k}^{k-1}$ with the graph  $\bcG_{k-1}^{k-1}$, thus getting the exact analogue to hypothesis $\mathbf{B_2}(n)$.
\end{rema*}

One checks easily that if one combines these hypotheses with hypotheses $\mathbf{B_{1,2}}$ on $(M_n)_{n\geq 1}$, one can replace the sequences $(\cE_n^0)_{n\geq 0}$, $(M_n)_{n\geq 1}$ and $(\Psi_n)_{n\geq 0}$ by the sequences $(\bcE_k^0)_{k\geq 0}$, $(H_k\circ\bM_k)_{k\geq 1}$ and $(\Phi_k)_{k\geq 0}$: for each hypothesis $\mathbf{C_i}$ which is satisfied by $(H_k)$, the corresponding hypothesis $\mathbf{B_i}$ is satisfied by $(H_k\circ \bM_k)$. 

In particular, one can apply propositions~\ref{p.extension}, \ref{p.extension2} and~\ref{p.recurrence} to these new sequences and obtain the following consequences.

\begin{prop}
\label{p.extension-bis}
Propositions~\ref{p.extension}, \ref{p.extension2} and~\ref{p.recurrence} still hold when one replaces $\Psi_{n}$, $\Psi$, $g_{n}$, $g$, and $\cE_{n}^{n_{0}}$ respectively with $\Phi_{k}$, $\Phi$, $f_{k}$, $f$, and $\bcE_{k}^{k_{0}}$.
\end{prop}

\section{Inserting the dynamics of $h$ on the Cantor set $K \times C$}
\label{s.insert-dynamics}

In this section, we assume hypotheses $\mathbf{A_{1,2,3}}$ for the sequence of collections of rectangles $(\cE_n^0)_{n\geq 0}$ and hypotheses $\mathbf{B_{1,2,5,6}}$ for the sequence of homeomorphisms $(M_n)_{n\geq 1}$. We also assume that the graph $\cG_0^0$ has no edge. Considering a fibered map $h$ as in  theorem~\ref{t.main}, we state two additional hypotheses $\mathbf{C_{5,6}}$ for the sequences $(n_k)$ and $(H_k)$, that are the counterpart (but not the exact analogues) of hypotheses $\mathbf{B_{5,6}}$: they will imply that the homeomorphism $f$ realises the dynamics of  $h$ on the set $\bigcup_{i\in \ZZ}R^i(K)\times C$.

Then we  explain  how to construct a sequence of homeomorphisms $(H_k)$ and a sequence of integers $(n_k)$ so that hypotheses $\mathbf{C_{1,2,3,4,5,6}}$ will be satisfied.
There is no simple way to define an extraction for the sequence $(H_k)$; hence, one has to  combine the inductive construction of the sequence $(H_k)$ with successive extractions of the sequences $(\cE_n^0)$ and $(M_n)$. This is structured as follows (in subsection~\ref{ss.bilan}): 
\begin{itemize}
\item[--] we assume that the homeomorphism $H_{\ell}$ and the integer $n_{\ell}$ are known for every $\ell\leq k-1$;
\item[--] then, we construct the integer $n_{k}$ in such a way that hypotheses $\mathbf{C_{3,4}}(k)$ will be satisfied whatever the homeomorphism $H_{k}$ might be (see proposition~\ref{p.existence-H-1}); 
\item[--] and then, we define the homeomorphism $H_{k}$ so that hypotheses $\mathbf{C_{1,2,5,6}}(k)$ will be satisfied (see proposition~\ref{p.existence-H-2}).
\end{itemize}

\subsection{Hypothesis $\mathbf{C_5}$ and some consequences}
\label{ss.hyp-C5}

We consider the following hypotheses:
\begin{itemize}
\item[]
\begin{itemize}
\item[$\mathbf{C_{5}}(k)$] \emph{(Cantor sets in the fibres of $\Phi$)}~\\
For each point $x \in K$, the map $H_{k}$ preserves $\bPsi_{k} (\{x\} \times C)$:
 $$H_{k}(\bPsi_{k} (\{x\} \times C)) = \bPsi_{k} (\{x\} \times C).$$
\end{itemize}
\end{itemize}

Under hypothesis $\mathbf{C_{5}}$, the map $\Phi_{k}$ satisfies a property  analogous to hypothesis $\mathbf{B_{5}}$ for $\Psi_{n}$. More precisely, we get the following consequences (which are illustrated by the topological flavour of figure~\ref{f.isomorphism} of the introduction, subsection~\ref{ss.outline}).

\begin{prop}
\label{p.H-rectangle}
Let $(n_{1}, \dots , n_{k})$ and $(H_{1},\dots,H_{k})$
be two finite sequences that satisfy hypotheses $\mathbf{C_{5}}(\ell)$ for every $\ell\leq k$.
\begin{enumerate}
\item For each $x\in K$, we have
$$\Phi_k(\{x\} \times C)=\bPsi_k(\{x\} \times  C)$$
and the fibre $\Phi^{-1}(x)$ contains $\{x\} \times C$.
\item If $X \ra \cdots \ra X'=R^p(X)$ is a path in  the graph $\bcG_k^k$
with $X,X'\in \bcE^0_{k}$, then, for each $x\in K_X$ we have
$$f_{k}^p(\{x\} \times C)  = \{R^p(x)\} \times C.$$
\end{enumerate}
\end{prop}

\begin{proof}
The equality $\Phi_k(\{x\} \times C)=\bPsi_k(\{x\} \times  C)$ is easily proved using an induction on $\ell$. 
Now $\Phi^{-1}(x) = \cap_{k \in\NN} \Phi_{k}^{-1}(X_{{k}})$ where $x \in X_{{k}} \in \bcE_{k}^0$ (see propositions~\ref{p.extension2} and~\ref{p.extension-bis}). Using the above equality and hypothesis $\mathbf{B_{5}}$ on $\bPsi_{k}$, one gets  $\Phi_k(\{x\} \times C) \subset X_{k}$ for every $k$. Thus  $\{x\} \times C$ is contained in $\Phi^{-1}(x)$. This completes the proof of the first property. 

For the second property, one uses hypothesis $\mathbf{B_{6}}$ for $\bPsi_{k}$, which implies the equality 
$$
\bPsi_{k}^{-1} \circ R^p \circ \bPsi_{k} (\{x\} \times C) = \{R^p(x)\} \times C.
$$
The first property shows that in this equality $\bPsi_k$ may be replaced by $\Phi_k$. This gives the desired equality.
\end{proof}

\subsection{The fibered map $h$}

We now consider a fibered dynamics (compare to the hypotheses of theorem~\ref{t.main}):
\begin{eqnarray*}
h:\bigcup_{i\in\ZZ} R^i(K)\times C & \longrightarrow   &  \bigcup_{i\in\ZZ} R^i(K)\times C \\
 (x,c)  & \longmapsto & (R(x), h_{x}(c))
\end{eqnarray*}
such that
\begin{enumerate}
\item $h$ is bijective;
\item  for every $\displaystyle x \in \bigcup_{i \in\ZZ} R^i(K)$, $h_{x}$ is a homeomorphism of $C$;
\item  for every integer $k$, the map $h^k$ is continuous on $K \times C$. 
\end{enumerate}
In particular, if $X \ra \cdots \ra X'=R^p(X)$ is a path in  the graph $\bcG_k^k$ with $X,X'\in \bcE^0_{k}$, then the map $h^p$ induces a homeomorphism between $K_{X} \times C$ and $K_{X'} \times C$. Note that by proposition~\ref{p.H-rectangle}, the map $f_{k}^p$ shares the same property (for the embedded versions of the Cantor sets $K_{X} \times C$ and $K_{X'} \times C$).

\subsection{Hypothesis $\mathbf{C_6}$ and some consequences}
\label{ss.hyp-C6}

\begin{itemize}
\item[]
\begin{itemize}
\item[$\mathbf{C_{6}}(k)$] \emph{(Embedding a non-trivial dynamics)} ~\\
If $X \ra \cdots \ra X'=R^p(X)$ is a path in  the graph $\bcG_k^k$ with $X,X'\in \bcE^0_{k}$, then
$$
f_{k}^p(x,c) = h^p (x,c)\quad\mbox{ for every }\quad(x,c) \in  K_{X} \times C.
$$
\end{itemize}
\end{itemize}

We get for $f$ a statement similar to  proposition~\ref{p.conjugacy} for $g$ (and the proof is the same). The proposition below is illustrated in figure~\ref{f.isomorphism} of the introduction.

\begin{prop}
\label{p.conjugacy-bis}
Assume  that hypotheses $\mathbf{C_{1,2,3,5,6}}$ are satisfied. Then, the map  $f$ on $\displaystyle\bigcup_{i\in \ZZ} f^i(K\times C)$ is isomorphic to the map $h$ on $\displaystyle \bigcup_{i\in \ZZ} R^i(K)\times C$.
\end{prop}

\subsection{Construction of the sequence $(n_{k})$: realisation of hypotheses~$\mathbf{C_{3,4}}$}
\label{ss.construc-n}

In this section, we explain how to get $\mathbf{C_{3}}$ and $\mathbf{C_{4}}$ from  hypotheses $\mathbf{A_3}$ and $\mathbf{B_{4}}$.

\begin{prop}
\label{p.existence-H-1}
Let us assume that $(M_n)$ satisfies the additional hypotheses $\mathbf{B_{4}}$. Let $(n_{1}, \dots , n_{k-1})$ and $(H_{0},\dots,H_{k-1})$ be two finite sequences. Then, for any integer $n_k>n_{k-1}$ large enough and for any homeomorphism $H_k$ that satisfies $\mathbf{C_1}(k)$, hypotheses $\mathbf{C_{3}}(k)$ and $\mathbf{C_4}(k)$ also hold.
\end{prop}

The idea of the proof is  to extract the sequence $(E_n^0)$ so that the rectangles of $\bcE_k^0$ have a very small diameter with respect to the modulus of continuity of the homeomorphisms $\bM_\ell$ and $H_\ell$ for $\ell<k$.

\begin{proof}
By hypothesis $\mathbf{A_3}$, for $n$ large enough the diameter of the connected components of $E^{k+1}_n$ is arbitrarily small. One deduces that for $n$ large enough the diameter of the connected components of $\Phi_{k-1}^{-1}(E^{k+1}_n\setminus E^{k-1}_n)$ is less than $\frac 1 k$. This gives hypothesis $\mathbf{C_3}(k)$.

By hypothesis $\mathbf{B_4}$, for any $\rho>0$ and for $n$ large enough the set $\Psi_n^{-1}(E^0_n)$ does not contain any ball of radius $\rho$. By definition of the extractions (see remark~\ref{r.extraction-fibres}), we have
$$
\Psi_n^{-1}(E^0_n)=\bPsi_{k-1}^{-1}\circ\left(M_n\circ\dots\circ M_{n_{k-1}+1}\right)^{-1}(E^0_n).
$$
Hence, if $\rho$ is small, the set
$$
\left(\Phi_{k-1}^{-1}\circ \bPsi_{k-1}\right) \circ  \Psi_n^{-1}(E^0_n) = (M_n\circ\dots\circ M_{n_{k-1}+1}\circ \Phi_{k-1})^{-1}(E^0_n)
$$
does not contain any ball of radius $\frac 1 k$. Choosing $n_k$ large enough, one gets the same property for the set $(\bM_k\circ\Phi_{k-1})^{-1}(\bE^0_k)$. Let us consider any homeomorphism $H_k$ satisfying hypothesis $\mathbf{C_1}(k)$: the support of $H_k$ is included in $\bE^{k-1}_k$. One deduces that $\Phi^{-1}_k(\bE^0_k)=(\bM_k\circ\Phi_{k-1})^{-1}(\bE^0_k)$. This gives hypothesis $\mathbf{C_4}(k)$.
\end{proof}

\subsection{Construction of the $H_{k}$'s:  realisation of hypotheses $\mathbf{C_{1,2,5,6}}$}
\label{ss.realis-C_1256}

We explain in this section how to build inductively the maps $H_{k}$ in order to satisfy
hypotheses $\mathbf{C_{1,2,5,6}}$.

\begin{prop}
\label{p.existence-H-2}
Let us assume that $\cG_0^0$ has no edge. Let $k \geq 1$ and  let $(n_{1}, \dots , n_{k})$ and   $(H_{1},\dots,H_{k-1})$ be two finite sequences that satisfy hypotheses $\mathbf{C_{1,2,5,6}}(\ell)$ for every $1 \leq \ell <k$. Then there exists a homeomorphism $H_{k}$ of $\cM$ such that hypotheses $\mathbf{C_{1,2,5,6}}(k)$ are also satisfied.
\end{prop}

\begin{proof}
As for the construction of the homeomorphisms $M_{n}$ in proposition \ref{p.existence-M}, we build $H_{k}$ independently on each connected component of the graph $\bcG_{k}^k$: we fix $\Gamma$ one of these components. We also denotes by $\Gamma_{1},\dots, \Gamma_{s}$ the connected components
of $\bcG_{k}^{k-1}$ contained in $\Gamma$; they are ordered by the dynamics on $\bcG_{k}^k$ (which has no cycle by hypothesis $\mathbf{A_{1}}$). For each connected component $\Gamma_{i}$, one chooses a rectangle $X_{i}\in \Gamma_{i}\cap \bcE^0_{k}$. The construction goes as follows. On step~1 we will  define inductively $H_{k}$ on each $X_{i}$, making sure that hypothesis $\mathbf{C_{5}}(k)$ is satisfied on each of these rectangles, and that  hypothesis $\mathbf{C_{6}}(k)$ is satisfied between each pair $X_{i}, X_{i+1}$. On step~2 we will extend $H_{k}$ along each connected component $\Gamma_i$ of the graph $\bcG_{k}^{k-1}$ so that the commutation with $R$ holds (hypothesis~$\mathbf{C_{2}}(k)$).
On step~3 we will check that hypothesis~$\mathbf{C_{6}}(k)$ holds between \emph{any} pair of rectangles in $\Gamma\cap \bcE_{k}^0$.

\bigskip

\noindent\textit{Step 1. Definition of $H_k$ on one vertex of each connected component of the graph $\bcG_{k}^{k-1}$.\;} One defines inductively $H_{k}$ on each $X_{i}$ in the following way.

On $X_{1}$, we set\footnote{We have some freedom here: $H_{1}$ should be the identity on the boundary of $X_{1}$ and preserve each Cantor set $\bPsi_{k}(C_{x})$ for $x\in K_{X_{1}}$.} $H_{k}=\id$.

Suppose that  $H_{k}$ has been defined on $X_{i}$ (so that $\mathbf{C_{5}}(k)$ is satisfied for every $x \in K_{X_{i}}$). Consider the integer $p_{i}\geq 1$ such that $X_{i+1}=R^{p_{i}}(X_{i})$. We want to define $H_{k}$ on $X_{i+1}$ in order to satisfy $\mathbf{C_{6}}(k)$, i.e. such that on $K_{X_{i}}\times C$ we have $f^{p_{i}}_{k}=h^{p_i}$. That is, we want the following diagram to commute:
$$
\xymatrix{ 
K_{X_{i}}\times C   \ar[dd]^{f_{k}^{p_{i}} = h^{p_{i}}} \ar[r]^(.4){\Phi_{k-1}} & 
\bPsi_{k-1}(K_{X_{i}}\times C)  \ar[r]^{\bM_{k}} & 
\bPsi_{k}(K_{X_{i}}\times C)  \ar[r]^{H_{k}} & 
\bPsi_{k}(K_{X_{i}}\times C)  \ar[dd]^{R^{p_{i}}}
 \\
 \\
K_{X_{i+1}}\times C  \ar[r]^(.4){\Phi_{k-1}}  &
\bPsi_{k-1}(K_{X_{i+1}}\times C)  \ar[r]^{\bM_{k}}  &
\bPsi_{k}(K_{X_{i+1}}\times C)  \ar[r]^{H_{k}}_{?}  &
\bPsi_{k}(K_{X_{i+1}}\times C)
}
$$

This allows to define $H_{k}$ on $\bPsi_{k}(K_{X_{i+1}}\times C)$ in a unique way as:
$$
H_{k}=R^{p_{i}}\circ H_{k}\circ \bM_{k}\circ \Phi_{k-1}\circ h^{-p_i}\circ \Phi_{k-1}^{-1}\circ \bM_{k}^{-1}.$$
Note that
\begin{itemize}
\item[--] Writing down the diagram makes use of
\begin{itemize}
\item proposition~\ref{p.H-rectangle}, which entails $\Phi_{k-1}(K_{X_{i}}\times C)= \bPsi_{k-1}(K_{X_{i}}\times C)$;
\item hypothesis~$\mathbf{C_{5}}(k)$ on $X_{i}$, which entails $H_{k}(\bPsi_{k}(K_{X_{i}}\times C)) = \bPsi_{k}(K_{X_{i}}\times C)$;
\item  hypothesis~$\mathbf{B_{6}}(k)$, which entails that $R^{p_{i}}(\bPsi_{k}(K_{X_{i}}\times C) ) = \bPsi_{k}(K_{X_{i+1}}\times C)$.
\end{itemize}
\item[--] The homeomorphism $H_{k}$ that appears in the top side of the diagram is the homeomorphism of$X_{i}$ that has already be defined (by induction on $i$). The other one is the homeomorphism of $X_{i+1}$ we want to define.
\item[--] By our continuity assumption on $h$, the restriction of the map $h^{p_i}$ on $K_{X_{i}}\times C$
is a homeomorphism, proving that $H_{k}$ is  homeomorphism of $\bPsi_{k}(K_{X_{i+1}}\times C)$.
\item[--] The fibered structures on the Cantor sets are preserved by each of these maps: one uses proposition~\ref{p.H-rectangle} for $\Phi_{k-1}$, the assumption that $h$ is a fibered map, hypothesis $\mathbf{B_{6}}(k)$ for $R^{p_{i}}$ and the induction assumption for $H_{k}$ on $\bPsi_{k}(K_{X_{i}}\times C)$. This proves that $\mathbf{C_{5}}(k)$ is satisfied for every $x \in K_{X_{i+1}}$
\end{itemize}

Note that $\bPsi_{k}(K_{X_{i+1}}\times C)$ is included in the interior of the rectangle $X_{i+1}$ (hypothesis $\mathbf{B_{5}}(k)$). One then extends the map $H_{k}$ to a homeomorphism of $X_{i+1}$
which is the identity on the boundary, using proposition \ref{p.extension-cantor} of appendix \ref{a.cantor}.

\bigskip

\noindent\textit{Step 2. Extension of $H_k$ to $\cM$.\;} 
We first define $H_{k}$ on the rectangles of each component $\Gamma_{i}$: for any rectangle $X\in \Gamma_{i}$, we consider the integer $p$ such that $X=R^p(X_{i})$ and set
$$
H_{k|X}=R^p\circ H_{k|X_{i}}\circ R^{-p}.
$$
The same construction is carried out for every connected component $\Gamma$ of $\bcG_{k}^k$. We extend $H_{k}$ by the identity elsewhere. One clearly obtains a homeomorphism $\cM$ which satisfies hypothesis $\mathbf{C_{1}}(k)$ on the support and $\mathbf{C_{2}}(k)$ on the commutation.

\bigskip

\noindent\textit{Step 3. Hypotheses $\mathbf{C_{5,6}}(k)$ are satisfied.\;}
For each rectangle $X=X_{i}$, the map $H_{k}$ preserves the Cantor set $\bPsi_{k}(K_{X}\times C)$ and its fibered structure by construction. Any other rectangle $X\in \bcG_{k}^{0}$ is the iterate of some rectangle $X_{i}$ that belongs to the same connected component of $\bcG_{k}^{k-1}$ as $X$.
Hence, the same property holds for any $X$ by using the equivariance given by $\mathbf{B_{6}}(k)$.
This proves hypothesis $\mathbf{C_{5}}(k)$.

Let us prove hypothesis $\mathbf{C_6}(k)$. We consider any path $X \ra \cdots \ra X'=R^p(X)$ in $\bcG^k_{k}$ with $X,X'\in \bcE^0_{k}$. One has to prove that $f_{k}^p$ and $h^p$ coincide on $K_{X}\times C$. The rectangles $X$ and $X'$ are in a same connected component of $\bcG_{k}^{k-1}$ as some rectangles $X_{i}$ and $X_{j}$ respectively; hence (by ``transitivity'' of the required property) it is sufficient to show the property in the two following particular cases:
\begin{itemize}
\item[--] $X=X_{i}$ and $X'=X_{i+1}$. Here the property is satisfied by construction (see step 1).
\item[--] $X$ and $X'$ belong to a same connected component $\Gamma_{i}$
of $\bcG_{k}^{k-1}$. If $k=1$ there is nothing to prove\footnote{The graph $\bcG_{0}^0$ is assumed to have no edge; by compatibility (hypothesis $\mathbf{A_{1.c}}$), this also holds for $\bcG_{1}^0$. Thus  when $k=1$  one has $X=X'$.}, so we assume $k \geq 2$. Hypotheses $\mathbf{B_{2}}(k)$ and $\mathbf{C_{2}}(k)$ imply $f_{k}^p=f_{k-1}^p$ on $\Phi_{k}^{-1}(X)$ (apply remark~\ref{r.gngnp1}) which contains $K_X\times C$ (apply hypothesis $\mathbf{C_{5}}(k)$ \emph{via} item 1 of proposition~\ref{p.H-rectangle} and hypothesis $\mathbf{B_{5}}(k)$). Hence the induction hypothesis $\mathbf{C_{6}}(k-1)$ gives the property.
\end{itemize}
\end{proof}

\subsection{First part of the proof of theorem~\ref{t.main}, and minimal homeomorphisms with positive topological entropy}
\label{ss.bilan}

Under the hypotheses of  theorem~\ref{t.main}, we are now able to construct a homeomorphism $f$ that is topologically semi-conjugate to $R$ and has a subsystem isomorphic to $h$. For some appropriate choice of maps $R$ and $h$, this will lead to a minimal homeomorphisms with positive topological entropy on the torus $\TT^d$ (Rees theorem). 

\paragraph{Construction of a homeomorphism $f$ with a subsystem isomorphic to $h$}
Assume that the hypotheses of theorem~\ref{t.main} are satisfied:
\begin{itemize}
\item[--] $R$ is a homeomorphism on a manifold $\cM$ of dimension $d\geq 2$;
\item[-- ] $\mu$ is an aperiodic ergodic measure for $R$, and $A\subset\cM$ is a measurable set which has positive measure for $\mu$ and has zero-measure for every other ergodic $R$-invariant measure; 
\item[--]  $h$ is a bijective map which is fibered over $R$ (see section~\ref{ss.more-statements} for the precise definition).
\end{itemize}

Note that we can assume that $A$ is included in the support of $\mu$. We first apply the results of section~\ref{s.rectangles}: there exists a Cantor set $K$ included in $A$ (and thus also in $\supp(\mu)$) which is dynamically meagre and dynamically coherent, and a family $(\cE^0_{n})$ of rectangles satisfying hypotheses $\mathbf{A_{1,2,3}}$ (proposition~\ref{p.exist-rectangles-1} et~\ref{p.exist-rectangles-2}). According to remark~\ref{r.single-rectangle}, we can also assume that the graph $\cG_{0}^0$ has no edge, and even that $\cE_{0}^0$ is reduced to a single rectangle (since $\mu$ is ergodic).
We now apply the results of section~\ref{s.Cantor-sets-in-the-fibres} to obtain a sequence $(M_{n})$ satisfying hypotheses $\mathbf{B_{1,2,4,5,6}}$ (see proposition~\ref{p.existence-M}).
Then we apply alternatively propositions~\ref{p.existence-H-1} and~\ref{p.existence-H-2}
to get the sequences $(n_{k})$ and $(H_{k})$ such that hypotheses $\mathbf{C_{1,\dots,6}}$ are satisfied. Proposition~\ref{p.extension} (and~\ref{p.extension-bis}) ensures the existence of a map $\Phi$ and a homeomorphism $f$ such that $\Phi f= R \Phi$. Item~2 of proposition~\ref{p.extension2} implies that $\Phi$ is one-to-one outside the set $\Phi^{-1}(\supp(\mu))$. Since hypothesis $\mathbf{C_{4}}$ is satisfied, proposition~\ref{p.recurrence} shows that, if $R$ is transitive (resp. minimal), then the homeomorphism $f$ is also transitive (resp. minimal). And according to proposition~\ref{p.conjugacy-bis}, the dynamics of $f$ on $\displaystyle \bigcup_{i \in\ZZ}f^i(K\times C)$ is isomorphic to the dynamics of $h$: thus $f$ has a subsystem isomorphic to $h$. 

In subsection~\ref{ss.proof} we will complete the proof of theorem~\ref{t.main}, describing precisely the measurable dynamics of $f$ from the dynamics of $R$ and $h$. The proof of addendum~\ref{t.main} will be given in appendix~\ref{a.recurrence}.

\paragraph{Minimal homeomorphisms with positive topological entropy.}
Let us specify the above construction to the case where $h$ is a product map $R \times h_{C}$ and $h_{C}$ is a homeomorphism of the Cantor set $C$ admitting an invariant measure $\nu$ with positive entropy  (for example, $h_{C}$ is conjugate to the shift map and $\nu$ is a Bernoulli measure). 
Then the measure $\mu \times \nu$ has positive entropy for $h$. Since the homeomorphism $f$ is isomorphic to $h$ on some subset, this implies that $f$ has a measure with positive entropy. Using the variational principle, this shows that $f$ has positive topological entropy.  We have obtained the following result.
\begin{theo}
Every manifold $\cM$ of dimension $d\geq 2$ which carries a minimal homeomorphism also carries a  minimal homeomorphism with positive topological entropy.
\end{theo}
In particular, if $\cM=\TT^d$ and $R$ is an irrational rotation, we get Rees  theorem: 
\emph{there exists a minimal homeomorphism $f$ on $\TT^d$ with positive topological entropy.}

\section{Suppressing the dynamics outside the Cantor set $K \times C$}
\label{s.control}

We assume here that hypotheses $\mathbf{A_{1,2,3}}$, $\mathbf{B_{1,2,4,5,6}}$ are satisfied and that the graph $\cG_0^0$ has no edge. We have also fixed a fibered map $h$ on the set $\bigcup_{i\in \ZZ}R^i(K)\times C$. In the previous section, we have constructed $f$ so that its restriction to the set $\bigcup_{i\in \ZZ}f^i(K\times C)$ is isomorphic to $h$. Getting the isomorphism required by theorem~\ref{t.main} now amounts to  tighten the construction so that the set $\bigcup_{i\in\ZZ} f^i\left(\Phi^{-1}(K)\setminus K\times C\right)$ has measure $0$ for every $f$-invariant probability measure.

This will be obtained by requiring two additional hypotheses $\mathbf{C_{7,8}}$ on the sequences $(n_k)$ and $(H_k)$, under which  the following fact holds:  if we denote by $F$ the first return map of $f$ in $\Phi^{-1}(K)$,  every $F$-orbit will accumulate only on $K\times C$.  For this purpose, we will consider a sequence of ``waste bins" $(P_i)$  in $\cM$. We loosely describe  the hypotheses in terms of ``waste collection and management''. Consider a point $\wt x$ in $\Phi^{-1}(K)\setminus (K\times C)$ that returns infinitely many times in $\Phi^{-1}(K)$ in the future. Then
\begin{itemize}
\item[--] the ``waste collection'' (hypothesis~$\mathbf{C_{7}}$)
will ensure that some positive iterate $F^r(\wt x)$ will fall into some waste bin $P_{i}$;
\item[--] the ``waste management'' (hypothesis~$\mathbf{C_{8}}$) will ensure that, given that $F^r(\wt x)$ belongs to a waste bin, the forward $F$-orbit of $F^r(\wt x)$  will accumulate only on $K \times C$.
\end{itemize}

In order to obtain these hypotheses we will explain how to modify the construction of the $H_k$'s (the second part for the proof of the existence of $(n_k)$ and $(H_k)$). One important point is that the previous hypotheses $\mathbf{C_{5,6}}$ only dealt with the dynamics of $f$  on the Cantor set $K \times C$, whereas the new ones $\mathbf{C_{7,8}}$ will only deal with the dynamics in the complement of this Cantor set. Therefore the new constraints on the $H_{k}$'s will be compatible with the previous ones.

\subsection{The waste bins $(P_i)_{i\in\NN}$ and the neighbourhoods $(V_k)_{k\in\NN}$}
\label{ss.def-waste-bins}

We choose a sequence of topological closed balls $(P_i)_{i\in\NN}$ tamely embedded\footnote{\label{f.tamely-embedded} That is, each $P_{i}$ is the image of the unit ball of $\RR^d$ under a continuous one-to-one map from $\RR^d$ into $\wt X_{0}$.}  in the interior of the rectangle $\wt X_0 =X_0$, that will play the role of ``waste bins". We suppose that:  
\begin{itemize}
\item[--] the $P_i$'s are pairwise disjoint, and disjoint from the Cantor set $K\times C$, 
\item[--] $\displaystyle\mathop{\limsup}_{i\rightarrow\infty}P_i=K \times C$.
\end{itemize}

We also introduce a decreasing sequence $(V_k)_{k\in \NN}$ of neighbourhoods of the Cantor set $K \times C\subset \cM$ such that   
$$
\bigcap_{k\in\NN} V_k=K \times C.
$$
The neighbourhoods $(V_k)_{k\in \NN}$ will be used to perform the ``waste collection". Roughly speaking, at the $k^{th}$ step of the construction, we will make sure that the orbit of any point that is not in $V_k$ fall in some waste bin after some time.

\subsection{Hypotheses $\mathbf{C_{7,8}}$}
\label{ss.hyp-C_78}

\begin{figure}[htbp]
\input{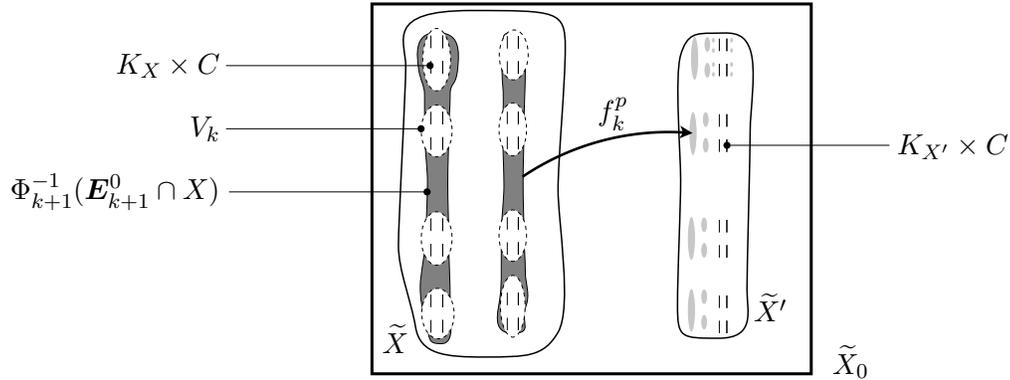}
\caption{Waste collection (hypothesis $\mathbf{C_{7}}(k)$). The waste to be collected are shaded in dark grey (in $\wt X$), the waste bins are in light grey (in $\wt X'$).}
\label{f.waste-collection}
\end{figure}

We consider the following hypotheses: 
\begin{itemize}
\item[]
\begin{itemize}
\item[$\mathbf{C_{7}}(k)$]  \emph{(Waste collection, see figure~\ref{f.waste-collection})}\\
Let $X \ra \cdots \ra X'=R^p(X)$ be a path in the graph $\bcG_k^k$ with $X,X'\in \bcE^0_{k}$,
which is not a path in the graph $\bcG_k^{k-1}$.
Then 
$$f_k^p\left(\Phi^{-1}_{k+1}(\bE^0_{k+1} \cap X)\setminus V_k\right)\subset \bigcup_{i \in \NN }P_i.$$

\item[$\mathbf{C_{8}}(k)$] \emph{(Waste management)}\\ 
Let $X \ra \cdots \ra X'=R^p(X)$ be a path in the graph $\bcG_k^k$ with $X,X'\in \bcE^0_{k}$. Then
for every $i\in\NN$,
$$f_k^p(\Phi^{-1}_{k+1}(\bE^0_{k+1} \cap X) \cap P_i)\subset \bigcup_{i'>i} P_{i'}.$$
\end{itemize} 
\end{itemize}

\begin{rema*}
If $X$ is an element of $\bcE_k^0$, then the set $\bE^0_{k+1} \cap X$ is the union of the ``immediate sub-rectangles" of $X$, and the set $\Phi^{-1}_{k+1}(\bE^0_{k+1} \cap X)$ is a (small) neighbourhood of $K_{X} \times C$.
\end{rema*}

\begin{rema}
\label{r.decalage-generationnel}
Remember that, because of hypothesis $\mathbf{C_{1}}$ on the support of $H_{k+1}$, we have
$$\Phi_{k+1}^{-1}(\bE_{k+1}^0)=(\bM_{k+1}\circ\Phi_k)^{-1}(\bE_{k+1}^{0}).$$
Hence hypothesis $\mathbf{C_7}(k)$ involves the integers 
$n_1,\dots,n_k$ and the homeomorphisms $H_1,\dots,H_k$ (in order to know the maps $\Phi_{k}$ and $f_k$), \emph{but it also involves the integer $n_{k+1}$}
(in order to know the set $\Phi_{k+1}^{-1}(\bE_{k+1}^0)$). It does not involve $H_{k+1}$.
The same remark holds for hypothesis $\mathbf{C_8}(k)$.
\end{rema}

\subsection{Consequences of hypotheses $\mathbf{C_{7,8}}$}
\label{ss.cons-C_78}

Let us assume that the sequences $(n_k)_{k\in\NN}$ and $(H_k)_{k\in\NN}$ were constructed in such a way that hypotheses $\mathbf{C_{1,2,3,5}}$ are satisfied.
Then, proposition~\ref{p.extension-bis} provides us with a map $\Phi:\cM\rightarrow\cM$ and a homeomorphism $f:\cM\rightarrow\cM$.
We denote by $\wt K$ the set $\Phi^{-1}(K)$ and by $F$ the first return map in $\wt K$ of $f$.
Let us recall that we have embedded the Cantor set $K\times C$ in $\cM$, that  $K\times C\subset \wt K$  and that   the Cantor set $K\times C$ is invariant under $F$ and $F^{-1}$:
if $\wt x\in\wt K$ has a forward iterate in $\wt K$, properties $\wt x\in K\times C$ and $F(\wt x)\in K\times C$ are equivalent (this is the content item 2 of proposition~\ref{p.H-rectangle} with the help of remark~\ref{r.gngnp1}).

\begin{prop}
\label{p.omega-limit}
Assume hypotheses $\mathbf{C_{1,2,3,5,7,8}}$ and consider a point $\wt x\in\wt K$ whose forward $f$-orbit returns infinitely many times in $\wt K$.
Then the $\omega$-limit set of the $F$-orbit of $\wt x$ is included in the Cantor set $K \times C$. 
\end{prop} 

\begin{lemma}
\label{l.waste-management}
Assume hypotheses $\mathbf{C_{1,2,3,5,8}}$ and consider a point $\wt y\in\wt K$ having a forward iterate by $f$ in $\wt K$.
If $\wt y$ is in a waste bin $P_i$, then $F(\wt y)$ is in a waste bin $P_{i'}$ with $i'>i$. 
\end{lemma}

\begin{proof}
Let $p$ be the return time of $\wt y$ in $\wt K$ (\emph{i.e.} the
integer such that $F(\wt y)=f^p(\wt y)$). Since $\wt K=\Phi^{-1}(K)=\Phi^{-1}(\bigcap_{k\in\NN} \bE_k^0)$,
the point $\wt y$ belongs to  $\wt X=\Phi^{-1}(X)$ for some rectangle $X\in\bcE_k^0$ and some $k>p$.
Since $k>p$, the rectangle $X'=R^p(X)$ is in $\bcE_k^k$, and belongs to the same connected component of the graph $\bcG_k^k$ as $X$.
Moreover, $X'$ must be also  an element of $\bcE_k^0$ since it contains the point $R^p(\Phi(\wt y))\in K$ (by iterability, hypothesis~$\mathbf{A_{1.a}}$).
Thus we can apply hypothesis $\mathbf{C_8}$. Note that $\wt y$ also belongs to the set 
$\Phi^{-1}(\bE_{k+1}^0 \cap X)$, which equals  $\Phi_{k+1}^{-1}(\bE_{k+1}^0 \cap X)$
 (by hypotheses $\mathbf{B_{1}}$ and $\mathbf{C_{1}}$ concerning the support), and that $f^p(\wt y ) = f_k^p(\wt y)$ (by proposition~\ref{p.extension-bis}).
Thus hypothesis~$\mathbf{C_8}$ implies that, if $\wt y$ is in a waste bin $P_i$, then
the point $F(\wt y)$ is in a waste bin $P_{i'}$ with $i'>i$.  
\end{proof}

\begin{proof}[Proof of proposition~\ref{p.omega-limit}]
If the point $\wt x$ is in $K\times C\subset \wt K$, then the proposition follows
immediately from the invariance of $K\times C$ under $F$.
Consider now the case where $\wt x$ is not in $K\times C$. Since $\bigcap_{k\in\NN}
V_k=K\times C$, there exists $k_0$ such that $\wt x\notin V_{k_0}$. Moreover,
since $(V_k)_{k\in\NN}$ is a \emph{decreasing} sequence of neighbourhoods,  
$\wt x$ is outside $V_k$ for every $k\geq k_0$. 

The point $\wt x$ is in $\wt K=\bigcap_{k\in\NN} \Phi^{-1}(\bE_k^0)$.
As a consequence, for every $k\in\NN$, the point $\wt x$ is in  $\wt X_k=\Phi^{-1}(X_k)$ for some rectangle $X_k\in\bcE_{k}^0$.
Let $r$ be a positive integer, and let $p$ be the integer such that $F^{r}(\wt x)=f^{p}(\wt x)$.
There exists a unique integer $k$ such that the rectangle $R^{p}(X_k)$ is in the same connected component of the graph of $\bcG_{k}^k$ as $X_k$
but not in the same connected component of the graph $\bcG_{k}^{k-1}$. 
Up to replacing $r$ by a bigger integer, we may assume that $k$ is bigger than $k_0$ (since there is no cycle in $\bcG_{k_{0}}^{k_{0}}$, hypothesis~$\mathbf{A_{1.a}}$). 
Then hypothesis~$\mathbf{C_7}(k)$ implies that the point $F^{r}(\wt x)$ is in the waste bin $P_{i_0}$ for some $i_0\in\NN$. 

Now, using lemma~\ref{l.waste-management} recursively, we obtain that $F^{r+s}(\wt x)$ is in a waste bin $P_{i_s}$, where $(i_s)_{s\in\NN}$ is an increasing sequence of integers. As a consequence, the $\omega$-limit set of the orbit of $\wt x$ under $F$ is included in $\displaystyle\mathop{\limsup}_{i\rightarrow\infty} P_i=K\times C$.  
\end{proof}

\begin{coro}
\label{c.full-measure}
Assume hypotheses $\mathbf{C_{1,2,3,5,7,8}}$ and let $\nu$ be an $f$-invariant probability measure  on $\cM$. Then
$$
\nu\left(\bigcup_{j\in\ZZ} f^j\left(\widetilde K\setminus(K\times C)\right)\right)=0.
$$
In particular, $f$ is isomorphic to the disjoint union
$$
\left(\bigcup_{j\in\ZZ} R^j(K)\times C,h\right) \bigsqcup \left(\cM\setminus\bigcup_{j\in\ZZ} R^j(K),R\right).
$$
\end{coro}

\begin{proof}
The restriction of $\nu$ to $\widetilde K$ is $F$-invariant. So, by proposition~\ref{p.omega-limit} and using Poincar\'e recurrence theorem, we have $\nu\left(\widetilde K\setminus (K\times C)\right)=0$. The first part of the corollary follows.

For the second part, one defines the bi-measurable map
$$
\Theta\colon \left(\cM\setminus \bigcup_{j\in \ZZ}f^j\left(\wt K\right)\right)
\bigsqcup \left(\bigcup_{j\in \ZZ} f^j(K\times C)\right) 
\to
\left(\cM\setminus\bigcup_{j\in\ZZ} R^j(K)\right)     \bigsqcup   \left(\bigcup_{j\in\ZZ} R^j(K)\times C\right)
 $$
given by $\Phi$ on the set $\displaystyle\cM\setminus \bigcup_{j\in \ZZ}f^j\left(\wt K\right)$ and by proposition~\ref{p.conjugacy-bis} on the set $\displaystyle\bigcup_{j\in \ZZ} f^j(K\times C)$.
By the first part of the corollary, the set 
$$
\displaystyle\cM_0:=\left(\cM\setminus \bigcup_{j\in \ZZ}f^j\left(\wt K\right)\right)\bigsqcup \left(\bigcup_{j\in \ZZ} f^j(K\times C)\right)
$$
 is full in $\cM$: for any invariant measure $\nu$ we have $\nu(\cM\setminus \cM_0)=0$. By propositions~\ref{p.extension-bis} and~\ref{p.extension2},  $\Phi$ is one-to-one on   the set $\displaystyle\cM\setminus \bigcup_{j\in \ZZ}f^j\left(\wt K\right)$, so that $\Theta$ is one-to-one. By propositions~\ref{p.extension}  and~\ref{p.conjugacy-bis}, $\Theta$ is a conjugacy. The second part of the corollary follows. 
\end{proof}

\subsection{Realisation of hypotheses $\mathbf{C_7}$ and $\mathbf{C_8}$}
\label{ss.control-existence}

We are left to prove the following proposition.

\begin{prop}
\label{p.existence-H-n}
There exists a sequence of integers $(n_k)_{k\in\NN}$ and a sequence of homeomorphisms $(H_k)_{k\in\NN}$ such that hypotheses $\mathbf{C_{1,\dots,8}}$ are satisfied. 
\end {prop}

\begin{proof}
We proceed by induction. We consider an integer $k_0$. We assume that the integers $n_{0},n_1,\dots,n_{k_{0}}$ and the homeomorphisms $H_{0}, \dots,H_{k_{0}-1}$ are constructed, and that hypotheses $\mathbf{C_{1,\dots ,8}}(k)$,  are satisfied for every $k\leq k_0-1$ (see remark~\ref{r.decalage-generationnel}). We will explain how to construct an integer $n_{k_0+1}$ and a homeomorphism $H_{k_0}$ such that hypotheses $\mathbf{C_{1, \dots ,8}}(k_{0})$ are satisfied.

\paragraph{Some explanations on the proof.}
As explained at the beginning of the section, we will first use proposition~\ref{p.existence-H-2} to get a homeomorphism  $H_{k_0}$ such that hypotheses $\mathbf{C_{1,2,5,6}}(k_{0})$ are satisfied, and then we will modify this homeomorphism outside the Cantor set $\bPsi_{k_0}(K\times C)$ in order to get hypotheses $\mathbf{C_{7,8}}(k_0)$.

One important difficulty is due to the interplay between the choices of the homeomorphism $H_{k_0}$ and the integer $n_{k_0+1}$: on the one hand, there are many reasons why we have to choose the integer $n_{k_0+1}$ after the homeomorphism $H_{k_0}$ (e.g. to get hypothesis $\mathbf{C_{3,4}}(k_0)$); on the other hand, it seems that we need to know the integer $n_{k_0+1}$ when we construct the homeomorphism $H_{k_0}$ in order to get hypothesis $\mathbf{C_{7}}(k_0)$ (see remark~\ref{r.decalage-generationnel}). To solve this problem, we will consider the set:
$$
G =   \bigcap_{i \geq 1} \left( M_{n_{k_{0}}+i} \circ \cdots \circ  M_{n_{k_{0}}+1}\right)^{-1}
\left(E^0_{n_{k_{0}}+i}\right).
$$
Note that, if we imagine (for sake of simplicity) that the sequence $(M_n)_{n\geq 1}$ satisfies hypothesis $\mathbf{B_3}$ (so that the sequence" $(\Psi_n)_{n\geq 1}$ converges towards a map $\Psi$), then we have 
$$
G =\Psi_{n_{k_{0}}} \left( \Psi^{-1}(K)\right).
$$
According to hypothesis~$\mathbf{B_{4}}$, the set $G$ has empty interior.
We will construct the homeomorphism ${H_{k_{0}}}$ and a (thin) neighbourhood $W$ of $\Phi_{k_{0}}^{-1}(G)$, in such a way that the set $W \setminus V_{k_{0}}$ will be mapped in a waste bin by the appropriate powers of $f_{k_0}$: more precisely,  hypothesis~$\mathbf{C_{7}}(k_{0})$ will hold with  the set $\Phi^{-1}_{k+1}(\bE^0_{k+1})$ replaced by $W$. Then we will choose the integer  $n_{k_{0}+1}$ big enough, so that the set 
$$
\Phi_{k_{0}+1}^{-1}\left(\bE^0_{k_{0}+1}\right)=\Phi_{k_{0}}^{-1}\circ \bM_{k_0+1}^{-1}\left(\bE_{k_0+1}^0\right)=\Phi_{k_{0}}^{-1}\circ \left(M_{n_{k_{0}+1}} \circ \cdots \circ  M_{n_{k_{0}}+1}\right)^{-1} \left(E^0_{n_{k_{0}+1}}\right)
$$
will be included in $W$. Thus we will get hypothesis $\mathbf{C_{7}}(k_0)$.
The same strategy will work to get hypothesis $\mathbf{C_{8}}$.

\paragraph{Some notations.}
For any rectangle $X \in \bcE_{k}^k$, the rectangle $\Phi_{k}^{-1}(X)$ will be denoted by $\wt X$.
Similarly, we set $\wt \bE_{k}^0 = \Phi_{k}^{-1}(\bE_{k}^0)$.

The different connected component $\Gamma$ of the graph $\bcG_{k_0}^{k_0}$ will be considered independently. Let $\Gamma$ be a connected component of the graph $\bcG_{k_0}^{k_0}$.  We denote by $\Gamma_1,\dots,\Gamma_s$ the connected components of the graph $\bcG_{k_0}^{k_0-1}$ included in $\Gamma$ (ordered  by the orientation of $\Gamma$). For $j=1,\dots,s$, we denote by $X_j^\mathrm{in}$ the first vertex of $\Gamma_j$ which is in $\cE_{k_0}^0$, and by $X_j^\mathrm{out}$ the last vertex of $\Gamma_j$ which is in $\cE_{k_0}^0$. Finally, for $j=1,\dots,s-1$, we denote by $q_j$ the positive integer such that $X_{j+1}^\mathrm{in}=R^{q_j}(X_j^\mathrm{out})$. See figure~\ref{f.graphe}.

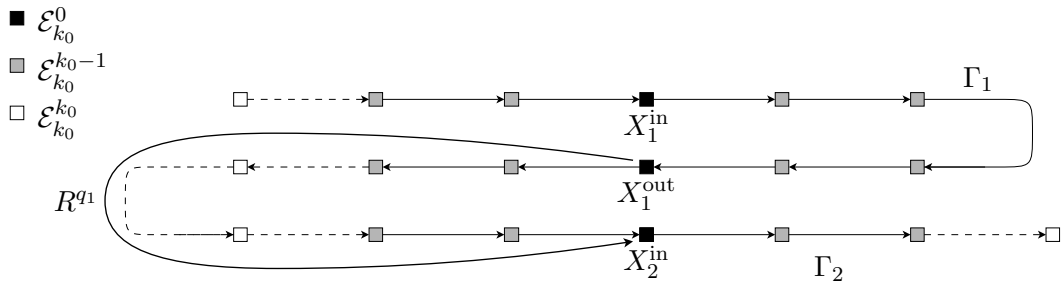
\begin{figure}[ht]
\def\JPicScale{0.9}
\ifx\JPicScale\undefined\def\JPicScale{1}\fi
\psset{unit=\JPicScale mm}
\psset{linewidth=0.3,dotsep=1,hatchwidth=0.3,hatchsep=1.5,shadowsize=1}
\psset{dotsize=0.7 2.5,dotscale=1 1,fillcolor=black}
\psset{arrowsize=1 2,arrowlength=1,arrowinset=0.25,tbarsize=0.7 5,bracketlength=0.15,rbracketlength=0.15}
\begin{pspicture}(0,0)(161,53)
\pspolygon[linewidth=0.1](39,21)(41,21)(41,19)(39,19)
\psline[linewidth=0.1]{->}(119,30)(101,30)
\psline[linewidth=0.1]{->}(99,30)(81,30)
\psline[linewidth=0.1,linestyle=dashed,dash=1 1]{->}(59,30)(41,30)
\psline[linewidth=0.1,linestyle=dashed,dash=1 1]{->}(41,20)(59,20)
\pspolygon[linewidth=0.1,fillstyle=solid](99,41)(101,41)(101,39)(99,39)
\pspolygon[linewidth=0.1,fillstyle=solid](99,31)(101,31)(101,29)(99,29)
\newrgbcolor{userFillColour}{0.71 0.71 0.71}
\pspolygon[linewidth=0.1,fillcolor=userFillColour,fillstyle=solid](139,31)(141,31)(141,29)(139,29)
\newrgbcolor{userFillColour}{0.71 0.71 0.71}
\pspolygon[linewidth=0.1,fillcolor=userFillColour,fillstyle=solid](139,41)(141,41)(141,39)(139,39)
\newrgbcolor{userFillColour}{0.71 0.71 0.71}
\pspolygon[linewidth=0.1,fillcolor=userFillColour,fillstyle=solid](59,41)(61,41)(61,39)(59,39)
\pspolygon[linewidth=0.1,fillstyle=solid](99,21)(101,21)(101,19)(99,19)
\pscustom[linewidth=0.1]{\psline(141,40)(152,40)
\psbezier(152,40)(157,40)(157,40)(157,35)
\psbezier(157,35)(157,30)(157,30)(152,30)
\psbezier(152,30)(147,30)(141,30)(141,30)
}
\psline[linewidth=0.1]{->}(150,30)(141,30)
\newrgbcolor{userFillColour}{0.71 0.71 0.71}
\pspolygon[linewidth=0.1,fillcolor=userFillColour,fillstyle=solid](119,41)(121,41)(121,39)(119,39)
\newrgbcolor{userFillColour}{0.71 0.71 0.71}
\pspolygon[linewidth=0.1,fillcolor=userFillColour,fillstyle=solid](119,31)(121,31)(121,29)(119,29)
\psline[linewidth=0.1]{->}(139,30)(121,30)
\newrgbcolor{userFillColour}{0.71 0.71 0.71}
\pspolygon[linewidth=0.1,fillcolor=userFillColour,fillstyle=solid](79,31)(81,31)(81,29)(79,29)
\psline[linewidth=0.1]{->}(79,30)(61,30)
\newrgbcolor{userFillColour}{0.71 0.71 0.71}
\pspolygon[linewidth=0.1,fillcolor=userFillColour,fillstyle=solid](59,31)(61,31)(61,29)(59,29)
\psline[linewidth=0.1]{->}(61,40)(79,40)
\newrgbcolor{userFillColour}{0.71 0.71 0.71}
\pspolygon[linewidth=0.1,fillcolor=userFillColour,fillstyle=solid](81.1,41)(79,41)(79,39)(81.1,39)
\psline[linewidth=0.1]{->}(81,40)(99,40)
\psline[linewidth=0.1]{->}(121,40)(139,40)
\psline[linewidth=0.1]{->}(101,40)(119,40)
\pspolygon[linewidth=0.1,fillstyle=solid](99,21)(101,21)(101,19)(99,19)
\newrgbcolor{userFillColour}{0.71 0.71 0.71}
\pspolygon[linewidth=0.1,fillcolor=userFillColour,fillstyle=solid](139,21)(141,21)(141,19)(139,19)
\newrgbcolor{userFillColour}{0.71 0.71 0.71}
\pspolygon[linewidth=0.1,fillcolor=userFillColour,fillstyle=solid](59,21)(61,21)(61,19)(59,19)
\newrgbcolor{userFillColour}{0.71 0.71 0.71}
\pspolygon[linewidth=0.1,fillcolor=userFillColour,fillstyle=solid](119,21)(121,21)(121,19)(119,19)
\psline[linewidth=0.1]{->}(61,20)(79,20)
\newrgbcolor{userFillColour}{0.71 0.71 0.71}
\pspolygon[linewidth=0.1,fillcolor=userFillColour,fillstyle=solid](81.1,21)(79,21)(79,19)(81.1,19)
\psline[linewidth=0.1]{->}(81,20)(99,20)
\psline[linewidth=0.1]{->}(121,20)(139,20)
\psline[linewidth=0.1]{->}(101,20)(119,20)
\pspolygon[linewidth=0.1](39,31)(41,31)(41,29)(39,29)
\pscustom[linewidth=0.1,linestyle=dashed,dash=1 1]{\psline(39,30)(27,30)
\psbezier(27,30)(23,30)(23,30)(23,25)
\psbezier(23,25)(23,20)(23,20)(27,20)
\psbezier(27,20)(31,20)(39,20)(39,20)
}
\psline[linewidth=0.1,linestyle=dashed,dash=1 1]{->}(31,20)(39,20)
\pspolygon[linewidth=0.1](39,41)(41,41)(41,39)(39,39)
\psline[linewidth=0.1,linestyle=dashed,dash=1 1]{->}(41,40)(59,40)
\psline[linewidth=0.1,linestyle=dashed,dash=1 1]{->}(141,20)(159,20)
\pspolygon[linewidth=0.1](159,21)(161,21)(161,19)(159,19)
\newrgbcolor{userFillColour}{0.86 0.86 0.86}
\rput(149,43){$\Gamma_1$}
\newrgbcolor{userFillColour}{0.86 0.86 0.86}
\rput(127,15){$\Gamma_2$}
\newrgbcolor{userFillColour}{0.86 0.86 0.86}
\rput(100,16){$X_2^\mathrm{in}$}
\newrgbcolor{userFillColour}{0.86 0.86 0.86}
\rput(100,36){$X_1^\mathrm{in}$}
\newrgbcolor{userFillColour}{0.86 0.86 0.86}
\rput(100,26){$X_1^\mathrm{out}$}
\pscustom[linewidth=0.2]{\psbezier{-}(98,31)(78.4,33.8)(62.5,35)(45,35)
\psbezier(45,35)(27.5,35)(20,32)(20,25)
\psbezier(20,25)(20,18)(27.65,15)(45.5,15)
\psbezier{->}(45.5,15)(63.35,15)(79.1,16.2)(98,19)
}
\rput[r](19,25){$R^{q_1}$}
\newrgbcolor{userFillColour}{0.86 0.86 0.86}
\rput(10,44){}
\newrgbcolor{userFillColour}{0.86 0.86 0.86}
\rput[l](10,44){$\cE_{k_0}^{k_0-1}$}
\newrgbcolor{userFillColour}{0.86 0.86 0.86}
\rput[l](10,37){$\cE_{k_0}^{k_0}$}
\newrgbcolor{userFillColour}{0.86 0.86 0.86}
\rput[l](10,51){$\cE_{k_0}^0$}
\pspolygon[linewidth=0.1,fillstyle=solid](6,53)(8,53)(8,51)(6,51)
\pspolygon[linewidth=0.1](6,39)(8,39)(8,37)(6,37)
\newrgbcolor{userFillColour}{0.71 0.71 0.71}
\pspolygon[linewidth=0.1,fillcolor=userFillColour,fillstyle=solid](6,46)(8,46)(8,44)(6,44)
\end{pspicture}
\caption{
The connected components $\Gamma_1,\dots,\Gamma_s$ in $\Gamma$, the rectangles $X_1^\mathrm{in},\dots,X_s^\mathrm{in}$ and $X_{1}^\mathrm{out},\dots,X_{s}^\mathrm{out}$, the integers  $q_1,\dots,q_{s-1}$}
\label{f.graphe}
\end{figure}

For every $j\leq s$, we will define the homeomorphism $H_{k_0}$ on the rectangle $X_j^\mathrm{in}$. Then hypotheses $\mathbf{C_{1,2}}(k_0)$ (support and commutation with $R$) will enforce the definition on  the remaining vertices of $\Gamma$.

\paragraph{Step 1. Construction of the restriction of $\bm{H_{k_0}}$ to $\bm{\Psi_{k_{0}}(K \times C)}$.} Using proposition~\ref{p.existence-H-2}, we can construct a homeomorphism $H_{k_{0}}^0$ such that hypotheses $\mathbf{C_{1,2,5,6}}(k_{0})$ are satisfied. We will modify this homeomorphism $H_{k_0}^0$ in order to get a homeomorphism $H_{k_0}$ such that hypotheses $\mathbf{C_{7,8}}(k_0)$ will also be satisfied. The new homeomorphism $H_{k_0}$ will coincide with $H_{k_0}^0$ on the set $\bPsi_{k_{0}}(K \times C).$ Note that hypotheses $\mathbf{C_{5,6}}(k_0)$ depend only on the restriction of $H_{k_0}$ to this set (see point 1 of proposition~\ref{p.H-rectangle}). Hence the replacement of $H_{k_0}^0$ by $H_{k_0}$ will not destroy hypotheses $\mathbf{C_{5,6}}(k_{0})$. We will use the notation 
$$
{f}_{k_{0}}^0:=\left(H_{k_0}^0\circ\bM_{k_0}\circ \Phi_{k_0-1}\right)^{-1} \circ R\circ \left(H_{k_0}^0\circ  \bM_{k_0}\circ \Phi_{k_0-1}\right).
$$
The map (yet to be constructed) 
$$
f_{k_{0}}=\left(H_{k_0}\circ  \bM_{k_0}\circ  \Phi_{k_0-1}\right)^{-1} \circ R\circ  \left(H_{k_0}\circ  \bM_{k_0}\circ  \Phi_{k_0-1}\right)
$$ 
will be seen as an alteration of the map $f_{k_{0}}^0$.

\paragraph{Step 2. Definition of $\bm{H_{k_0}}$ on $\bm{\Gamma_1}$.} 
On the first component $\Gamma_{1}$, we do not modify the homeomorphism $H_{k_0}^0$ given by proposition~\ref{p.existence-H-2} : we set $H_{k_0}:=H_{k_0}^0$.
From now on, the formulae
$$
\Phi_{k_{0}} :=  H_{k_{0}} \circ \bM_{k_{0}} \Phi_{k_{0}-1} \mbox{ and } 
f_{k_{0}} := \Phi_{k_{0}}^{-1} \circ R \circ \Phi_{k_{0}}
$$
define the map $\Phi_{k_{0}}$ above the vertex of $\Gamma_{1}$ and the map
$f_{k_{0}}$ along the edges of the graph $\Phi_{k_{0}}^{-1}(\Gamma_{1})$.

\paragraph{Step 3. Definition of $\bm{H_{k_0}}$ on $\bm{\Gamma_2}$.}
We will now define ${H_{k_0}}$ on $X_2^\mathrm{in}$, in such a way that the map
${f}_{k_0}^{q_1}:\wt X_1^\mathrm{out} \to \wt X_2^\mathrm{in}$
will have the desired properties (i.e. will map the appropriate sets into a waste bin). 

Note that if $X=R^{-p}(X_1^\mathrm{out})$ is any vertex of $\Gamma_1$, then 
${f}_{k_0}^{p}(K_{X} \times C ) = K_{X_1^\mathrm{out}} \times C$ (point 2 of proposition~\ref{p.H-rectangle}). 
Now we choose a set  $U_1 \subset \wt X_{1}^\mathrm{out}$ which is a small neighbourhood of $K_{X_1^\mathrm{out}} \times C$ such that:
\begin{itemize}
\item[--] if $X=R^{-p}(X_1^\mathrm{out})$ is a vertex of $\Gamma_1$ in $\bcE_{k_0}^0$, then $U_1$ is contained in ${f}_{k_0}^{p}(V_{k_0}\cap\wt X)$;
\item[--] the boundary of $U_1$ is disjoined from all the waste bins;
\item[--] $U_{1}$ is a finite union of pairwise disjoint topological balls tamely embedded  in $\cM$ (this is possible since the Cantor set $K\times C$ is tamely embedded).
\end{itemize}
Let
$$
G =  \bigcap_{i \geq 1} \left( M_{n_{k_{0}}+i} \circ \cdots \circ  M_{n_{k_{0}}+1}\right)^{-1} \left(E^0_{n_{k_{0}}+i}\right).
$$ 
The set $\Phi_{k_{0}}^{-1} (G \cap X_1^\mathrm{out})$ is contained in $\wt X_1^\mathrm{out}$ and has empty interior (hypothesis $\mathbf{B_{4}}$). Moreover, from the very definition of $G$ we see that the set $\Phi_{k_{0}}^{-1} (G \cap X_1^\mathrm{out})$ can be obtained as a decreasing intersection of finite union of pairwise disjoint topological balls. Hence, we can find a neighbourhood $W_1$ of $\Phi_{k_{0}}^{-1} (G \cap X_1^\mathrm{out})$, such that $W_1$ is a finite union of pairwise disjoined  topological balls, tamely embedded in $\inte(\wt X_{1}^\mathrm{out})$, with arbitrarily small internal radii. In particular, we can assume that the internal radii of these balls is smaller than the infimum of the internal radii of the connected components of $U_1$, so that no connected component of $U_{1}$ is contained in $W_{1}$ (see figure~\ref{f.waste-collection} of subsection~\ref{ss.hyp-C_78}, where $U_{1}$ is replaced by $V_{k}$ and $W_{1}$ by $\Phi^{-1}_{k+1}(\bE^0_{k+1} \cap X)$).  Let $Z$ be the union of  $K_{X_1^\mathrm{out}} \times C$ and all the waste bins contained in $U_1$; thanks to the properties of $U_{1}$ and of the waste bins, $Z$ is a compact set inside the interior of $U_{1}$.
We will use the following lemma (since everything takes place inside the interior of the rectangle $\wt X_{1}^\mathrm{out}$, we can assume that the ambient space is $\RR^d$).

\begin{lemma}
Let $W_1$ and $U_1$ be finite union of pairwise disjoint closed topological balls tamely embedded in $\RR^d$. Assume that every connected component of $U_1$ meets $\RR^d \setminus W_1$. Let $Z$ be a compact set inside the interior of $U_{1}$.

Then  there exists a set $W'_{1}$ which is again a finite union of pairwise disjoint closed topological balls tamely embedded in $\RR^d$, which contains $W_1\setminus U_1$ and which does not meet $Z$.
\end{lemma}

\begin{proof}
Using a homeomorphism $L$ whose support is included in  $U_{1}$, one can push $Z$ outside $W_1$ (the map $L$ can be constructed Fibrely on each connected component of $U_{1}$). Then one takes  $W'_{1} = L^{-1}(W_1)$.
\end{proof}

Denote by $B_{1}, \dots , B_{r}$ the pairwise disjoint closed topological balls given by the previous lemma. Let  $\cB$ be the collection of these balls and of all the waste bins $P_i$ contained in $U_{1}$. 
According to the lemma, $\cB$ is a family of pairwise disjoint sets. Note that $\cB$ contains all but a finite number of waste bins $P_{i}$; in particular, each  ball $B_{1}, \dots , B_{r}$  meets only a finite number of the waste bins $P_i$. Thus for every $B \in \cB$, we can define an integer  $i(B)$ such that no waste bin $P_{i}$ with $i > i(B)$ meets $B$. Note that every ball  of $\cB$ is tamely embedded. Now this family $\cB$ contains everything that we want to dispose of: each ball $B$ will be put inside a waste bin $P_{i}$ with $i > i(B)$. This will be done by the way of the following extension lemma (note that this is the place where we use the inclusion $ K \times C \subset\displaystyle\mathop{\limsup}_{i\rightarrow\infty}P_i $).

\begin{fact}
\label{f.chirurgie}
There exists a homeomorphism $g : \wt X_1^\mathrm{out}\to \wt X_2^\mathrm{in}$  such that 
\begin{enumerate}
\item $g$ coincides with $\left({f}_{k_0}^0\right)^{q_{1}}$ on $K_{X_1^\mathrm{out}}\times C$ and on the boundary of  $\wt X_1^\mathrm{out}$~;
\item for every $B\in\cB$, there exists $i > i(B)$ such that $g(B) \subset P_{i}$.
\end{enumerate}
\end{fact}

\begin{proof}
This fact is contained in corollary~\ref{c.cantor-ball} in appendix~\ref{a.cantor}.
Here is the translation: the sets $X,X'$ of the corollary are $\wt X_1^\mathrm{out}, \wt X_2^\mathrm{in}$; $Q,Q'$ are $K_{X_1^\mathrm{out}}\times C, K_{X_2^\mathrm{in}}\times C$; $\alpha$ and $\beta$ are the restrictions of $\left({f}_{k_0}^0\right)^{q_{1}}$ respectively  on the boundary of  $\wt X_1^\mathrm{out}$ and on $K_{X_1^\mathrm{out}}\times C$; the sequence $(B_{j})$ is any indexation of the family $\cB$; the sequence $(B'_{i})$ is any indexation of the family of those waste bins $P_{i}$ that are included in the interior of $\wt X_2^\mathrm{in}$. Note that since $\displaystyle\mathop{\limsup}_{i\rightarrow\infty}P_i=K \times C$ and $Q' = K_{X_2^\mathrm{in}}\times C \subset \inte(\wt X_2^\mathrm{in})$, we have $\displaystyle\mathop{\limsup}_{i\rightarrow\infty}B'_{i}=Q'$.
\end{proof}

We want to define the homeomorphism ${H_{k_{0}}}$ on $X_2^\mathrm{in}$ in such a way that 
the homeomorphism ${f}_{k_0}^{q_1}:\wt X_1^\mathrm{out}\to\wt X_2^\mathrm{in}$  will coincide with the homeomorphism $g$ given by the above fact. For this purpose, recall that  ${f}_{k_0}^{q_1}$ will be such that
$$
{f}_{k_0}^{q_1}|_ {\wt X_1^\mathrm{out}}= ( {\Phi_{k_{0}-1}} \circ {\bM_{k_{0}}} \circ {H_{k_{0}}} |_{X_2^\mathrm{in}} ) ^{-1}   \circ R^{q_1} \circ  ( {\Phi_{k_{0}-1}} \circ {\bM_{k_{0}}} \circ {H_{k_{0}}} |_{X_1^\mathrm{out}} ).
$$
In the right-hand term of this equality, all the maps are already defined except ${H_{k_{0}}} |_{X_2^\mathrm{in}}$.
As a consequence, replacing ${f}_{k_0}^{q_1}|_ {\wt X_1^\mathrm{out}}$ by $g$ in the above equality,
we obtain a formula that tells us how to define the map ${H_{k_{0}}} |_{X_2^\mathrm{in}}$.
 Observe that:
\begin{itemize}
\item[--] The first property in  fact~\ref{f.chirurgie} implies that $H_{k_0} |_{X_2^\mathrm{in}}$ will coincide with $H_{k_0}^0$ on $\bPsi_{k_0}(K_{X_2^\mathrm{in}} \times C)$ and on the boundary of $X_2^\mathrm{in}$, and in particular it is the identity on this boundary. 
\item[--] The second  property in the fact implies that ${f}_{k_0}^{q_1}$ maps $W_1\setminus U_1$ into a union of waste bins $P_i$.
\item[--] The second property in the fact also implies that, for every $i$, the homeomorphism  ${f}_{k_0}^{q_1}$ maps $W_{1} \cap P_i$ into a union of waste bins $P_{i'}$ with $i'>i$ (this also makes use of the properties of the family $\cB$). 
\end{itemize}
The last two points can be reformulated as follows:  hypotheses $\mathbf{C_{7,8}}(k_{0})$ are now satisfied in the special case $X=X_1^\mathrm{out}, X' = X_2^\mathrm{in}$ with the set 
$\Phi^{-1}_{k_{0}+1}(\bE^0_{k_{0}+1} \cap X)$ replaced by the set $W_{1}$ and the set $V_{k_{0}}$ replaced by the set $U_1$.

The homeomorphism ${H_{k_0}}$ is now defined on the rectangle $X_2^\mathrm{in}$. Since 
${H_{k_0}}$ has to commute with $R$ along the edges of $\bcG_{k_0}^{k_0-1}$ (hypothesis $\mathbf{C_2}(k_0)$), this forces automatically the definition of ${H_{k_0}}$ on all the vertices of $\Gamma_2$.

\paragraph{Step 4. End of the definition of $\bm{H_{k_0}}$.} 
We repeat the same procedure as in step~3 to define ${H_{k_0}}$ on the vertices of $\Gamma_3,\dots,\Gamma_s$. Of course, we extend ${H_{k_0}}$ by the identity outside ${\bE_{k_0}^{k_0-1}}$ so that hypothesis $\mathbf{C_1}(k_0)$ is satisfied.

\paragraph{Step 5. (Temporary) choice of the integer $\bm{n_{k_{0}+1}}$.}
For every $j\leq s-1$, we have defined above (in step 3 for $j=1$, in step 4 for $j>1$) a neighbourhood $W_{j}$ of $\Phi_{k_0}^{-1}(G\cap X_j^\mathrm{out})$. By definition, $G$ is the decreasing intersection of the sets:
$$
\left( M_{n_{k_{0}}+i} \circ \cdots \circ  M_{n_{k_{0}}+1}\right)^{-1} \left(E^0_{n_{k_{0}}+i}\right)
$$
for $i\geq 1$. So, we can choose an integer $n_{k_{0}+1}$ big enough, so that (for every $j$) the set 
$$
\begin{array}{lll}
\Phi_{k_{0}+1}^{-1}\left(\bE^0_{k_{0}+1}\cap X_j^\mathrm{out}\right) & = & \Phi_{k_{0}}^{-1}\circ \bM_{k_0+1}^{-1}\left(\bE_{k_0+1}^0 \cap X_j^\mathrm{out}\right)\\
 & = & \Phi_{k_{0}}^{-1}\circ \left(M_{n_{k_{0}+1}} \circ \cdots \circ  M_{n_{k_{0}}+1}\right)^{-1} \left(E^0_{n_{k_{0}+1}}\cap X_j^\mathrm{out}\right)
 \end{array}
 $$
is included in $W_j$.

\paragraph{Hypothesis $\mathbf{C_{8}}(k_0)$ is satisfied.} 
Let $X \ra \cdots \ra X'=R^p(X)$ be a path in the graph $\bcG_k^k$ with $X,X'\in \bcE^0_{k}$.
Fix an integer $i$. We have to check that the homeomorphism $f_{k_{0}}^p$ maps the set $\wt \bE^0_{k_{0}+1} \cap \wt X \cap P_i=\Phi_{k_0+1}^{-1}(\bE^0_{k_{0}+1} \cap X) \cap P_i$ in a union of waste bins $P_{i'}$ with $i' > i$. This property behaves well under composition, so it is enough to prove it in the two following particular cases: 

\medskip

\noindent\textit{Case~1 (transition between two successive components of the graph $\bcG_{k_0}^{k_0-1}$): there exists an integer $j\leq s-1$ such that $X=X_j^\mathrm{out}$ and $X'=X_{j+1}^\mathrm{in}$.}  In this case,  the hypothesis follows from the construction of the homeomorphism $H_{k_0}$ (see the end of  step 3 for $j=1$) and the choice of the integer $n_{k_{0}+1}$.

\medskip

\noindent\textit{Case~2: the rectangles $X,X'$ are in the same connected component of the graph $\bcG_{k_0}^{k_0-1}$.}  In this case, the hypothesis follows easily from our induction hypothesis $\mathbf{C_{8}}(k_0-1)$, since ${f}_{k_0}^p$ coincides with ${f}_{k_0-1}^p$ on $\wt X$ (by commutation, hypotheses $\mathbf{B_{2}}$ and $\mathbf{C_{2}}$).
 
\paragraph{Hypothesis $\mathbf{C_7}(k_0)$ is satisfied}
Let $X \ra \cdots \ra X'=R^p(X)$ be a path in the graph $\bcG_k^k$ with $X,X'\in \bcE^0_{k}$,
which is not a path in the graph $\bcG_k^{k-1}$. We have to prove that the homeomorphism ${f}_{k_0}^p:\wt X\rightarrow\wt X'$ maps $(\wt \bE_{k_0+1}^0 \cap \wt{X} ) \setminus V_{k_0}=\Phi_{k_0+1}^{-1}(\bE^0_{k_{0}+1} \cap X)\setminus V_{k_0}$ into a union of waste bins $P_i$. We consider the integer $j$ such that $X$ is a vertex of $\Gamma_j$, and the non-negative integer $r$ such that $R^r(X)=X_j^\mathrm{out}$. We see ${f}_{k_0}^p:\wt X\rightarrow\wt X'$ as the composition of ${f}_{k_0}^r:\wt X\to\wt X_j^\mathrm{out}$, ${f}_{k_0}^{q_j}:\wt X_j^\mathrm{out} \to\wt X_{j+1}^\mathrm{in}$ and ${f}_{k_0}^{p-r-q_j}:\wt X_{j+1}^\mathrm{in}\to \wt X'$.
\begin{itemize}
\item[a.]  Let us first look at the homeomorphism ${f}_{k_0}^r:\wt X\to\wt X_j^\mathrm{out}$. Recall that the rectangles $\wt X$ and $\wt X_j^\mathrm{out}$ are in the same connected component of the graph $\bcG_{k_0}^{k_0-1}$. Hence by construction (see step 3 for the case $j=1$), the neighbourhood $U_j$ is contained in ${f}_{k_0}^r(\wt X \cap  V_{k_0})$. As a consequence, ${f}_{k_0}^r$ maps $\left(\wt \bE_{k_0+1}^0 \cap\wt X\right)\setminus V_{k_0}$ into  $\left(\wt \bE_{k_0+1}^0 \cap\wt X_j^\mathrm{out} \right)\setminus U_j$. 
\item[b.] The homeomorphism $H_{k_0}$ was constructed in such a way  that  ${f}_{k_0}^{q_j}:\wt X_j^\mathrm{out}\to\wt X_{j+1}^\mathrm{in}$
 maps $W_j\setminus U_j$ into a union of waste bins $P_i$ (see the end of step 3). And the integer $n_{k_0+1}$ was chosen in such a way that $\wt \bE_{k_0+1}^0\cap\wt X_j^\mathrm{out}$ is contained in $W_j$ (see step 5). As a consequence,  ${f}_{k_0}^{q_j}$ maps $\left(\wt \bE_{k_0+1}^0 \cap\wt X_j^\mathrm{out}\right)\setminus U_j$ into a set $\left(\wt \bE_{k_0+1}^0 \cap\wt X_{j+1}^\mathrm{in}\right) \cap S$ where $S$ is   a union of waste bins $P_i$.
 \item[c.] Finally, since we have already checked hypothesis $\mathbf{C_{8}}(k_0)$, we know that the homeomorphism ${f}_{k_0}^{p-r-q_j}:\wt X_{j+1}^\mathrm{in}\to \wt X'$ maps any
set $\left(\wt \bE_{k_0+1}^0 \cap\wt X_{j+1}^\mathrm{in}\right) \cap P_{i}$, with $P_{i}$ a waste bin,  into a union of waste bins.
\end{itemize}
Putting a, b and c  together, we obtain the desired property.

\paragraph{Step 6. Getting hypotheses $\mathbf{C_{3,4}}(k_{0})$ (convergence and thinness of the fibres)}
On the one hand, increasing the integer $n_{k_0+1}$ does not destroy hypotheses $\mathbf{C_{7,8}}(k_0)$ (nor, of course, the other hypotheses which do not involve the integer $k_{0}+1$). On the other hand, proposition~\ref{p.existence-H-1} shows that we can get hypotheses $\mathbf{C_{3,4}}(k_{0})$ by increasing if necessary the integer $n_{k_0+1}$. This completes the proof.
\end{proof}

\subsection{Proof of theorem~\ref{t.main}}
\label{ss.proof}

Let us assume that the hypotheses of theorem~\ref{t.main} hold. We modify  the construction of subsection~\ref{ss.bilan} in the following way. The sequences $(n_{k})$ and $(H_{k})$ are provided by proposition~\ref{p.existence-H-n} so that hypotheses $\mathbf{C_{1,\dots,8}}$ are satisfied (whereas in subsection~\ref{ss.bilan}, we had only obtained hypotheses $\mathbf{C_{1,\dots,6}}$). 
Proposition~\ref{p.extension} (and~\ref{p.extension-bis}) ensures the existence of the map $\Phi$ (which is one-to-one outside $\supp(\mu)$) and the homeomorphism $f$ such that $\Phi f= R \Phi$. And corollary~\ref{c.full-measure} provides an isomorphism between  $f$ and 
$$\left(\bigcup_{j\in\ZZ} R^j(K)\times C,h\right) \bigsqcup \left(\cM\setminus\bigcup_{j\in\ZZ} R^j(K),R\right).$$
as desired. Since hypothesis~ $\mathbf{C_4}$ is satisfied, proposition~\ref{p.recurrence} shows that, if $R$ is transitive (resp. minimal), then  $f$ is also transitive (resp. minimal). This completes the proof of theorem~\ref{t.main}.

\newpage
\appendix
\part*{Appendices}
\addcontentsline{toc}{part}{Appendices}
\renewcommand{\thesection}{\Roman{section}}

\section{Extension of homeomorphisms between Cantor sets}
\label{a.cantor}

The following proposition is needed in the paper to extend homeomorphisms between Cantor sets. The techniques involved in the proof are very classical (see~\cite{Bin,Hom}). 
Nevertheless the needed statement does not appear in the literature, firstly because 
the classical theorems are not written as extension theorems, and secondly because we deal with totally disconnected sets and not only Cantor sets. Thus we provide a proof.

\begin{prop}
\label{p.extension-cantor}
Let $X,X'$ be two copies of the unit cube $[0,1]^d$, and $\alpha$ a homeomorphism between $X$ and $X'$. Let $Q,Q'$ be two totally discontinuous tamely embedded compact sets in $\inte(X),\inte(X')$,
and $\beta$ be a homeomorphism between $Q$ and $Q'$.

Then there exists a homeomorphism $\gamma$ between $X$ and $X'$, which coincides with  $\beta$ on $Q$ and with $\alpha$ on the boundary of $X$.
\end{prop}

Since the totally discontinuous compact set $Q$  is tamely embedded in $X=[0,1]^d$ (see definition~\ref{d.tame-Cantor}), there exists  a strictly decreasing sequence  $(Q_{n})_{n \geq 0}$ of subsets of $\inte([0,1]^d)$ such that:
\begin{enumerate}
\item $\displaystyle Q = \bigcap_{n \geq 0}Q_{n}$; 
\item each $Q_{n}$ is the union of a finite family of pairwise disjoint closed topological balls (homeomorphic to $[0,1]^d$).
\end{enumerate}
 Note that the supremum of the diameters of the connected components of $Q_{n}$ goes to $0$ when $n$ tends to infinity (since $Q$ is totally disconnected). We also consider an analogous sequence of sets  $(Q'_{n})_{n \geq 0}$ for $Q'$. 

\begin{proof}
The proof consists in  successive reductions of the problem to an almost trivial case.

\paragraph{Reduction 1}
\emph{One can assume that $X=X'=[0,1]^d$ and $\alpha=\mathrm{Id}$.}
Indeed, one can bring back $Q'$ in $X$ \emph{via} $\alpha^{-1}$, solve the reduced problem in $X$, then compose the solution with $\alpha$.

\paragraph{Reduction 2}
We let $I_{1}$ be the arc $[0,1] \times \{(0,\dots , 0)\}$ in $X=X'=[0,1]^d$. 
\emph{One can further  assume that $Q,Q' \subset I_{1}$.}
To see this it suffices to find a homeomorphism $F$ of $[0,1]^d$ such that $F$ is the identity on the boundary and $F(I_{1})$ contains $Q$ (and do the same for $Q'$).
The homeomorphism $F$ is constructed by a technique of successive approximations (see for example \cite{Osb}, Theorem 1): $F=\lim F_{n}$ when $F_{n}(I_{1})$ meets every connected components of $\inte(Q_{n})$ and $F_{n+1} = G_{n} \circ F_{n}$ with $G_{n}$ supported in $Q_{n}$. The existence of $G_{n}$ comes from the fact that the group of homeomorphisms of $[0,1]^d$ which are the identity on the boundary acts transitively on the $n$-uplet of distinct points in $[0,1]^d$.

\paragraph{Reduction 3}
\emph{One can further assume that $X=X'=[0,1]^2$ (that is, $d=2$).}
Note that $Q,Q'$ are obviously tamely embedded in $[0,1]^2 \times  \{(0,\dots , 0)\}$
 (thanks to reduction 2).  So solving the problem in dimension 2  enables us first to construct $\gamma$ as a homeomorphism
of $[0,1]^2 \times  \{(0,\dots , 0)\}$ and then extend $\gamma$ to $[0,1]^d$ using an isotopy from $\gamma$ to the identity (Alexander trick).

\paragraph{Reduction 4} \emph{One can additionally assume that  for every $n \geq 0$, and every connected component $R$ of $Q_{n}$,  there exists a connected component $R'$ of $Q'_{n}$ such that $\beta(R \cap Q) = R' \cap Q'$.}
We start with the initial sequences $(Q_{n}), (Q'_{n})$ and we explain how to construct two decreasing sequences $(\wh Q_{n}), (\wh Q'_{n})$ that satisfy the additional assumption. 
The construction will satisfy $\wh Q_{2p}= Q_{m_{p}}$ and $\wh Q'_{2p+1}= Q'_{n_{p}}$ 
for some increasing sequences $(m_{p})$ and $(n_{p})$, which will guarantee that 
$\cap \wh Q_{n}=Q$ and $\cap \wh Q'_{n}=Q'$.

We let $\wh Q_{0}= Q_{0}$, and denote by $R_{1}, \dots R_{i_{0}}$ the connected components of $Q_{0}$.  
The family $\cF = \{\beta(R_{i} \cap Q), i=1 \dots i_{0}\}$ is a finite covering of $Q'$ by closed open sets.  Let $n_{0}$ be a sufficiently big integer so that the covering of $Q'$ by the connected components of $Q'_{n_{0}}$ is finer than $\cF$. A repeated use of the first item in   lemma~\ref{l.topologie} below (with $R = [0,1]^2$) provides a family  $\cF' = \{R'_{1} ,\dots R'_{i_{0}}\}$ of pairwise disjoint closed topological discs whose union contains $Q'_{n_{0}}$, and 
such that $R'_{i} \cap Q' = \beta(R_{i} \cap Q)$ for each $i$
 (each $R'_{i}$ contains all the connected  components $R'$ of $Q'_{n_{0}}$ satisfying $R' \cap Q' \subset \beta(R_{i} \cap Q)$). We let $\wh Q'_{0}$ be the union of the elements of $\cF'$.

Now let $\wh Q'_{1}= Q'_{n_{0}}$. We construct $\wh Q_{1}$ from $\wh Q'_{1}$ by
imitating the above construction of $\wh Q'_{0}$ from $\wh Q_{0}$ (note that this time  lemma~\ref{l.topologie} is applied with $R$ equal to the appropriate connected component of $\wh Q_{0}$, so that $\wh Q_{1} \subset \inte(\wh Q_{0})$).
We proceed to construct both sequences $(\wh Q_{n}), (\wh Q'_{n})$, exchanging at each step the role played by $Q$ and $Q'$.

\paragraph{Final proof}
We now prove the special case of the theorem corresponding to the successive reductions.
The wanted homeomorphism $\gamma$ is obtained as the limit of a sequence $(\gamma_{n})$ where $\gamma_{n}$ sends each component $R_{n}$ of $Q_{n}$ onto the corresponding component $R'_{n}$ of $Q'_{n}$ (see reduction 4). This sequence is constructed recursively with $\gamma_{n+1} = \delta_{n} \circ \gamma_{n}$ and the support of the homeomorphism $\delta_{n}$ included in $Q_{n}$.
The homeomorphism $\delta_{n}$ is given by the second item of lemma~\ref{l.topologie}.
\end{proof}

\begin{lemma}[in dimension 2]
\label{l.topologie}~
\begin{enumerate}
\item Let $R$ be a closed topological disc (\emph{i. e.} a set homeomorphic to $[0,1]^2$), and let 
$R_{1}, \dots R_{i_{0}}$ be a family of pairwise disjoint closed topological discs inside $\inte(R)$. Then there exists a closed topological disc $R_{1,2}$ inside $\inte(R)$ that contains $R_{1}$ and $R_{2}$ and is disjoint from all the $R_{i}$, $i \geq 3$.
\item Let $R'_{1}, \dots R'_{i_{0}}$ be another such family. Then there exists a homeomorphism of $R$, which is the identity on the boundary, and that maps each $R_{i}$ onto $R'_{i}$.
\end{enumerate}
\end{lemma}

The first item is proved the following way: choose one point $x_{i}$ in each disc $R_{i}$, choose a disc $\wh R_{1,2}$ that contains $x_{1}, x_{2}$ (in its interior) but not the $x_{i}, i \geq 3$. Then 
set $R_{1,2} = F^{-1}(\wh R_{1,2})$ where $F$ is a homeomorphism that is supported on a neighbourhood of the union of the $R_{i}$'s and contract each $R_{i}$ to a sufficiently small disc around $x_{i}$. The second item is a variation on Schoenflies theorem.
Here is one way to prove it. Firstly, the problem can easily be solved if the last word ``onto'' is replaced by the word ``into''. Secondly, by considering discs $R''_{i}$ slightly larger than the $R'_{i}$'s, the  problem is reduced to the case $i_{0}=1$. Then this is the classical Schoenflies theorem with compact support, see for example~\cite{Bin2}, Theorem II.6.C on page 31.

\normalsize

\begin{coro}
\label{c.cantor-ball}
Assume the hypotheses of proposition~\ref{p.extension-cantor}, together with the following additional data:  $(B_{j})_{j \geq0}, (B'_{i})_{i \geq0}$ are two sequences of pairwise disjoint topological closed balls respectively in $\inte(X), \inte(X')$ such that
\begin{enumerate}
\item the $B_{j}$'s are disjoint from $Q$,  the $B'_{i}$'s are disjoint from $Q'$;
\item $\displaystyle\mathop{\limsup}_{j\rightarrow\infty}B_j \subset  Q$;
\item $\displaystyle\mathop{\limsup}_{i\rightarrow\infty}B'_i = Q'$;
\item each $B_{j}$ is tamely embedded\footnote{See footnote~\ref{f.tamely-embedded} in subsection~\ref{ss.def-waste-bins}.}
 in $X$.
\end{enumerate}
Let $\varphi: \NN \ra \NN$ be any function.
Then there exists a homeomorphism $\gamma$ 
such that the conclusion of proposition~\ref{p.extension-cantor} holds, and 
for every $j \geq 0$ there exists $i > \varphi(j)$ such that $\gamma(B_{j}) \subset B'_{i}$.
\end{coro}

\begin{proof}
For every $j$, let $x_j$ be any point in the interior of $B_{j}$. 
We first extend continuously the map $\beta$ to the set $\{x_{j}\}$, with values in the $B'j$'s, in  the following way:
for every  $j\geq0$, pick a point $q$ in $Q$ which is among the closest to $x_{j}$, then pick some point $\beta(x_{j})$ in the interior of a ball $B'_{i}$ with $i > \varphi(j)$ and  such that $d(\beta(x_{j}), \beta(q)) \leq d(x_{j},q)$ (such a point exists thanks to item 3 of the hypotheses).  This can be done so that $\beta$ is still one-to-one, and since $\wh Q := Q \cup \{ x_{j} \}$ is a compact set (hypothesis 2), $\beta$ is still a homeomorphism on its image.

Now we apply proposition~\ref{p.extension-cantor} to the sets $\wh Q$ and $\wh Q' = \beta(\wh Q)$. Denote by $\wh \gamma$ the resulting homeomorphism.
Note that for every $j$, $x_{j}$ is a point of $\inte(B_{j})$ such that $\wh \gamma (x_{j}) \subset \inte(B'_{i_{j}})$ for some $i_{j} > \varphi(j)$.

Let $O_{j}$ be a topological closed ball containing $\wh \gamma(B_{j})$ in its interior. It is easy to find a  homeomorphism $\delta_{j}$  supported in $O_{j}$ such that $\delta_{j}(\wh \gamma(B_{j})) \subset \inte(B'_{i_{j}})$. Note that since the balls $B_{j}$ are tamely embedded, we can assume that the $O_{j}$ are pairwise disjoint. Let $\delta$ be the infinite (commutative) composition of all the $\delta _{j}$. Thanks to hypothesis 2, and since $Q$ is totally disconnected, the diameter of $B_{j}$ goes to $0$ when $j$ goes to infinity, thus we can assume that the same happens to the $O_{j}$'s, and then $\delta$ is a homeomorphism. It remains to  set $\gamma = \delta \circ \wh \gamma$.
\end{proof}

\section{Transitivity, minimality: proof of addendum~\ref{a.main}}
\label{a.recurrence}

In section~\ref{s.transitivity-minimality} and~\ref{s.general-H}, we have explained how to get some recurrence properties for the homeomorphism $f$ in the case where the dynamics of the initial homeomorphism $R$ is minimal (or transitive) on the whole manifold $\cM$. 
This appendix deals with the general case, when the dynamics of $R$ is not supposed to be minimal nor transitive on $\cM$. Assuming some recurrence properties for the homeomorphism $R$ on the subset
$$
\Lambda =\mbox{Cl}\left(\bigcup_{j\in \ZZ}R^j(K)\right),
$$
we explain how to obtain some recurrence properties for the homeomorphism $f$ on the set $\Phi^{-1}(\Lambda)$. In particular, we will prove addendum~\ref{a.main} (see subsection~\ref{ss.more-statements}). 

All along the appendix, we suppose we are given a sequence of collections of rectangles $(\cE_n^0)$ satisfying hypotheses $\mathbf{A_{1,3}}$, and such that the graph $\cG_0^0$ has no edge. 

\subsection{Hypotheses $\wt{\mathbf{B_4}}$ and $\wt{\mathbf{C_4}}$}
\label{ss.hyp-tilde}

Given a sequence of homeomorphisms  $(M_n)_{n\geq 1}$, we consider the following hypothesis.
\begin{itemize}
\item[]
\begin{itemize}
\item[$\widetilde{\mathbf{B_4}}$] \emph{(Fibres are thin: strong form)}\\
For any $\delta>0$, there exists $r\geq 1$ such that for any $n\geq r$,
the set $\Psi_n^{-1}(E^r_n\setminus E^0_n)$ is $\delta$-dense in the set
$\Psi_n^{-1}(E^r_n)$.
\end{itemize}
\end{itemize}

\begin{rema*}
Hypothesis $\mathbf{B_4}$ (see section~\ref{s.transitivity-minimality}) can be reformulated as follows: for every $\delta>0$ there exists $n_0$ such that, for any $n\geq n_0$, the set $\Psi_n^{-1}(\cM\setminus E_n^0)$ is $\delta$-dense in $\cM$. This shows that hypothesis $\wt{\mathbf{B_4}}$ is a stronger than $\mathbf{B_4}$.
\end{rema*}

If the sequence of homeomorphisms $(M_n)_{n\geq 1}$ satisfies hypotheses $\mathbf{B_{1,2,3}}$, then we can reformulate hypothesis $\wt{\mathbf{B_4}}$ using the map $\Psi=\lim\Psi_n$.
\begin{prop}
\label{p.reformulation-ter}
Hypothesis $\widetilde{\mathbf{B_4}}$ is satisfied if and only if the set $\Psi^{-1}(K)$ has empty interior in the set $\Psi^{-1}(\Lambda)$. 
\end{prop}

\begin{proof}
As in the proof of proposition~\ref{p.reformulation}, for any fixed $r>0$, we get that $\Psi^{-1}(\bigcup_{|i| \leq r} R^i(K))$
is the decreasing intersection of the compact sets $\Psi_n^{-1}(E_n^r)$ and that $\Psi^{-1}(\bigcup_{|i| \leq r} R^i(K)\setminus K)$ is the decreasing intersection of the sets $\Psi_n^{-1}(E^r_n\setminus E_n^0)$. This implies that $\wt{\mathbf{B_{4}}}$ is satisfied if and only if for every $\delta>0$, there exists $r>0$ such that the set $\Psi^{-1}(\bigcup_{|i| \leq r} R^i(K)\setminus K)$ is $\delta$-dense in $\Psi^{-1}(\bigcup_{|i| \leq r} R^i(K))$. Hence, $\wt{\mathbf{B_{4}}}$ is satisfied if and only if $\Psi^{-1}(\Lambda\setminus K)$ is dense in $\Psi^{-1}(\La)$, \emph{i.e.} if and only if $\Psi^{-1}(K)$ has empty interior in $\Psi^{-1}(\Lambda)$. 
\end{proof}  

\begin{coro}
\label{c.equiv-B4}
If $\Lambda=\cM$, hypotheses $\mathbf{B_4}$ and $\wt{\mathbf{B_4}}$ are equivalent. 
\end{coro}

\begin{proof}
This follows immediately from propositions~\ref{p.reformulation} and~\ref{p.reformulation-ter}.
\end{proof}

Replacing $\cM$ by $\Lambda$ in the proof of proposition~\ref{p.recurrence}, we obtain the following result.
\begin{prop}
\label{p.recurrence-ter}
If the dynamics of $R$ on $\Lambda$ is transitive (resp. minimal) and hypothesis $\widetilde{\mathbf{B_4}}$ is satisfied, then the dynamics of $g$ on $\Psi^{-1}(\Lambda)$ is also transitive (resp. minimal).
\end{prop}

Now suppose we are also given a sequence of integers $(n_k)_{k\in\NN}$ and a sequence of homeomorphisms $(H_k)_{k\in\NN}$. Then we can consider the following hypothesis.
\begin{itemize}
\item[]
\begin{itemize}
\item[$\widetilde{\mathbf{C_4}}$] \emph{(Fibres are thin: strong form)}\\
For any $\varepsilon>0$, there exists $r\geq 1$ such that for any $\ell\geq r$,
the set $\Phi_\ell^{-1}(\bE^r_\ell\setminus \bE^0_\ell)$ is $\varepsilon$-dense in the set
$\Phi_\ell^{-1}(\bE^r_\ell)$.
\end{itemize}
\end{itemize}
Of course, proposition~\ref{p.reformulation-ter}, corollary~\ref{c.equiv-B4} and proposition~\ref{p.recurrence-ter} are still valid if one replaces the maps $\Psi, g$ by the maps $\Phi, h$, and hypothesis  $\widetilde{\mathbf{B_4}}$ by hypothesis $\widetilde{\mathbf{C_4}}$.

\subsection{Realisation of hypotheses $\wt{\mathbf{B_4}}$ and $\wt{\mathbf{C_4}}$}
\label{ss.realis-tildeB4}

A important technical problem arises when one tries to realise simultaneously hypothesis $\wt{\mathbf{B_4}}$ (or $\wt{\mathbf{C_4}}$) and some other hypotheses:
\begin{itemize}
\item[--] the extraction process defined in subsection~\ref{ss.def-extraction} is crucial for obtaining some hypotheses (as e.g. $\mathbf{B_{3}}$);
\item[--] hypothesis $\wt{\mathbf{B_4}}$ is not preserved by the extraction process.
\end{itemize}
To overcome this problem, we will have to realise some hypotheses that are stronger than $\wt{\mathbf{B_{4}}}$, and preserved by the extraction process.

Given a sequence of homeomorphisms $(M_{n})_{n\geq 1}$ satisfying hypotheses~$\mathbf{B_{1,2}}$, we consider the sequence of homeomorphisms $(\Psi_{n}^0)_{n\geq 1}$ defined by:
$$
\Psi_{n}^0 = {M_{n}}_{\mid E_{n-1}^0} \circ {M_{n-1}}_{\mid E_{n-2}^0} \circ \cdots \circ {M_{1}}_{\mid E_{0}^0}
$$
where $M_{i|E_{i-1}^{0}}$ denotes the homeomorphism which is equal to $M_i$ on $E_{i-1}^{0}$ and equal to the identity elsewhere. The same proof as in section~\ref{ss.cons-B_123} shows that sequence of homeomorphisms $(\Psi_{n}^0)$ converges to a continuous map $\Psi^0$.

\begin{rema}
The map $\Psi_{n}^0$ can be thought of as ``what would be the extracted map $\bPsi_{1}$ if $n_{1}$ was equal to $n$.'' Also note that this definition is ``stable under the extraction process'': more precisely, if $(n_{k})$ is some sequence of integers as in section~\ref{s.extraction}, and $(\bM_{k})$ is the associated extracted sequence as in subsection~\ref{ss.def-extraction}, one can define a sequence 
$$
\bPsi_{k}^0 = {\bM_{k}}_{\mid \bE_{k-1}^0} \circ {\bM_{k-1}}_{\mid \bE_{k-2}^0} \circ \cdots \circ {\bM_{1}}_{\mid \bE_{0}^0}.
$$
Then one has $\bPsi_{k}^0 = \Psi_{n_{k}}^0$ for every $k$, so that the sequence 
$(\bPsi_{k}^0)$ converges towards the same map $\Psi^0$.
\end{rema}

We introduce a hypothesis that is a first step towards hypothesis $\wt{\mathbf{B_4}}$:
\begin{itemize}
\item[]
\begin{itemize}
\item[$\wt{\bm{\beta_4}}$]
The set $(\Psi^0)^{-1}(K)$ has empty interior in $(\Psi^0)^{-1}(\Lambda)$, where
$\displaystyle\Lambda = \adhe\left(\bigcup_{i \in \ZZ}R^i(K)\right)$.
\end{itemize}
\end{itemize}
This hypothesis should be compared to the definition of dynamical meagreness 
(definition~\ref{d.meagre}).
According to the previous remark, this hypothesis is stable under the extraction process.
We now prove that hypothesis~$\wt{\bm{\beta_4}}$ can be realised.

\begin{prop}
\label{p.existence-M-bis}
Assume that $K$ is dynamically meagre. Then there exists a sequence $(M_{n})_{n\in\NN}$ of homeomorphisms of $\cM$ such that hypotheses $\mathbf{B_{1,2,5,6}}$ and $\wt{\bm{\beta_4}}$ are satisfied.   
\end{prop}

\begin{proof}
We adapt the proof of proposition~\ref{p.existence-M}. Hypothesis $\wt{\bm{\beta_4}}$ is a straightforward consequence of the following points.

\begin{itemize}
\item[--] Let $x$ be a point of $K$, and for every $n$ let $\wh X_{n}$ be the rectangle of $\cE_{n}^0$ containing $x$. Since $K$ is dynamically meagre, for infinitely many values of $n$, the set $\Lambda \cap \inte(\wh X_{n-1} \setminus E_{n}^0)$ is not empty. Also note that the dynamical meagreness implies that $\Lambda$ has no isolated point, so that
$\Lambda \cap \inte(\wh X_{n-1} \setminus E_{n}^0)$ is infinite as soon as it is non-empty.
\item[--] In the inductive construction of the sequence $(M_{n})_{n\in\NN}$, we may require the following additional property to be satisfied:
\begin{description}
\item For every rectangle $\wh X \in \cE_{n-1}^0$ such that the set $\Lambda \cap \inte(\wh X \setminus E_{n}^0)$ is not empty, the set $(\Psi_{n}^0)^{-1}(\Lambda \cap \inte(\wh X \setminus E_{n}^0))$ is $\frac{1}{n}$ dense in $(\Psi_{n}^0)^{-1}(\wh X)$.
\end{description}
This property is easily obtained by the following modification in the step 2 of the proof of proposition~\ref{p.existence-M}:
we choose the set $A'$ included in the (infinite) set $\Lambda \cap \inte(\wh X \setminus E_{n}^0)$, and the constant $\delta_{n}$ smaller than $\frac{1}{n}$.
\item[--] Hypothesis $\mathbf{B_{1}}$ implies that $(\Psi^0)^{-1}(\Lambda \cap \inte(\wh X \setminus E_{n}^0)) = (\Psi_{n}^0)^{-1}(\Lambda \cap \inte(\wh X \setminus E_{n}^0))$.
\end{itemize}
\end{proof}

The next step is to show that hypothesis~$\wt{\mathbf{B_4}}$ can be obtained from hypothesis~$\wt{\bm{\beta_4}}$ by extracting. 
For this we need a quantitative version of hypothesis~$\wt{\mathbf{B_4}}$.
\begin{itemize}
\item[]
\begin{itemize}
\item[$\widetilde{\mathbf{B_4}}(\delta)$]
There exists $r\geq 1$ such that for any $n\geq r$,
the set $\Psi_n^{-1}(E^r_n\setminus E^0_n)$ is $\delta$-dense in the set
$\Psi_n^{-1}(E^r_n)$.
\end{itemize}
\end{itemize}

\begin{prop}
\label{p.B4-tilde}
Assume that the sequence $(M_{n})$ satisfies hypotheses~$\mathbf{B_{1,2}}$ and $\wt{\bm{\beta_4}}$. Let $(n_{1}, \dots, n_{k-1})$ be an increasing sequence. Fix some positive $\delta$. Then for any increasing sequence $(n_{\ell})_{\ell \geq k}$ with $n_{k}$ large enough, 
the hypothesis $\widetilde{\mathbf{B_4}}(\delta)$ is satisfied by the extracted sequences
$(\bcE^0_{\ell})$ and $(\bPsi_{\ell})$.
\end{prop}

During the proof of proposition~\ref{p.B4-tilde}, we will need to compare the maps $\bPsi_{k}^{-1}$ and $(\Psi^0_{n_{k}})^{-1}$. This is the purpose of the following lemma:

\begin{lemma}
\label{l.yenamar}
For every increasing sequence $(n_{1}, \dots, n_{k})$, the maps $\bPsi_{k}^{-1}$ and $(\Psi^0_{n_{k}})^{-1}$ coincide on the set $E_{k-1}^0$.
\end{lemma}

The proof is an easy induction, and makes use of the equality $E_{i-1}^{k-1} \cap E_{k-1}^0 = E_{i-1}^0$ for every $i=n_{k-1}+1, \dots , n_{k}$  that follows from hypothesis $\mathbf{A_{1.c}}$
(see subsection~\ref{ss.consequences-A}).

\begin{proof}[Proof of the proposition]
We will prove the following fact.
 There exists $r \geq 1$ such that for any $n_{k}$ large enough,
\begin{enumerate}
\item  the set $\bPsi_k^{-1}((E^r_{n_{k}}\setminus E^0_{n_{k}}) \cap E_{k}^0)$ is $(2\delta)$-dense\footnote{That is, every point $\bPsi_k^{-1}(E^0_{n_{k}})$ is at distance less than $2\delta$ from a point of the set $\bPsi_k^{-1}((E^r_{n_{k}}\setminus E^0_{n_{k}}) \cap E_{k}^0)$.} in the set $\bPsi_k^{-1}(E^0_{n_{k}})$;
\item the diameter of each connected component of $\bPsi_k^{-1}((E^r_{n_{k}}\setminus E^0_{n_{k}}) \cap E_{k}^0)$ is less than $\delta$.
\end{enumerate}
For this, we first choose an $\varepsilon > 0$ which is less than $\delta $ and such that the $\varepsilon$-neighbourhood of $(\Psi^0)^{-1}(K)$ is included in $(\Psi^0)^{-1}(E_{k}^0)$.
We now  apply $\widetilde{\bm{\beta_4}}$. Since the map $(\Psi^0)^{-1}$ is a homeomorphism outside $K$, this provides us with a positive integer $r$ (we may assume $r \geq k$) such that 
$(\Psi^0)^{-1} (K_{r}) $ is $\varepsilon$-dense in $(\Psi^0)^{-1}(K)$, where 
$K_{r} = (\cup_{\mid i \mid \leq r} R^i (K)) \setminus K$.
By the choice of $\varepsilon$ this implies
that
\begin{itemize}
\item[]
\begin{itemize}
\item[$(\star)$]   the set $(\Psi^0)^{-1} (K_{r} \cap E_{k}^0) $ is $\delta$-dense in $(\Psi^0)^{-1}(K)$. 
\end{itemize}
\end{itemize}
Using hypothesis $\mathbf{A_{3}}$, we choose an integer $N \geq r$ large enough so that 
\begin{itemize}
\item[]
\begin{itemize}
\item[$(\star\star)$]  the diameter of each connected component of  the set $(\Psi^0)^{-1}((E^r_{N}\setminus E^0_{N})\cap E_{k}^0)$ is less than $\delta$;
\item[$(\star\star\star)$] the set $(\Psi^0)^{-1}(K)$ is $\delta $-dense in the set $(\Psi^0)^{-1}(E^0_{N})$.
\end{itemize}
\end{itemize}
Note that for every $s \geq r$ we have $(\Psi_{s}^0)^{-1} = (\Psi^0)^{-1} $ on $\cM\setminus E_{r}^0$, in particular on $E^r_{N}\setminus E^0_{N}$ (by compatibility, see  subsection~\ref{ss.consequences-A}), so that $(\Psi^0)^{-1}((E^r_{N}\setminus E^0_{N}) \cap E_{k}^0) = (\Psi_{s}^0)^{-1}((E^r_{N}\setminus E^0_{N}) \cap E_{k}^0)$.

Now let us check that, for any $n_{k} \geq N$, we have properties 1 and 2 stated at the beginning of the proof. Using  lemma~\ref{l.yenamar} and the above equality (for $s=n_{k}$), we get
$$
\begin{array}{rcl}
\bPsi_k^{-1}((E^r_{n_{k}}\setminus E^0_{n_{k}}) \cap E_{k}^0)  & =   & 
(\Psi^0_{n_{k}})^{-1}((E^r_{n_{k}}\setminus E^0_{n_{k}}) \cap E_{k}^0) \\
 & =   & (\Psi^0)^{-1}((E^r_{n_{k}}\setminus E^0_{n_{k}}) \cap E_{k}^0) \\
 & \subset   & (\Psi^0)^{-1}((E^r_{N}\setminus E^0_{N}) \cap E_{k}^0).
\end{array}
$$
This, together with property $(\star\star)$, implies property~2.
For property~1, let us note that the set $K_{r}$ is contained in $E^r_{n_{k}}\setminus E^0_{n_{k}}$. Thus the set $\bPsi_k^{-1}((E^r_{n_{k}}\setminus E^0_{n_{k}}) \cap E_{k}^0)$ contains a set that is $\delta$-dense in $(\Psi^0)^{-1}(K)$ (by property $(\star)$), which in turn is $\delta $-dense in $(\Psi^0)^{-1}(E^0_{N})$ (by property $(\star\star\star)$), which contains $(\Psi^0)^{-1}(E^0_{n_{k}})$: as a consequence,
 $\bPsi_k^{-1}((E^r_{n_{k}}\setminus E^0_{n_{k}}) \cap E_{k}^0)$ is $(2\delta)$-dense in 
$(\Psi^0)^{-1}(E^0_{n_{k}})$. It remains to note that
$$
\begin{array}{rcl}
(\Psi^0)^{-1}(E^0_{n_{k}})  & = & (\Psi_{n_{k}}^0)^{-1}(E^0_{n_{k}})  \\
 & = & (\bPsi_{{k}})^{-1}(E^0_{n_{k}}).
\end{array}
$$
The first equality follows from the definition of $\Psi^0$, and the second one from lemma~\ref{l.yenamar}.

To complete the proof of the proposition,  let  $(n_{\ell})_{\ell \geq k}$ be an increasing sequence with $n_{k} \geq N$ as above, and let $\ell \geq r$ (thus $\ell \geq k$). Then
\begin{itemize}
\item[--] by hypothesis $\mathbf{A_{1.b}}$, the set $\bPsi_\ell^{-1}((E^r_{n_{\ell}}\setminus E^0_{n_{\ell}}) \cap E_{k}^0)$ meets every connected component of the set $\bPsi_\ell^{-1}(((E^r_{n_{k}}\setminus E^0_{n_{k}}) \cap E_{k}^0)$;
\item[--]  by definition of the map $\bPsi_{\ell}$, one has $\bPsi_\ell^{-1}((E^r_{n_{k}}\setminus E^0_{n_{k}}) \cap E_{k}^0) = \bPsi_k^{-1}((E^r_{n_{k}}\setminus E^0_{n_{k}}) \cap E_{k}^0)$.
\end{itemize}
Thus properties 1 and 2 above imply that the set $\bPsi_\ell^{-1}((E^r_{n_{\ell}}\setminus E^0_{n_{\ell}}) \cap E_{k}^0)$  is $(3\delta)$-dense in the set $\bPsi_k^{-1}(E^0_{n_{k}})$.
Finally, observe that 
$$\bPsi_k^{-1}(E^0_{n_{k}}) = \bPsi_\ell^{-1}(E^0_{n_{k}})\supset\bPsi_\ell^{-1}(E^0_{n_{\ell}}).$$ Hence, the set $\bPsi_\ell^{-1}((E^r_{n_{\ell}}\setminus E^0_{n_{\ell}}) \cap E_{k}^0)$  is $(3\delta)$-dense in the set $\bPsi_\ell^{-1}(E^0_{n_{\ell}})$. Hence hypothesis $\widetilde{\mathbf{B_4}}(3\delta)$ is satisfied by the extracted sequences.
\end{proof}

Finally we will get hypothesis $\widetilde{\mathbf{C_4}}$ in the same way as for $\widetilde{\mathbf{B_4}}$.
When given two sequences $(n_{k})$ and $(H_{k})$, we consider: 
\begin{itemize}
\item[]
\begin{itemize}
\item[$\widetilde{\mathbf{C_4}}(\varepsilon)$]
There exists $r\geq 1$ such that for any $\ell\geq r$,
the set $\Phi_\ell^{-1}(\bE^r_\ell\setminus \bE^0_\ell)$ is $\varepsilon$-dense in the set
$\Phi_\ell^{-1}(\bE^r_\ell)$.
\end{itemize}
\end{itemize}

\begin{prop}
\label{p.existence-H-ter}
Assume that the sequence $(M_{n})$ satisfies hypotheses~$\mathbf{B_{1,2}}$ and $\wt{\bm{\beta_4}}$. Let $(n_{1}, \dots, n_{k-1})$  and $(H_{1}, \dots, H_{k-1})$ be two finite   sequences. Fix some positive $\varepsilon$. Then for any increasing sequence $(n_{\ell})_{\ell \geq k}$ with $n_{k}$ large enough and any sequence $(H_{\ell})_{\ell \geq k}$ such that hypotheses~$\mathbf{C_{1,2}}$ are satisfied, the hypothesis $\widetilde{\mathbf{C_4}}(\varepsilon)$ is also satisfied.
\end{prop}

\begin{proof}
We consider the homeomorphism $Z = \Phi_{k-1}^{-1} \circ \bPsi_{k-1}$.
Choose some $\delta  >0$ such that $d(x,y) \leq \delta $ implies $d(Z(x),Z(y)) \leq \varepsilon$.
According to the proof of proposition~\ref{p.B4-tilde}, we can find integers $r \geq k$ and $N \geq r$ such that for any choice of  $n_{k} \geq N$, 
\begin{enumerate}
\item  the set $\bPsi_k^{-1}((E^r_{n_{k}}\setminus E^0_{n_{k}}) \cap E_{k}^0)$ is $\delta$-dense in the set
$\bPsi_k^{-1}(E^0_{n_{k}})$;
\item the diameter of each connected component of $\bPsi_k^{-1}((E^r_{n_{k}}\setminus E^0_{n_{k}}) \cap E_{k}^0)$ is less than $\delta$.
\end{enumerate}
Consider some  sequences $(n_{\ell})_{\ell \geq k}$ and  $(H_{\ell})$ as in the statement of proposition~\ref{p.existence-H-ter}. Assume
 $n_{k} \geq N$.
By definition we have $\Phi_{k}^{-1} =  \Phi_{k-1}^{-1} \circ \bM_{k}^{-1} \circ H_{k}^{-1} = Z \circ \bPsi_{k}^{-1} \circ H_{k}^{-1}$.
Since the support of $H_{k}$ is $E_{n_{k}}^{k-1}$, the sets $(E^r_{n_{k}}\setminus E^0_{n_{k}}) \cap E_{k}^0$ and $E^0_{n_{k}}$ are preserved by $H_{k}^{-1}$. Consequently, points 1 and 2 above are still valid when $\bPsi_{k}$ is replaced by $\Phi_{k}$ and $\delta $ by $\varepsilon$.

It remains to prove that this implies that: for any $\ell \geq r$, the set $\Phi_\ell^{-1}(E^r_{n_\ell}\setminus E^0_{n_\ell})$ is $(2\varepsilon)$-dense in the set
$\Phi_\ell^{-1}(E^r_{n_\ell})$. This is entirely analogous to the end of the proof of proposition~\ref{p.B4-tilde}.
\end{proof}

\subsection{Proof of addendum~\ref{a.main}}

We modify the constructions of sections~\ref{ss.bilan} and~\ref{ss.proof} in the following way. Using proposition~\ref{p.existence-M-bis}, we can construct a sequence $(M_n)$ such that hypotheses $\mathbf{B_{1,2,4,5,6}}$ and $\wt{\bm{\beta_4}}$ are satisfied. Then, using proposition~\ref{p.existence-H-ter} together with proposition~\ref{p.existence-H-n}, we can construct a sequence of integers $(n_k)$ and a sequence of homeomorphisms $(H_k)$ such that hypotheses $\mathbf{C_{1,\dots,8}}$ and $\widetilde{\mathbf{C_4}}$ are satisfied\footnote{More precisely, we have to modify the proof of proposition~\ref{p.existence-H-n} so that the choice of the integer $n_{k_{0+1}}$ at step 6 involves not only proposition~\ref{p.existence-H-1} but also proposition~\ref{p.existence-H-ter}.}.
 Since the measure $\mu$ is ergodic, the set  $\Lambda$ is equal to the support of $\mu$. Hence, $R$ is transitive on $\Lambda$. Hence, proposition~\ref{p.recurrence-ter} (applied to $\Phi$, $f$ and $\widetilde{\mathbf{C_4}}$ instead of $\Psi$, $g$ and $\widetilde{\mathbf{B_4}}$) shows that:
\begin{itemize}
\item[--] in any case, $f$ is transitive on $\Phi^{-1}(\supp(\mu))$;
\item[--] if $R$ is minimal on $\Phi^{-1}(\supp(\mu))$, then $R$ is also minimal on $\supp(\mu)$.
\end{itemize}
This completes the proof of addendum~\ref{a.main}.

\section{Some examples}
\label{s.Denjoy}

In this appendix, we want to illustrate some of the results of the paper (mainly proposition~\ref{p.extension}, proposition~\ref{p.recurrence-ter} and our main theorem~\ref{t.main}) by a few examples.

\subsection{Denjoy counter-examples}

The simplest setting for proposition~\ref{p.extension} is when the collection $\cE_n^0$ is made of a single rectangle for every $n$ (or equivalently, when $K$ is a single point). This yields various generalisations of the classical Denjoy counter-examples on $\SS^1$.

\begin{prop}
\label{p.Denjoy}
Let $R$ be a homeomorphism on  a compact manifold $\cM$, and $x$ a point of $\cM$ which is not periodic under $R$. Consider a compact subset $D$ of $\cM$ which can be written as the intersection of a strictly decreasing sequence $(\wt X_{n})_{n \geq 0}$ of tamely embedded topological closed balls.
Then there exist a homeomorphism $f:\cM\to\cM$ and a continuous onto map $\Phi:\cM\to\cM$ such that $\Phi \circ  f = R \circ \Phi$, and such that 
\begin{itemize}
\item[--] $\Phi^{-1}(x) = D$;
\item[--] $\Phi^{-1}(y)$ is a single point if $y$ does not belong to the $R$-orbit of $x$. 
\end{itemize}
\end{prop}

\begin{remas*}~
\begin{itemize}
\item[--] The properties of $\Phi$ and $f$ imply that, if $R$ is minimal, then the set $\cM \setminus \bigcup_{n \in  \ZZ}f^n(\inte(D))$ is the only minimal closed  invariant set for $f$.
\item[--] So, if $R$ is an irrational rotation on  $\cM=\SS^1$ and $D$ is a non-trivial interval of $\SS^1$, then $f$ is a classical Denjoy counter-example. 
\item[--] In any case, if the interior of $D$ is non empty, then it is an open wandering set for $f$. In particular, if $R$ is minimal and $D$ has non-empty interior, then the dynamics of $f$ is very similar to the dynamics of Denjoy counter-examples on the circle.
\item[--] If $R$ is minimal and $D$ has empty interior, then $f$ is minimal. In this case, we obtain a kind of ``Denjoy-counter-example" whose dynamical behaviour is actually quite different from those of the classical Denjoy counter-examples on the circle. 
\end{itemize}
\end{remas*}

\begin{proof}[Proof of proposition~\ref{p.Denjoy}]
First note that we may assume that the point $x$ belongs to the interior of $\wt X_{0}$.\footnote{Choose some homeomorphism $\Psi$ of $\cM$ such that $x \in \inte(\Psi(\wt X_{0}))$; solve the problem with $(\wt X_{n})$ and $D$ replaced by $(\Psi(\wt X_{n}))$ and $\Psi(D)$; then replace $f$ by $\Psi^{-1}\circ f \circ \Psi$ and $\Phi$ by $\Phi \circ \Psi$.}
We choose  a decreasing sequence of rectangles $(X_n)_{n\in\NN}$ such that $X_0=\wt X_0$ and $\bigcap_{n\in\NN} X_n=\{x\}$ and such that, for every $n\in\NN$,  the rectangles $f^{-(n+1)}(X_n),\dots,X_n,\dots,f^{n+1}(X_n)$ are pairwise disjoint. Observe that to get the last property, it suffices to choose the rectangle $X_n$ small enough. We set $\cE_n^0:=\{X_n\}$. It is very easy to check that the sequence of collection rectangles $(\cE_n^0)_{n\in\NN}$ satisfies hypotheses $\mathbf{A_{1,2,3}}$. 

Secondly, we construct  a sequence of homeomorphisms $(M_n)_{n\geq 1}$ satisfying hypotheses $\mathbf{B_{1,2}}$ and such that $\Phi_n^{-1}(X_n)=\wt X_n$ for every $n$. We proceed as follows. Assume that $M_1,\dots,M_{n-1}$ have been constructed. Then $\Phi_{n-1}(\wt X_{n})$ is a strict sub-rectangle of the rectangle  $\Phi_{n-1}(\wt X_{n-1})=X_{n-1}$. So we can construct $M_n$ on the rectangle $X_{n-1}$ such that $M_n$ is the identity on the boundary of $X_{n-1}$ and such that $M_n(\Phi_{n-1}(\wt X_{n}))=X_{n}$  (i.e. $\Phi_n(\wt X_n)=X_n$). Then hypothesis~$\mathbf{B_1}$ and $\mathbf{B_2}$ do not leave any freedom for  the construction of $M_n$ on $\cM\setminus X_{n-1}$ (note that hypothesis $\mathbf{B_2}$ does not cause any problem since the rectangles $f^{-(n+1)}(X_n),\dots,X_n,f^{n+1}(X_n)$ are pairwise disjoint). 

Thirdly, using proposition~\ref{p.get-B_3}, we can ``extract a sub-sequence" in order to obtain hypothesis $\mathbf{B_3}$. Then, we can apply proposition~\ref{p.extension} in order to get a map $\Phi$ and a homeomorphism $f$ satisfying all the desired properties (the equality $\Phi^{-1}(x)=D$ follows from the equalities $\Phi_n^{-1}(X_n)=\wt X_n$ and $\Phi^{-1}(x)=\bigcap_{n\in\NN}\Phi_n^{-1}(x)$, see lemma~\ref{l.rectangles-pieges}).
\end{proof}

\subsection{Different ways of blowing-up an invariant  circle}

Now, we would like to illustrate hypothesis $\wt{\mathbf{B_{4}}}$ on some simple examples. For this purpose, we consider an irrational rigid rotation $R$ of the sphere $\SS^2$ (fixing the two poles $N$ and $S$). We denote by $\Lambda$ the equatorial circle of $\SS^2$ (which is invariant under $R$), and we pick a point $x\in\Lambda$. Using proposition~\ref{p.Denjoy}, we can construct a homeomorphism $f$ and a map $\Phi$ such that  $\Phi \circ  f = R \circ \Phi$, such that $\Phi^{-1}(x)$ is a non-trivial ``vertical'' segment $J$ and  such that $\Phi^{-1}(y)$ is a single point if $y$ does not belong to the $R$-orbit of $x$. It follows that $\wt\Lambda=\Phi^{-1}(\Lambda)$ is a one-dimensional (connected with empty interior) $f$-invariant compact set which separates $\SS^2$ into two connected open sets. 

\begin{figure}[htbp]
\includegraphics[width=15cm,height=7cm]{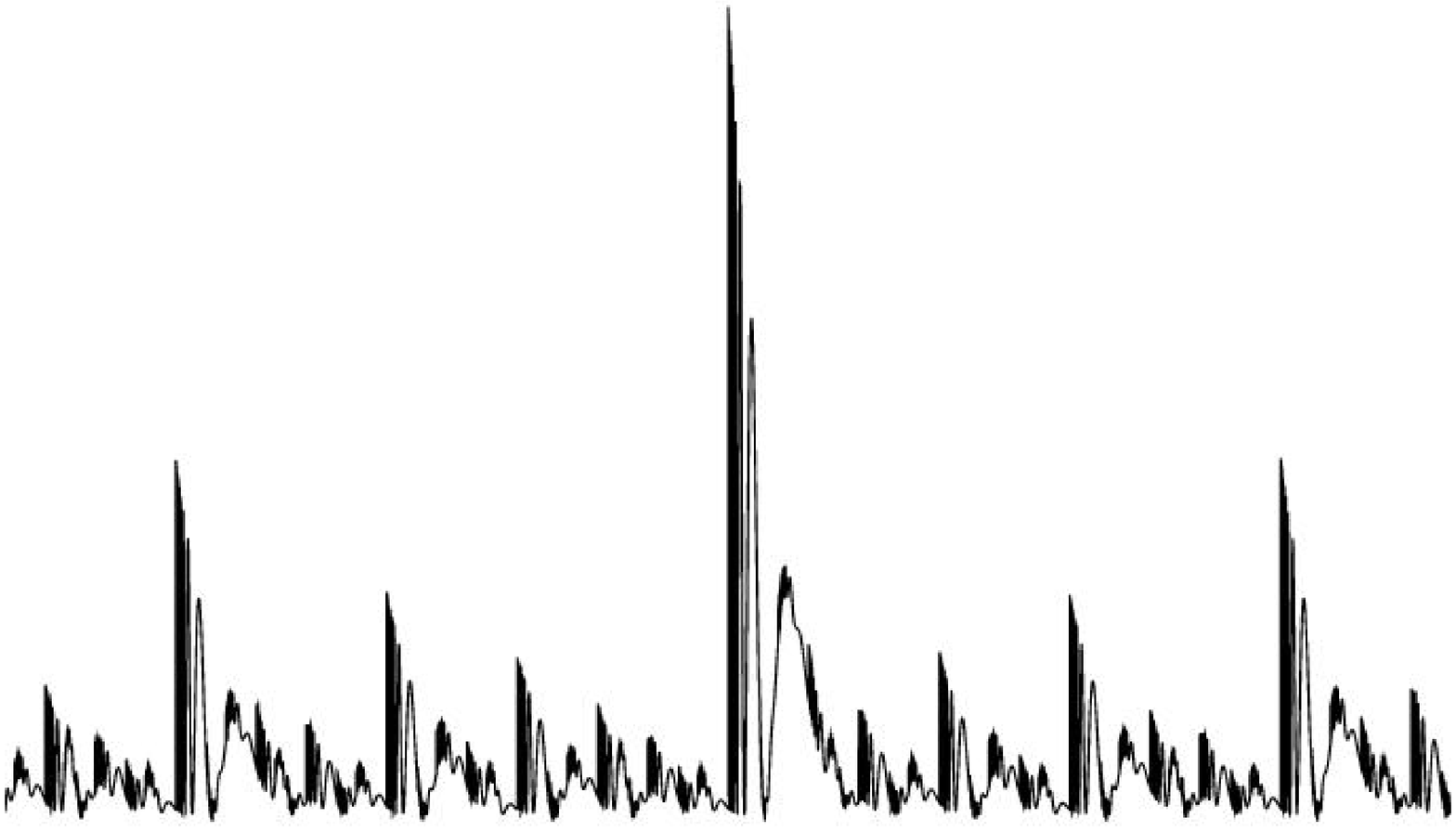}
\includegraphics[width=15cm,height=7cm]{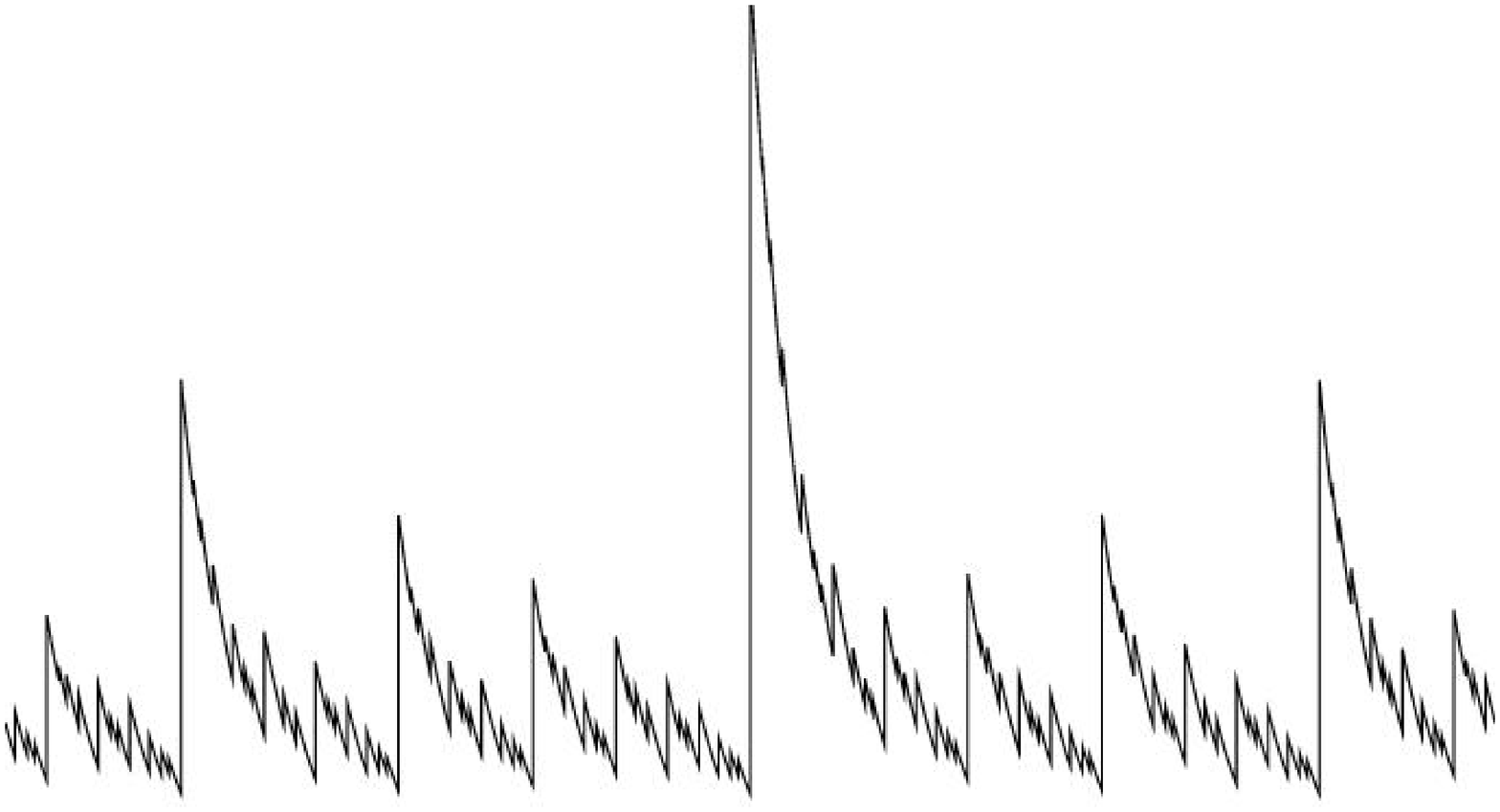}
\includegraphics[width=15cm,height=7cm]{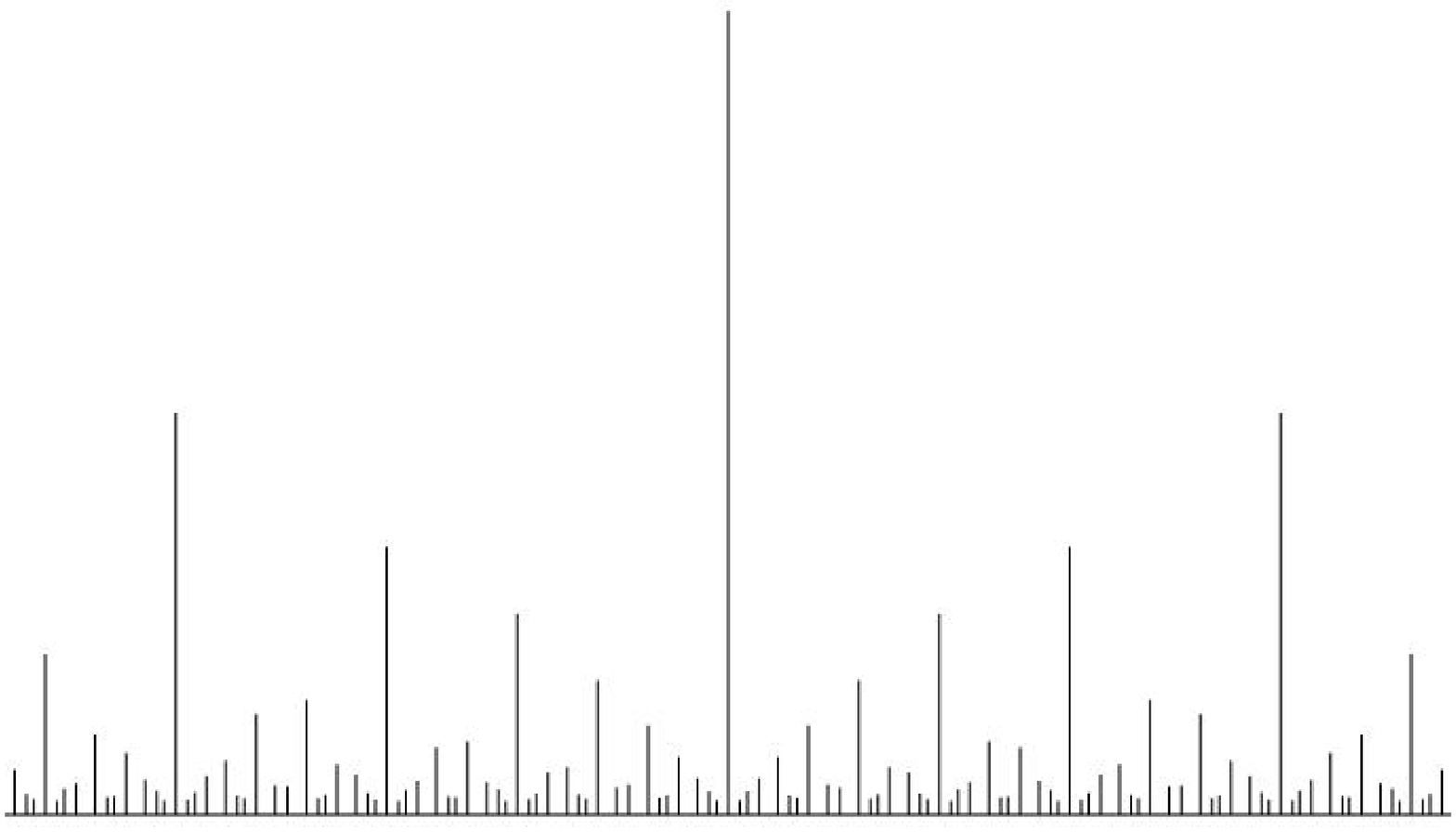}
\caption{Different ways of blowing-up an invariant  circle}
\label{f.denjoy-examples}
\end{figure}

Moreover, according to the way we choose the $M_{n}$'s, we can get quite different topologies for the set $\wt\Lambda$ and quite different dynamics for the restriction of $f$ to $\wt\Lambda$. Here are three possible types of behaviours: 
\begin{itemize}
\item[--] $\wt \Lambda$ is a non-arcwise connected set which is minimal for $f$ (figure~\ref{f.denjoy-examples}, I); 
\item[--] $\wt \Lambda$ is a topological circle which is not minimal for $f$: the restriction of $f$ to $\wt\Lambda$ is a Denjoy counter-example on the circle, the vertical segment $J$ is wandering (figure~\ref{f.denjoy-examples}, II);
\item[--] $\wt \Lambda$ contains a circle which is a minimal set for $f$, but is not equal to this circle (figure~\ref{f.denjoy-examples}, III, where the minimal set is the equatorial circle).
\end{itemize}
We will not explain precisely how to obtain such examples. We just note that an important point in the construction is the fact that hypothesis $\wt{\mathbf{B_{4}}}$ is satisfied in the first example, but not in the other two.

\begin{rema*}
The construction of the above examples can be made in such a way that $\Phi$ is $C^\infty$ on $\SS^2\setminus\wt\Lambda$. Moreover, if we identify $\SS^2\setminus\{N,S\} $ to the annulus $\SS^1\times\RR$ and see $f$ as a homeomorphism of $\SS^1\times\RR$, then all the constructions can be made in such a way that $f$ is a fibered homeomorphism (i.e. is of the form $f(x,y)=(x+\alpha,f_x(y))$). In this context, it is interesting to compare the three above examples with the classification of fibered homeomorphisms of T. Jaeger and J. Stark (see~\cite{StaJae}). This will be one of the purposes of a forthcoming paper~\cite{BegCroJaeLeR}.
\end{rema*}

\subsection{Pseudo-rotations with positive topological entropy on the $2$-sphere}
\label{ss.pseudo-rotation}

To end up this appendix, we would like to apply theorem~\ref{t.main} to obtain a more sophisticated example. An \emph{irrational pseudo-rotation} of the sphere $\SS^2$ is a  homeomorphism which preserves the orientation and the Lebesgue measure, has two fixed points $N,S$ and no other periodic point. The rotation set of an irrational  pseudo-rotation is reduced to a single irrational number (the \emph{angle} of the pseudo-rotation, see~\cite{LeC,BegCroLeR0}).

We choose any irrational angle $\alpha \in \SS^1$ and denote by $R_{\alpha}$ the rigid rotation of angle $\alpha$. We denote by $\Lambda$ be the equatorial circle invariant by $R_{\alpha}$.

\begin{prop}
\label{p.pseudo-rotation}
For every $\alpha\in\RR\setminus\QQ$, there exists an irrational pseudo-rotation $f$ on $\SS^2$ of angle $\alpha$  with positive topological entropy. 

Furthermore, there exists a continuous onto semi-conjugacy  $\Phi$ between  $f$ and the rigid rotation $R_{\alpha}$. If $\Lambda$ is the equatorial circle of $\SS^2$ (invariant under $R_\alpha$), then the set $\wt \Lambda=\Phi^{-1}(\Lambda)$ is a one-dimensional (connected with empty interior) minimal closed $f$-invariant set which carries all the entropy of $f$. It separates the sphere into two connected open sets. The map $\Phi$ is smooth on $\SS^2\setminus\wt\Lambda$; thus the restriction $f$ to ${ \SS^2 \setminus \wt \Lambda}$ is  $C^\infty$-conjugate to the restriction of $R_{\alpha}$ to  ${\SS^2 \setminus \Lambda}$.
\end{prop}

\begin{proof}
The proposition is almost a corollary of theorem~\ref{t.main} applied in the case where the manifold $\cM$ is the sphere $\SS^2$, the homeomorphism $R$ is the rigid rotation $R_{\alpha}$, the measure $\mu$ is the unique $R$-invariant measure supported by the equatorial circle $\Lambda$, the set $A$ is the equatorial circle $\Lambda$ and the map $h$ is the product of $R_{\alpha\mid\Lambda}$ by a Cantor homeomorphism with positive topological entropy. The only point which does not follow from theorem~\ref{t.main} is the fact that $\Phi$ is $C^\infty$ on $\SS^2\setminus \wt \Lambda$. This is an immediate consequence of the two following remarks:
\begin{itemize}
\item[--] in the construction of the homeomorphisms $M_n$'s and $H_k$'s explained in subsections~\ref{ss.realis-B_12456},~\ref{ss.realis-C_1256} and~\ref{p.existence-H-n}, it is not difficult to ensure that all these homeomorphisms are $C^\infty$ outside the iterates of the Cantor set $K\times C$;
\item[--] for every point $\wt x\in\SS^2\setminus\wt \Lambda$, all but a finite number of the homeomorphisms $M_n$'s and $H_k$'s are equal to the identity in a neighbourhood of $\wt x$; hence  $\Phi$ is locally the composition of a finite number of $C^\infty$ maps.
\end{itemize}
\end{proof}


\end{document}